\documentclass[12pt]{amsart}
\usepackage{amssymb}
\usepackage{amsmath}
\usepackage{amsbsy}
\usepackage{amscd}
\usepackage{url}
\usepackage[dvips]{graphicx}
\usepackage{here}
\advance\textheight by \topskip
\makeatletter
\def\Bbb{\mathbb}

\newenvironment{pf*}[1]{\proof[#1]}{\endproof}
\renewcommand{\thesubsection}{\thesection(\@roman\c@subsection)}
\makeatother
\newtheorem{Theorem}[equation]{Theorem}
\newtheorem{Corollary}[equation]{Corollary}
\newtheorem{Lemma}[equation]{Lemma}
\newtheorem{Proposition}[equation]{Proposition}

\theoremstyle{definition}

\makeatletter
\renewcommand\section{\@startsection{section}{1}%
  {\z@}{.7\linespacing\@plus\linespacing}{.5\linespacing}%
  {\reset@font\normalfont\bfseries\centering}}
\makeatother

\theoremstyle{remark}

\newtheorem*{Acknowledgements}{Acknowledgements}
\numberwithin{equation}{section}
\numberwithin{figure}{section}

\newcommand{\thmref}[1]{Theorem~\ref{#1}}

\newcommand{\Romnum}[1]{\expandafter\uppercase\expandafter{\romannumeral #1}} 
\newcommand{\C}{{\Bbb C}}

\newcommand{\CP}{\operatorname{\C P}}

\begin{document}
\title[On homology 3-spheres defined by two knots]{On homology 3-spheres defined by two knots}
\author{Masatsuna Tsuchiya}
\address{Department of mathematics, Gakushuin University, 5-1, Mejiro 1-chome, Toshima-ku, Tokyo, 171-8588, Japan}
\email{tsuchiya@math.gakushuin.ac.jp}
%
\begin{abstract}
We show that if each of $K_1$ and $K_2$ is a trefoil knot or a figure eight knot, the homology $3$-sphere defined by the Kirby diagram which is a simple link of $K_1$ and $K_2$ with framing ($0$, $n$) is represented by an $n$-twisted Whitehead double of $K_2$ .
\end{abstract}
\maketitle
%
%
\section{Introduction}\label{sec:intro}
%
%
We define $W_n(K_1, K_2)$ to be the $4$-dimensional handlebody represented by the following Kirby diagram, and define $M_n(K_1, K_2)$ to be $\partial(W_n)$, where $K_1$ and $K_2$ are knots. Note that $M_n(K_1, K_2)$ is a homology $3$-sphere.
\begin{figure}[h]
 \begin{center}
  \includegraphics[height=25mm]{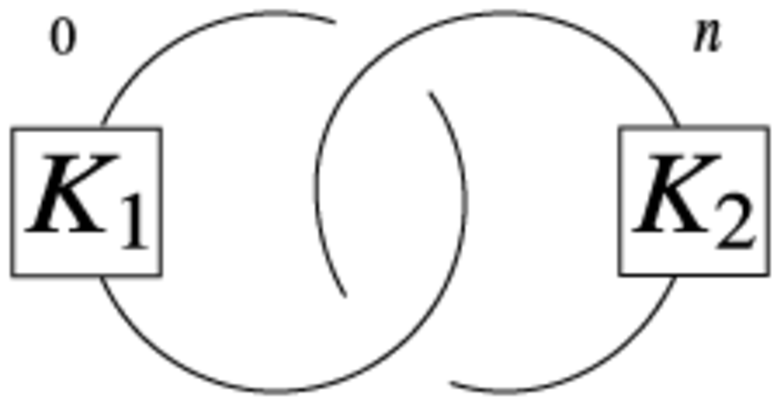}
 \end{center}
 \caption{$W_n(K_1, K_2)$}
 \label{fig1}
\end{figure}

When $K_1$ and $K_2$ are right handed trefoil knots $T_{2. 3}$, Y. Matsumoto asked in \cite{K} whether $M_0(T_{2. 3}, T_{2. 3})$ bounds a contractible $4$-manifold or not. By Gordon's result \cite{G}, if $n$ is odd, $M_n(T_{2. 3}, T_{2. 3}) $ does not bound any contractible $4$-manifold. If $n$ is $6$, N. Maruyama \cite{M} proved that $M_6(T_{2. 3}, T_{2. 3}) $ bounds a contractible $4$-manifold. If $n$ is $0$, S. Akbulut \cite{A} proved that $M_0(T_{2. 3}, T_{2. 3}) $ does not bound any contractible $4$-manifold.

In this note, we show that if each of $K_1$ and $K_2$ is a trefoil knot or a figure eight knot, the homology $3$-sphere defined by Figure \ref{fig1} is represented by an $n$-twisted Whitehead double of $K_2$ .\\

\noindent\textbf{Notations.}

\begin{figure}[!h]
 \begin{minipage}{0.5\hsize}
  \begin{center}
   \includegraphics[height=30mm]{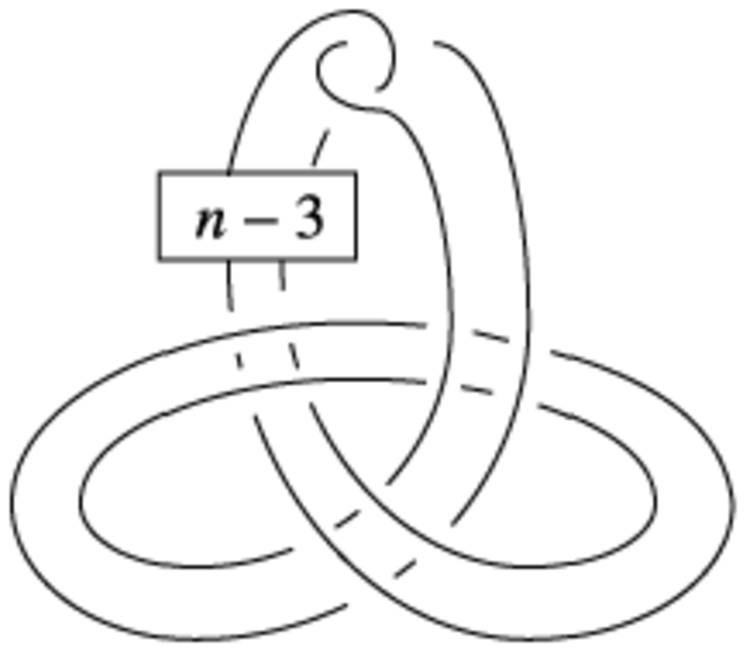}
  \end{center}
  \caption{$D_+(T_{2. 3}, n )$}
  \label{fig2}
 \end{minipage}%
 \begin{minipage}{0.5\hsize}
  \begin{center}
   \includegraphics[height=30mm]{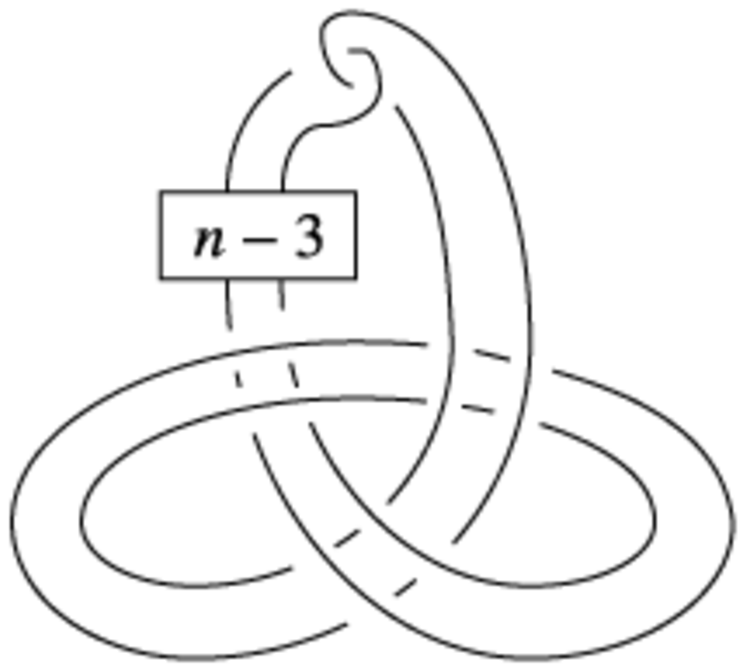}
  \end{center}
  \caption{$D_-(T_{2. 3}, n )$}
  \label{fig3}
 \end{minipage}
\end{figure}

(i).  Let $K$ be a knot, we define $D_+(K, n)$ ( or $D_-(K, n)$ ) to be the $n$-twisted Whitehead double of $K$ with a positive hook ( or a negative hook ). For example, when $K$ is a right handed trefoil knot $T_{2. 3}$, $D_+(T_{2. 3}, n )$ is the knot represented by Figure \ref{fig2}, and $D_-(T_{2. 3}, n )$ is the knot represented by Figure \ref{fig3}.

\begin{figure}[!h]
 \begin{center}
  \includegraphics[height=35mm]{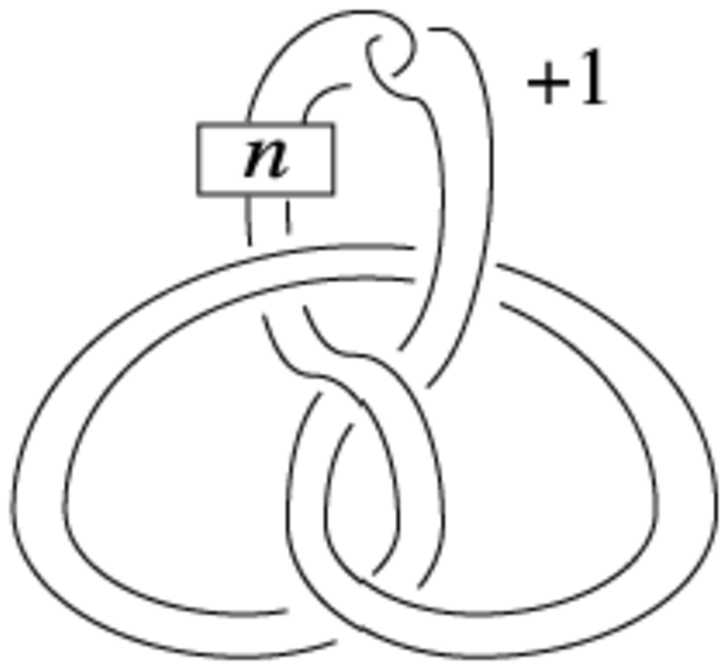}
 \end{center}
 \caption{$S^3_{+1}(D_+(K, n))$}
\label{fig4}
\end{figure}

(ii). We define $S^3_{\pm1}(K)$ to be the $\pm1$-surgery along a knot $K$. For example, when $K$ is a figure eight knot, $S^3_{+1}(D_+(K, n))$ is represented by Figure \ref{fig4}.

\begin{Theorem}\label{thm:main}
If each of $K_1$ and $K_2$ is a trefoil knot or a figure eight knot, $M_n(K_1, K_2)$ is represented by the second column on the following table. $\lambda(S^3_{\pm1}(D_\pm(K, n)))$ is the Casson invariant of $S^3_{\pm1}(D_\pm(K, n))$.

\begin{table}[H]
 \begin{center}
  \begin{tabular}{|c|c|c|} \hline 
$M_n(K_1, K_2)$ & $S^3_{\pm1}(D_\pm(K, n))$& $\lambda(S^3_{\pm1}(D_\pm(K, n)))$ \\ \hline
$K_1 :$ right handed trefoil, $K_2 :$ right handed trefoil & $S^3_{+1}(D_+(K_2, n))$ &$-n$ \\ \hline
$K_1 :$ left handed trefoil, $K_2 :$ right handed trefoil & $S^3_{-1}(D_-(K_2, n))$ &$-n$ \\ \hline
 & $S^3_{-1}(D_+(K_2, n))$ &\\
 $K_1 :$ figure eight knot, $K_2 :$ right handed trefoil & $\cong$ & $n$ \\
 & $S^3_{+1}(D_-(K_2, n))$ & \\ \hline
$K_1 :$ right handed trefoil, $K_2 :$ figure eight knot & $S^3_{+1}(D_+(K_2, n))$ &$-n$ \\ \hline
$K_1 :$ figure eight knot, $K_2 :$ figure eight knot & $S^3_{+1}(D_-(K_2, n))$ &$n$ \\ \hline
  \end{tabular}
 \end{center}
\end{table}
\end{Theorem}

We will prove \thmref{thm:main} in Section 2.\\

\noindent\textbf{Remark.}
When $n$ is 0, S. Akbulut \cite{A} shows essentially the same result of the first row on the table by a different method.

\begin{Corollary}[Gordon \cite{G}, cf. Y. Matsumoto \cite{Y} \S 3.1.]\label{corG}
Let $M_n(K_1, K_2)$ be one of the manifolds in the above table. If $n$ is odd, $M_n(K_1, K_2)$ does not bound any contractible $4$-manifold.
\end{Corollary}

\proof A short proof of this result goes as follows:

\noindent The Casson invariant, when reduced modulo $2$, is the Rohlin invariant:

\begin{center} 
$\lambda(M_n(K_1, K_2)) \equiv \mu(M_n(K_1, K_2)) \mod 2$
\end{center}

By \thmref{thm:main}, $\lambda(M_n(K_1, K_2))$ is $n$ or $-n$. Therefore if $n$ is odd, we have $\mu(M_n(K_1, K_2)) \equiv 1 \mod 2$, and so $M_n(K_1, K_2)$ does not bound any contractible $4$-manifold.
\endproof

\begin{Corollary}[N. Maruyama \cite{M}]\label{corM}
If $K_1$ and $K_2$ are right handed trefoil knots $T_{2. 3}$, $M_6(T_{2. 3}, T_{2. 3})$ bounds a contractible $4$-manifold.
\end{Corollary}

\proof
By the first row on \thmref{thm:main}'s table, $M_n(T_{2. 3}, T_{2. 3})$ is represented by $S^3_{+1}(D_+(T_{2. 3}, n))$. If $n$ is $6$, $D_+(T_{2. 3}, 6)$ is known to be a slice knot (\cite{R}, p226). Therefore by \cite{G}, $M_6(T_{2. 3}, T_{2. 3})$ bounds a contractible $4$-manifold.
\begin{figure}[H]
 \begin{minipage}{0.33\hsize}
  \begin{center}
   \includegraphics[height=25mm]{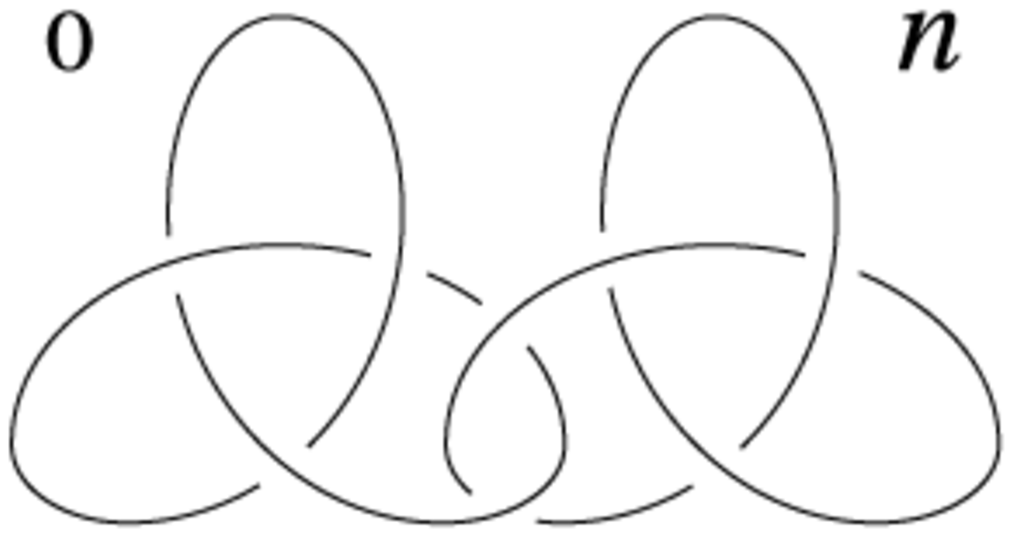}
  \end{center}
  \caption{$M_n(T_{2. 3}, T_{2. 3})$}
  \label{fig5}
 \end{minipage}%
 \begin{minipage}{0.25\hsize}
  \begin{center}
   \includegraphics[height=7mm]{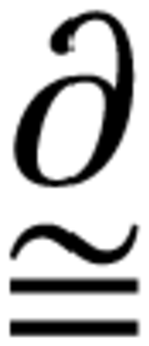}
  \end{center}
 \end{minipage}%
 \begin{minipage}{0.42\hsize}
  \begin{center}
   \includegraphics[height=25mm]{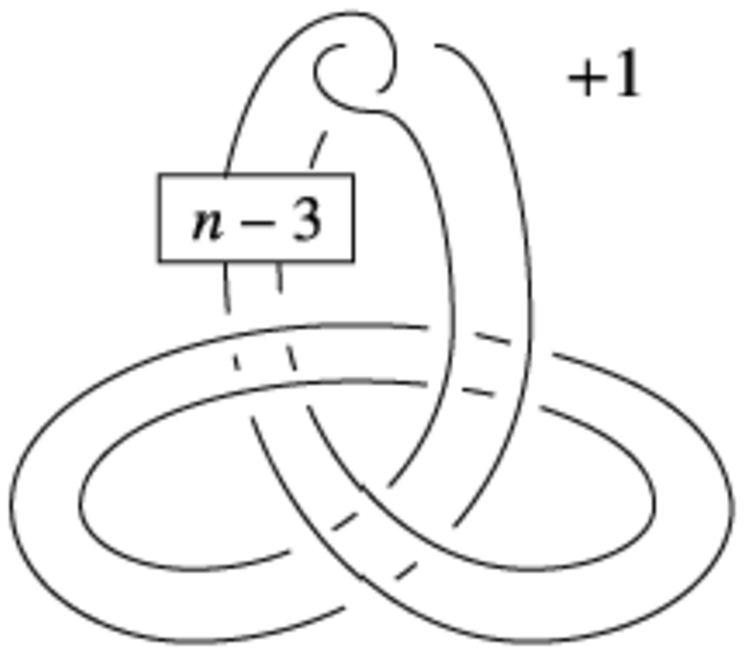}
  \end{center}
  \caption{$S^3_{+1}(D_+(T_{2. 3}, n))$}
  \label{fig6}
 \end{minipage}
\end{figure}

\endproof

\begin{Corollary}\label{cor1}
 If $n$ is $0$, $D_+(T_{2. 3}, 0)$ is not a slice knot.
\end{Corollary}

\proof
By \cite{A},  $M_0(T_{2. 3}, T_{2. 3})$ does not bound any contractible $4$-manifold. Therefore $D_+(T_{2. 3}, 0)$ is not a slice knot.
\endproof

\noindent\textbf{Remark.}
M. Hedden \cite{H} showed that if $n$ is smaller than $2$, $D_+(T_{2. 3}, n)$ is not a slice knot. 

\begin{Corollary}\label{cor3}
Let  $T_{2. 3}$ be a right handed trefoil knot and $4_1$ be a figure eight knot. The homology $3$-spheres $S^3_{-1}(D_+(T_{2. 3}, 0))$, $S^3_{+1}(D_-(T_{2. 3}, 0))$ and $S^3_{+1}(D_+(4_1, 0))$ are pairwise diffeomorphic.

\begin{figure}[H]
 \begin{minipage}{0.5\hsize}
  \begin{center}
   \includegraphics[height=25mm]{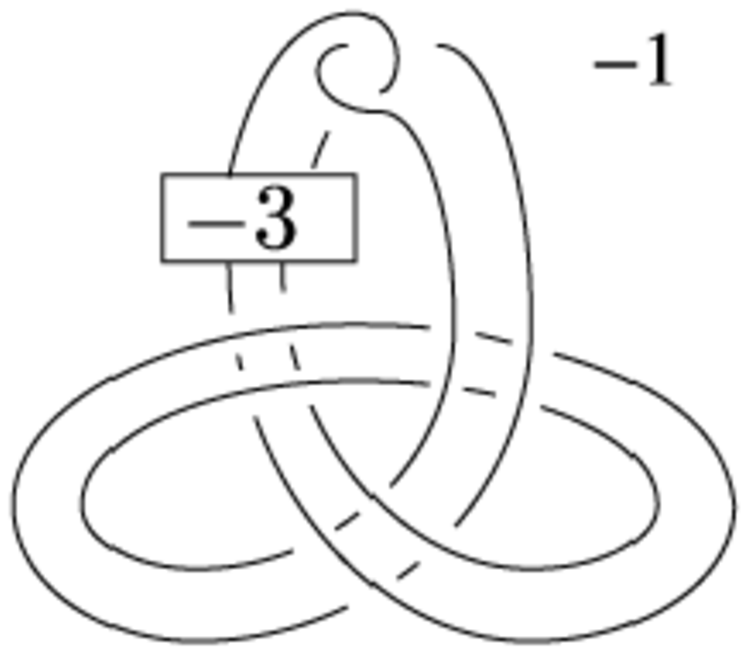}
  \end{center}
  \caption{$S^3_{-1}(D_+(T_{2. 3}, 0))$}
   \label{fig9-3}
 \end{minipage}%
 \begin{minipage}{0.5\hsize}
  \begin{center}
   \includegraphics[height=25mm]{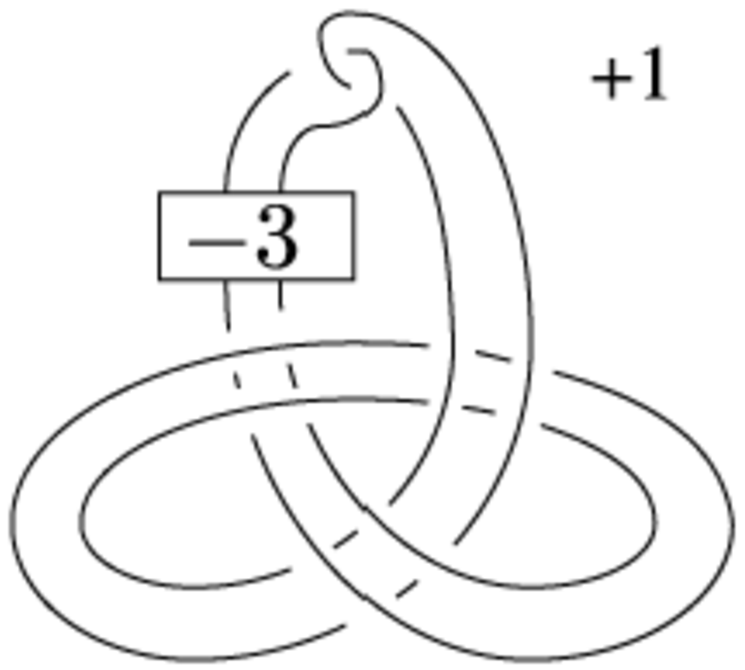}
  \end{center}
 \caption{$S^3_{+1}(D_-(T_{2. 3}, 0))$}
   \label{fig9-4}
 \end{minipage}%
\end{figure}

\begin{figure}[H]
 \begin{center}
  \includegraphics[height=25mm]{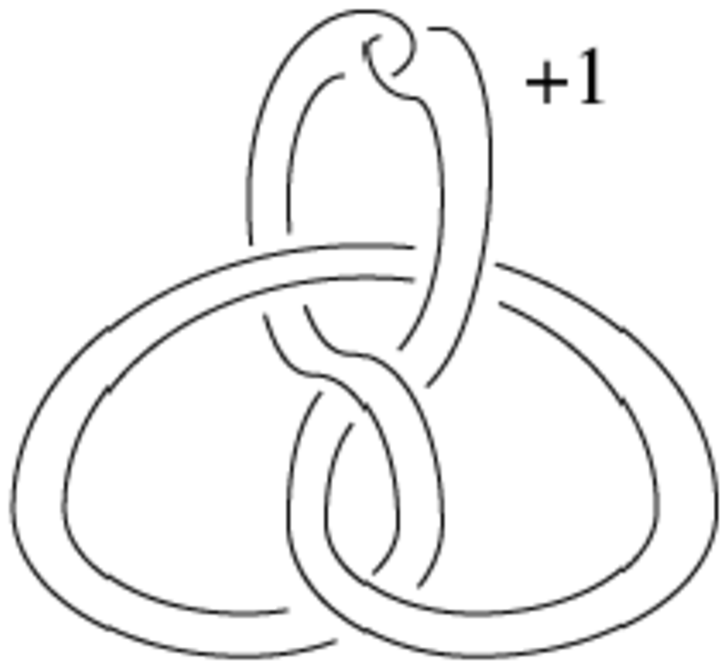}
 \end{center}
 \caption{$S^3_{+1}(D_+(4_1, 0))$}
 \label{fig9-5}
\end{figure}

\end{Corollary}

\proof
 By the third row and the fourth row on \thmref{thm:main}'s table, if $n = 0$, the $4$-dimensional handlebodies defined by Figures \ref{fig9-3}, \ref{fig9-4} and \ref{fig9-5} have the same boundaries. Therefore the homology $3$-spheres $S^3_{-1}(D_+(T_{2. 3}, 0))$, $S^3_{+1}(D_-(T_{2. 3}, 0))$ and $S^3_{+1}(D_+(4_1, 0))$ are pairwise diffeomorphic.
\endproof

\begin{figure}[H]
 \begin{minipage}{0.33\hsize}
  \begin{center}
   \includegraphics[height=25mm]{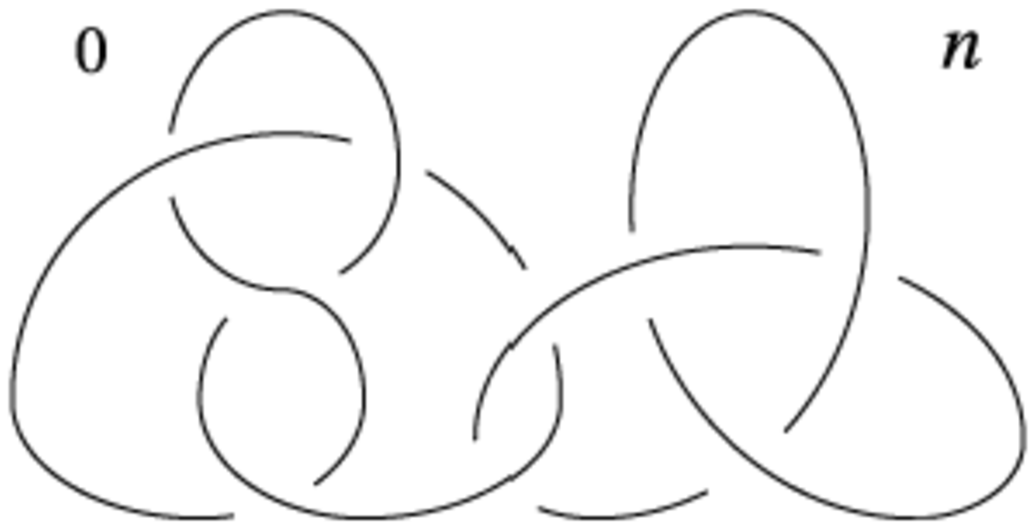}
  \end{center}
  \caption{$M_n(K_1, K_2)$}
  \label{fig7}
 \end{minipage}%
 \begin{minipage}{0.25\hsize}
  \begin{center}
   \includegraphics[height=7mm]{bdydiff.eps}
  \end{center}
 \end{minipage}%
 \begin{minipage}{0.42\hsize}
 \begin{center}
  \includegraphics[height=25mm]{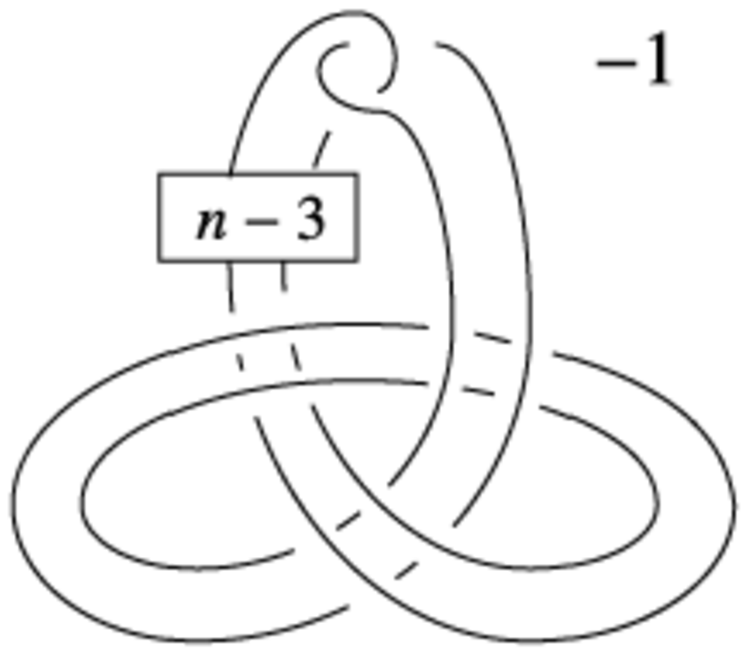}
 \end{center}
  \caption{$S^3_{-1}(D_+(K_2, n))$}
  \label{fig8}
 \end{minipage}%
\end{figure}

\begin{Corollary}\label{cor2}
If $K_1$ is a figure eight knot and $K_2$ is a right handed trefoil knot (see Figure \ref{fig7}), then $M_6(K_1, K_2)$ bounds a contractible $4$-manifold.
\end{Corollary}

\proof
 By the third row on \thmref{thm:main}'s table, $M_n(K_1, K_2)$ is represented by  $S^3_{-1}(D_+(K_2, n))$ and also by $S^3_{+1}(D_-(K_2, n))$. If $n$ is $6$, $D_+(K_2, 6)$ is known to be a slice knot (\cite{R}, p226). Therefore by \cite{G}, $M_6(K_1, K_2)$ bounds a contractible $4$-manifold.

\begin{figure}[H]
 \begin{minipage}{0.3\hsize}
  \begin{center}
   \includegraphics[height=7mm]{bdydiff.eps}
  \end{center}
 \end{minipage}%
 \begin{minipage}{0.4\hsize}
 \begin{center}
  \includegraphics[height=25mm]{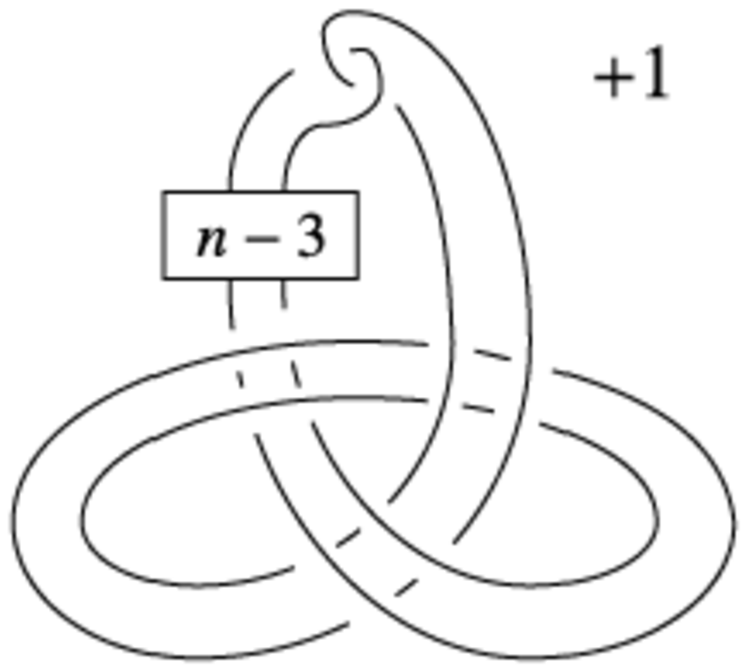}
 \end{center}
  \caption{$S^3_{+1}(D_-(K_2, n))$}
  \label{fig9}
 \end{minipage}
\end{figure}

\endproof

\noindent\textbf{Remark.}
By Corollary \ref{cor2}, the homology $3$-sphere $S^3_{+1}(D_-(T_{2,3}, 6))$ bounds a contractible $4$-manifold. The author does not know whether the knot $D_-(T_{2,3}, 6)$ is a slice knot or not.
\begin{figure}[H]
 \begin{center}
  \includegraphics[height=25mm]{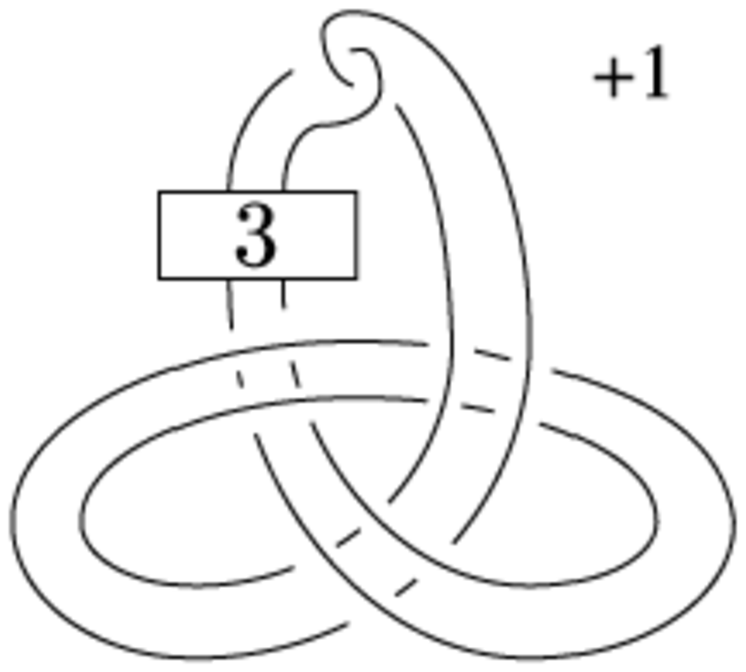}
 \end{center}
 \caption{$S^3_{+1}(D_-(T_{2,3}, 6))$}
 \label{fig9-2}
\end{figure}

\noindent\textbf{Question.}
Let $V^1_n$ be the $4$-dimensional handlebody defined by Figure \ref{fig9}, and $V^2_n$ be the $4$-dimensional handlebody defined by Figure \ref{fig8}. Since $\partial(V^1_n)$ is diffeomorphic to $\partial(V^2_n)$ by \thmref{thm:main}, we have a closed $4$-manifold $V^1_n \cup_{\partial} (-V^2_n)$. Because $D_+(T_{2,3}, 6)$ is a slice knot, we have a smooth $S^2$ with self intersection $-1$ in $V^2_6$ representing a generator of $H_2(V^2_6)$. Blow down this smooth $S^2$ from the $V^1_6 \cup_{\partial} (-V^2_6)$. Then we are left with a closed smooth $4$-manifold homotopy equivalent to $\CP^2$. Is this $4$-manifold diffeomorphic to $\CP^2$ ?

\begin{Proposition}\label{prop2}
$V^1_n \cup_{\partial} (-V^2_n)$ is diffeomorphic to $\CP^2 \sharp \CP^2$. 
\end{Proposition}

We show this fact in Section 3.\\

It seems that \thmref{thm:main} is related to \cite{M} Corollary 8 (3), but the author could not understand the relationship clearly.

The author does not know whether there is an even number $n$ $\neq$ $0, 6$, such that $M_n(T_{2,3}, T_{2,3})$ bounds a contractible $4$-manifold or not. M. Tange \cite{T} proved that if $n$ is smaller than $2$, $M_n(T_{2,3}, T_{2,3})$ does not bound any contractible $4$-manifold by computing the Heegaard Floer homology $HF^+(M_n(T_{2,3}, T_{2,3}))$ and the correction term $d(M_n(T_{2,3}, T_{2,3}))$.  

\begin{Acknowledgements}
The author would like to thank Yukio Matsumoto and Nobuhiro Nakamura for their useful comments and encouragement. The author is very grateful to the referee for his/her careful reading and useful comments.
\end{Acknowledgements}
%
%
%
\section{Proof of \thmref{thm:main}}\label{sec:constraint}
%
%
In this section, first we show that $M_n(K_1, K_2)$ is represented by $S^3_{\pm1}(D_\pm(K_2, n))$. Next we compute the Casson invariant $\lambda(S^3_{\pm1}(D_\pm(K_2, n)))$.
\subsection{Proof of the first row on \thmref{thm:main}'s table}
$K_1$ and $K_2$ are right handed trefoil knots. 

\proof
We show that the $4$-manifolds represented by Figures \ref{fig10} and \ref{fig26} have the same boundaries by following Kirby Calculus:
\begin{figure}[H]
 \begin{minipage}{0.45\hsize}
  \begin{center}
   \includegraphics[height=20mm]{Figure5.eps}
  \end{center}
  \caption{$M_n(K_1, K_2)$}
   \label{fig10}
 \end{minipage}%
 \begin{minipage}{0.1\hsize}
  \begin{center}
   \includegraphics[height=10mm]{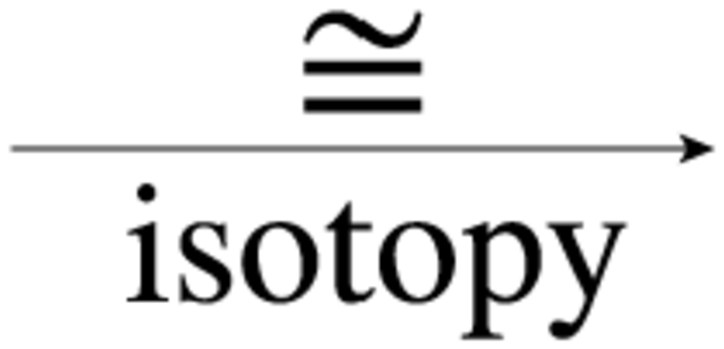}
  \end{center}
 \end{minipage}%
 \begin{minipage}{0.45\hsize}
  \begin{center}
   \includegraphics[height=20mm]{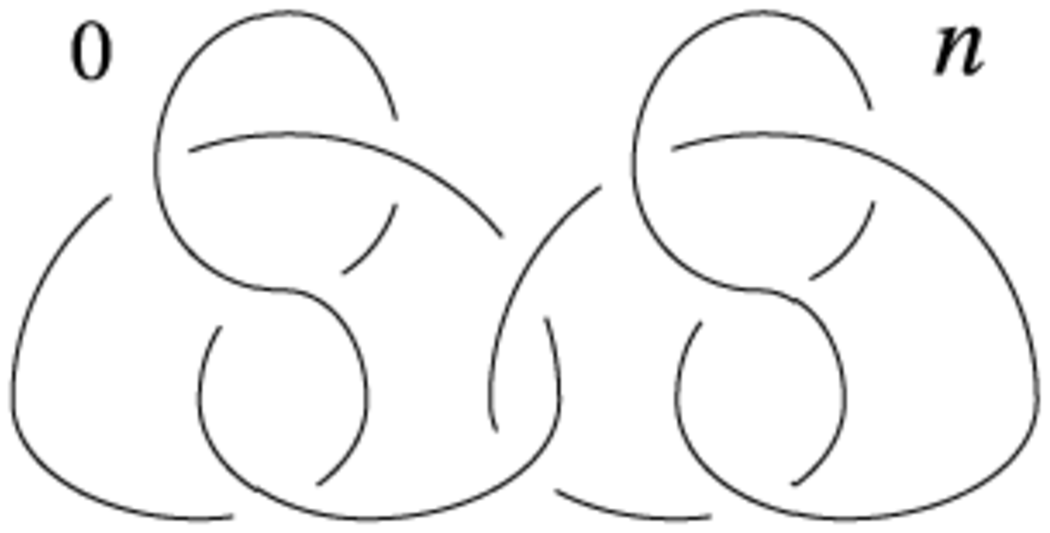}
  \end{center}
 \caption{}
   \label{fig11}
 \end{minipage}
\end{figure}

\begin{figure}[H]
 \begin{minipage}{0.1\hsize}
  \begin{center}
   \includegraphics[height=10mm]{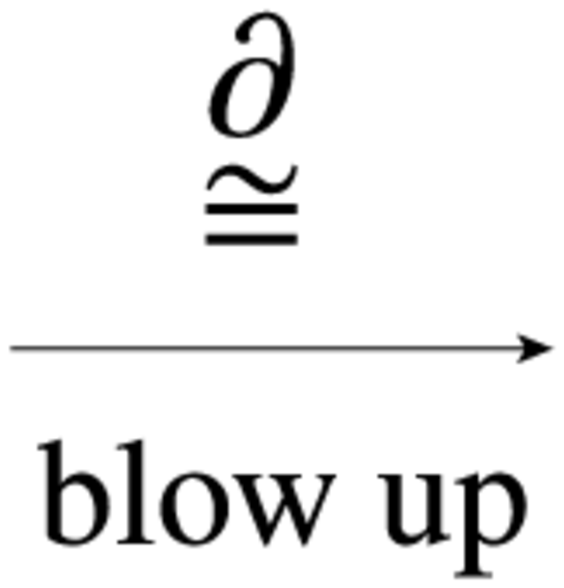}
  \end{center}
 \end{minipage}%
 \begin{minipage}{0.4\hsize}
  \begin{center}
   \includegraphics[height=20mm]{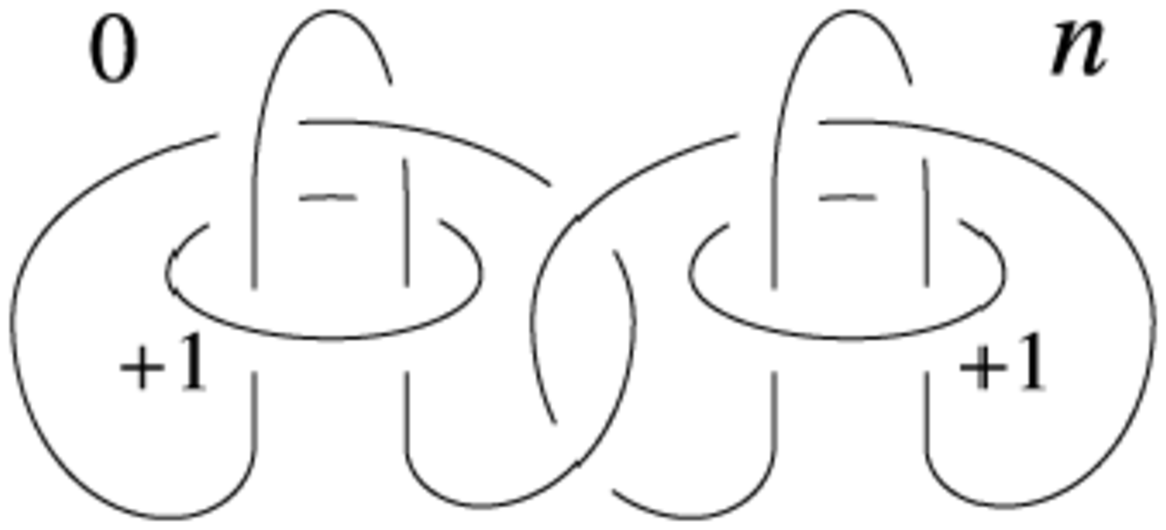}
  \end{center}
  \caption{}
   \label{fig12}
 \end{minipage}%
 \begin{minipage}{0.1\hsize}
  \begin{center}
   \includegraphics[height=10mm]{blowup.eps}
  \end{center}
 \end{minipage}%
 \begin{minipage}{0.4\hsize}
  \begin{center}
   \includegraphics[height=20mm]{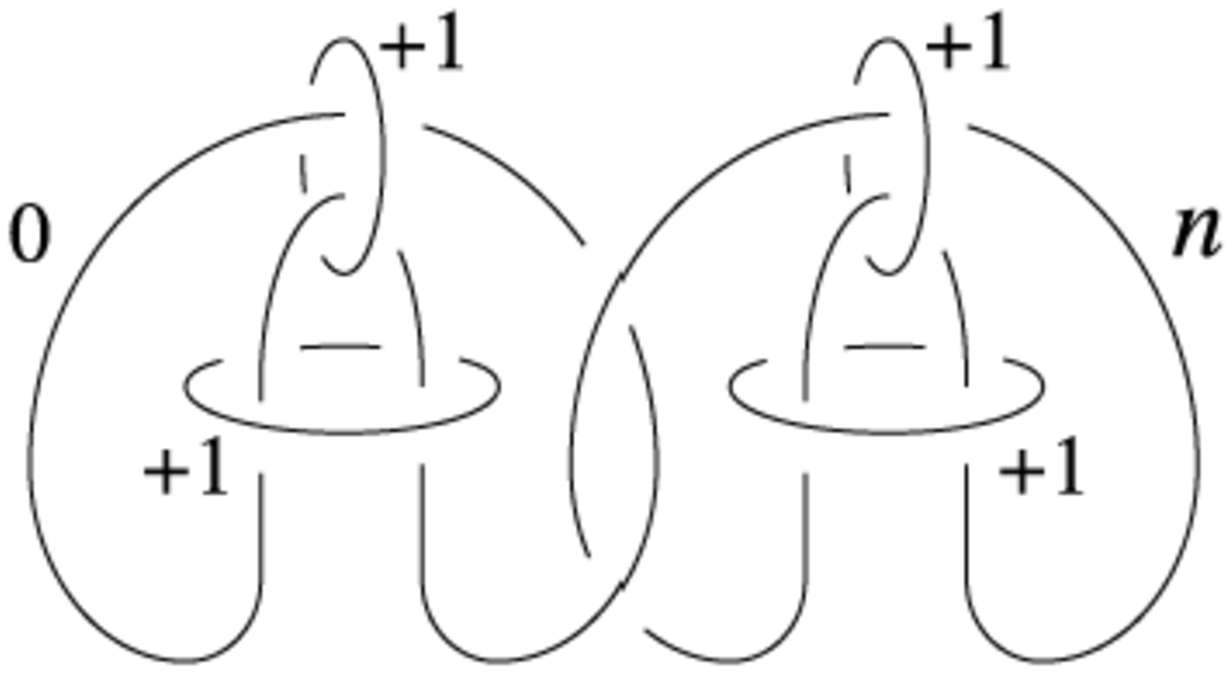}
  \end{center}
  \caption{}
   \label{fig13}
 \end{minipage}
\end{figure}

\begin{figure}[H]
 \begin{minipage}{0.1\hsize}
  \begin{center}
   \includegraphics[height=10mm]{isotopy.eps}
  \end{center}
 \end{minipage}%
 \begin{minipage}{0.4\hsize}
  \begin{center}
   \includegraphics[height=30mm]{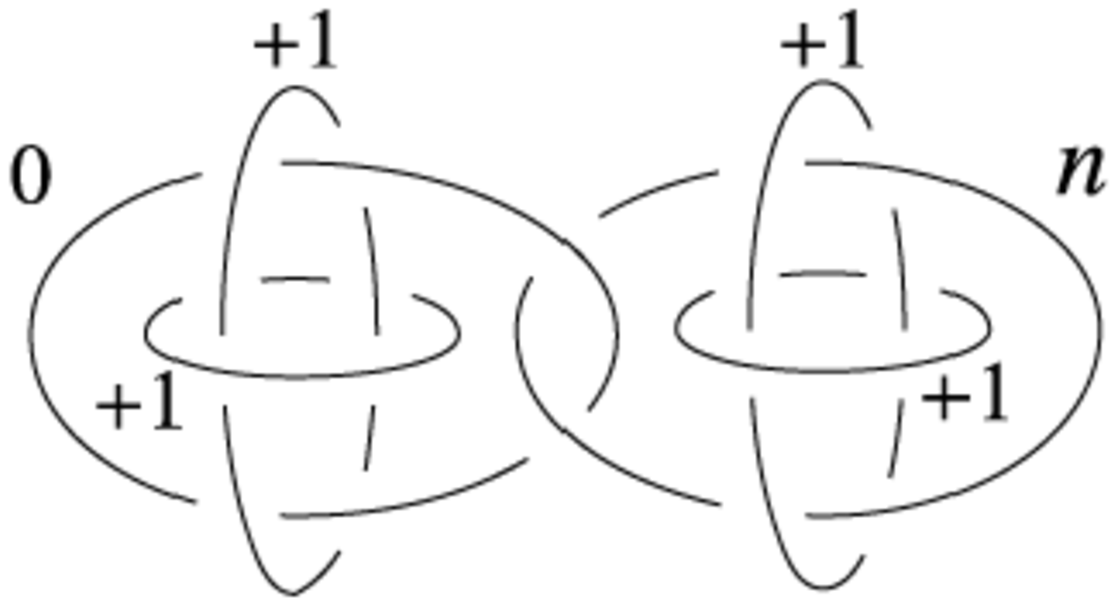}
  \end{center}
  \caption{}
   \label{fig14}
 \end{minipage}%
 \begin{minipage}{0.1\hsize}
  \begin{center}
   \includegraphics[height=10mm]{isotopy.eps}
  \end{center}
 \end{minipage}%
 \begin{minipage}{0.4\hsize}
  \begin{center}
   \includegraphics[height=30mm]{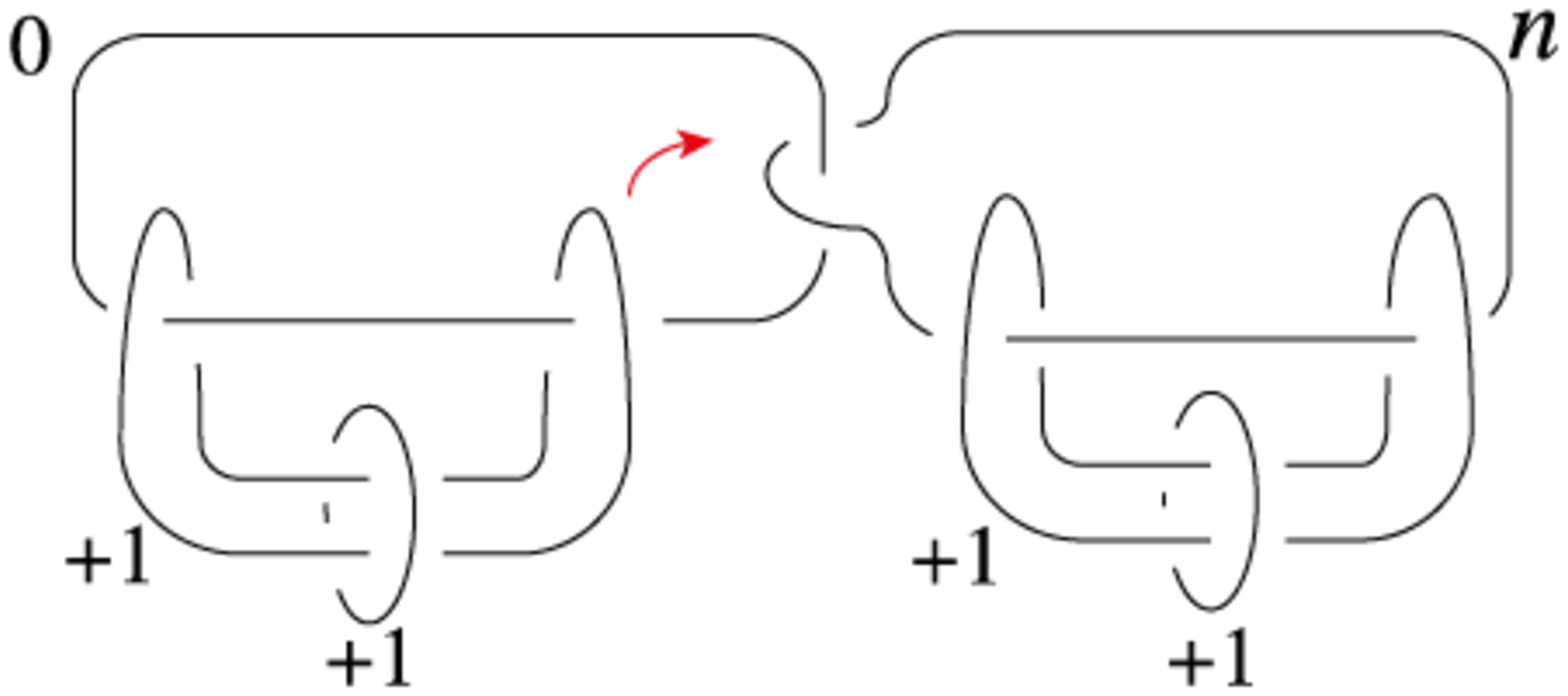}
  \end{center}
  \caption{}
   \label{fig15}
 \end{minipage}
\end{figure}

\begin{figure}[H]
 \begin{minipage}{0.5\hsize}
  \begin{center}
   \includegraphics[height=15mm]{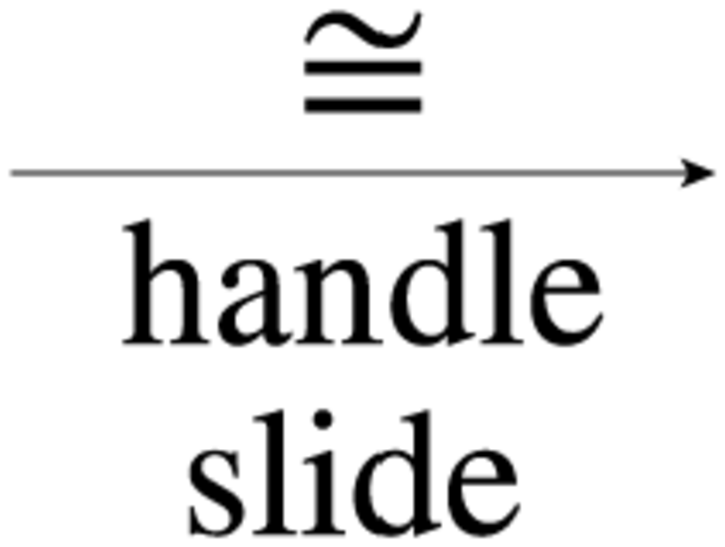}
  \end{center}
 \end{minipage}%
 \begin{minipage}{0.5\hsize}
  \begin{center}
   \includegraphics[height=35mm]{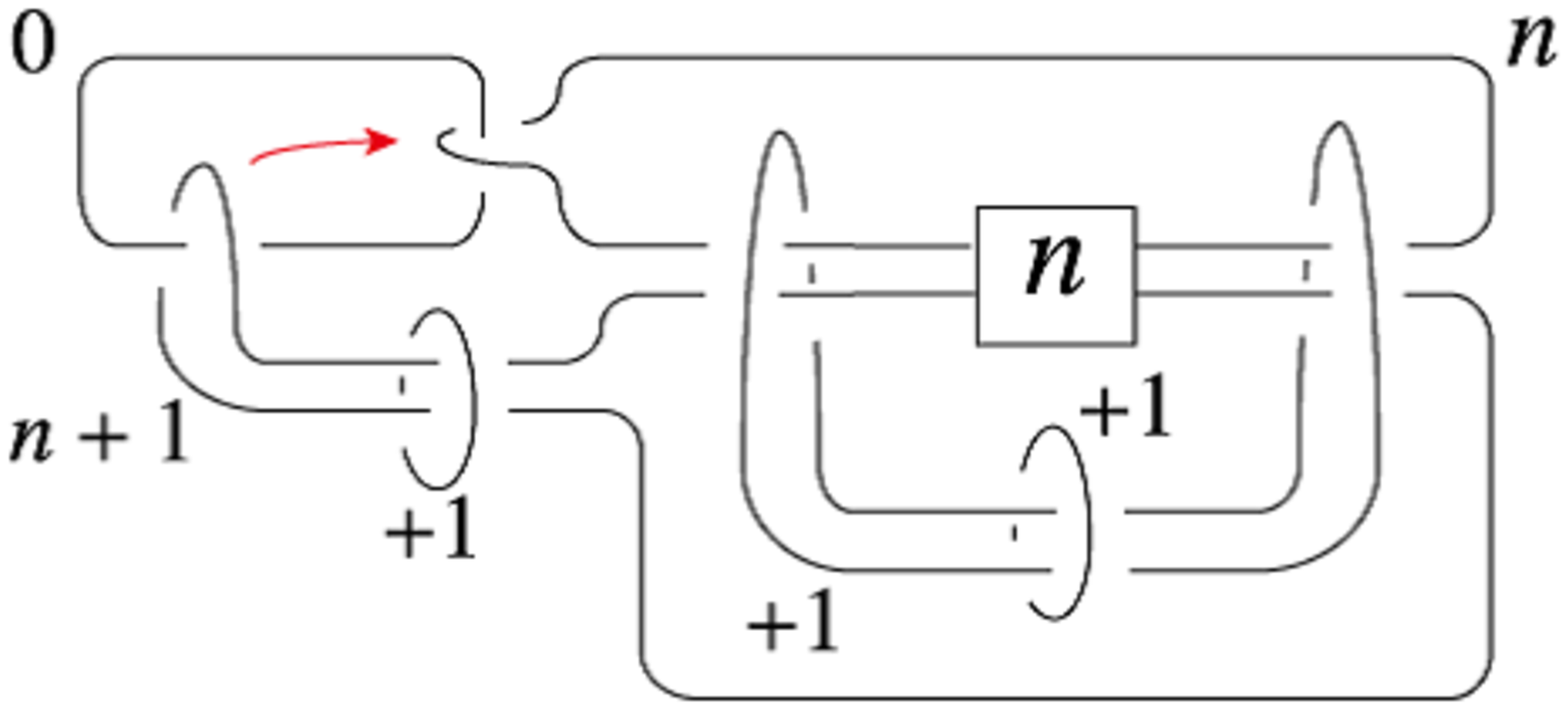}
  \end{center}
  \caption{}
  \label{fig16}
 \end{minipage}
\end{figure}

\begin{figure}[H]
 \begin{minipage}{0.5\hsize}
  \begin{center}
   \includegraphics[height=15mm]{hds.eps}
  \end{center}
 \end{minipage}%
 \begin{minipage}{0.5\hsize}
  \begin{center}
   \includegraphics[height=35mm]{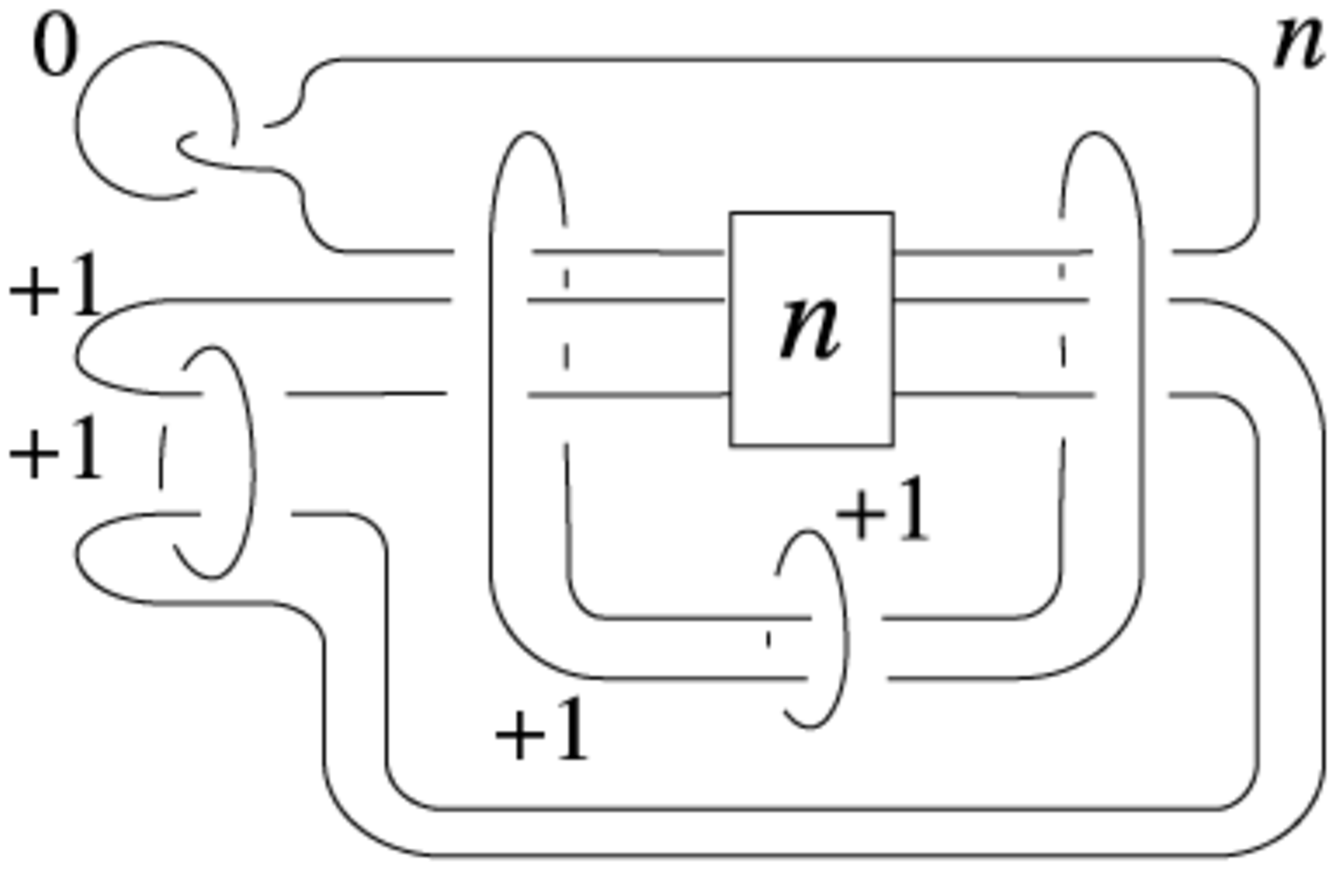}
  \end{center}
  \caption{}
  \label{fig17}
 \end{minipage}
\end{figure}

\begin{figure}[H]
 \begin{minipage}{0.1\hsize}
  \begin{center}
   \includegraphics[height=15mm]{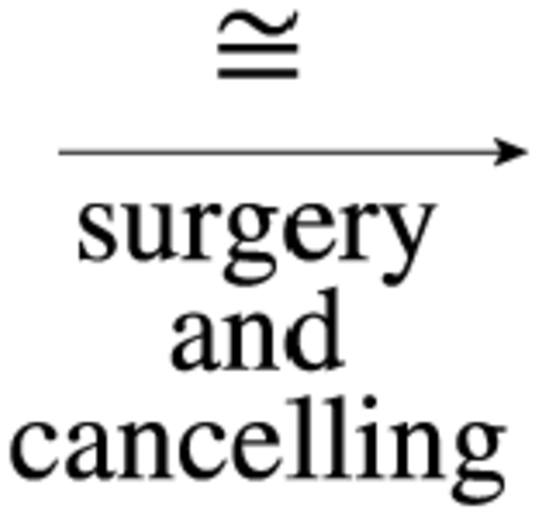}
  \end{center}
 \end{minipage}%
 \begin{minipage}{0.4\hsize}
  \begin{center}
   \includegraphics[height=30mm]{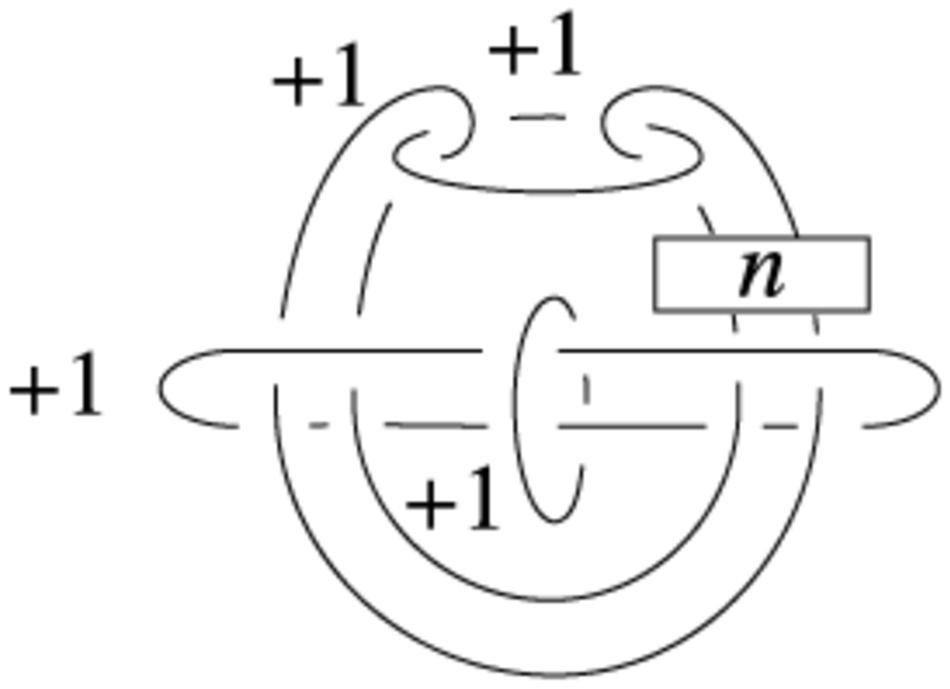}
  \end{center}
  \caption{}
   \label{fig18}
 \end{minipage}%
 \begin{minipage}{0.15\hsize}
  \begin{center}
   \includegraphics[height=10mm]{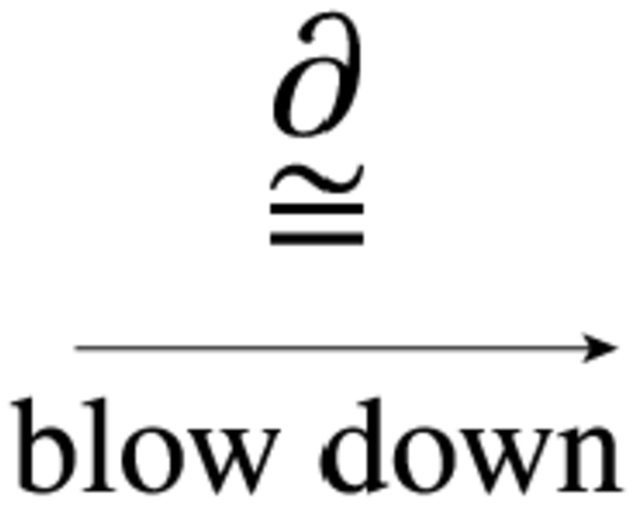}
  \end{center}
 \end{minipage}%
 \begin{minipage}{0.35\hsize}
  \begin{center}
   \includegraphics[height=30mm]{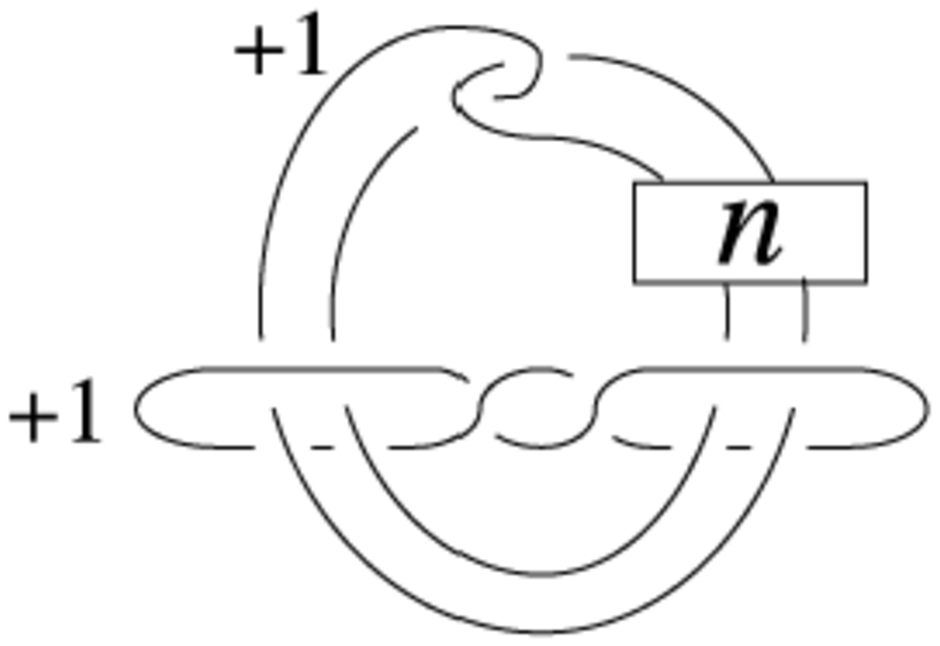}
  \end{center}
  \caption{}
   \label{fig19}
 \end{minipage}
\end{figure}

\begin{figure}[H]
 \begin{minipage}{0.1\hsize}
  \begin{center}
   \includegraphics[height=10mm]{isotopy.eps}
  \end{center}
 \end{minipage}%
 \begin{minipage}{0.4\hsize}
  \begin{center}
   \includegraphics[height=30mm]{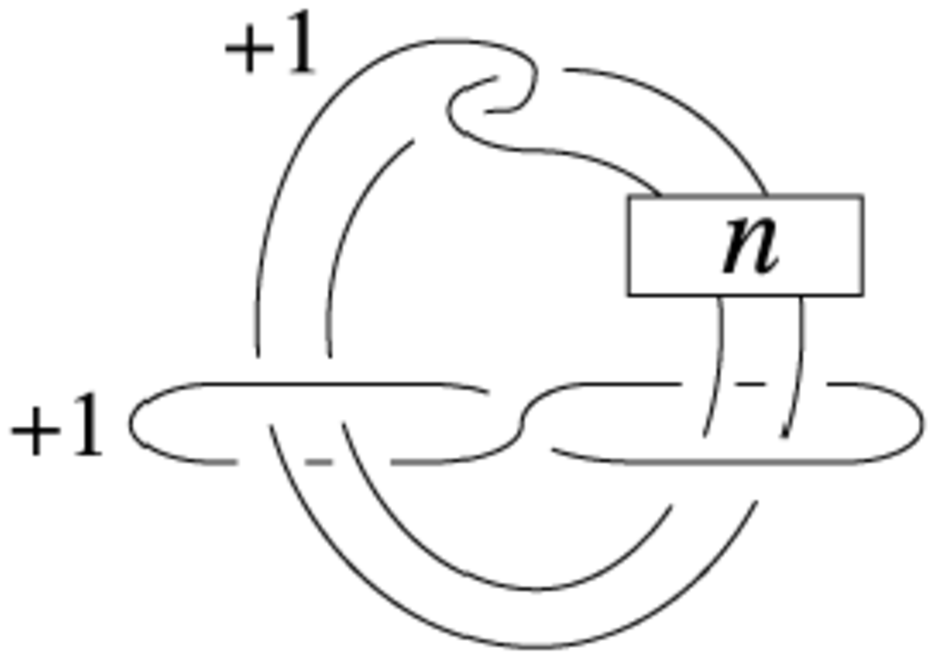}
  \end{center}
  \caption{}
   \label{fig20}
 \end{minipage}%
 \begin{minipage}{0.1\hsize}
  \begin{center}
   \includegraphics[height=10mm]{isotopy.eps}
  \end{center}
 \end{minipage}%
 \begin{minipage}{0.4\hsize}
  \begin{center}
   \includegraphics[height=35mm]{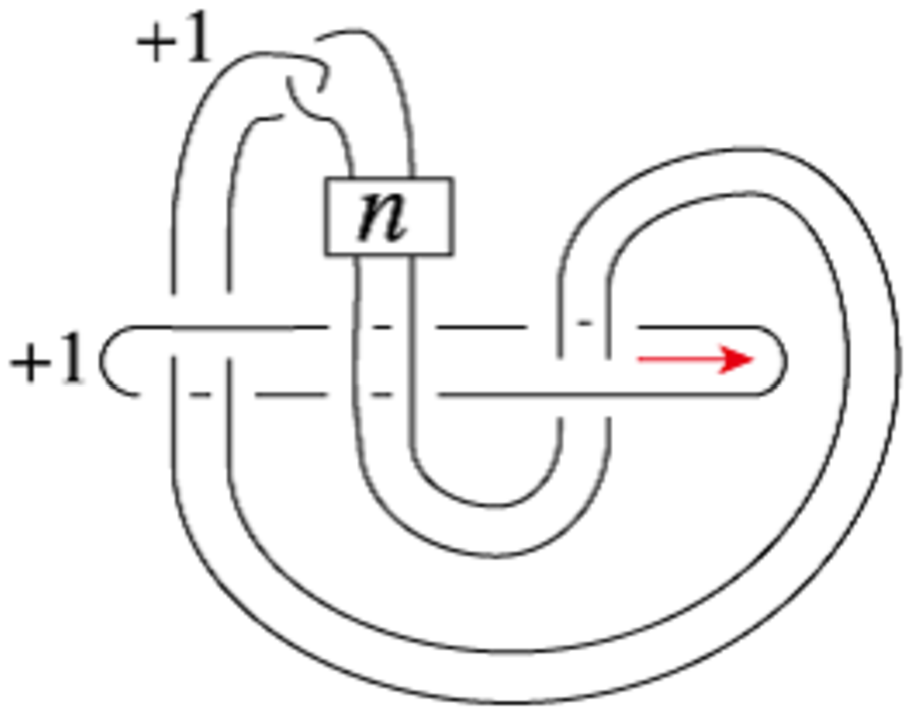}
  \end{center}
  \caption{}
   \label{fig21}
 \end{minipage}
\end{figure}

\begin{figure}[H]
 \begin{minipage}{0.1\hsize}
  \begin{center}
   \includegraphics[height=10mm]{hds.eps}
  \end{center}
 \end{minipage}%
 \begin{minipage}{0.4\hsize}
  \begin{center}
   \includegraphics[height=35mm]{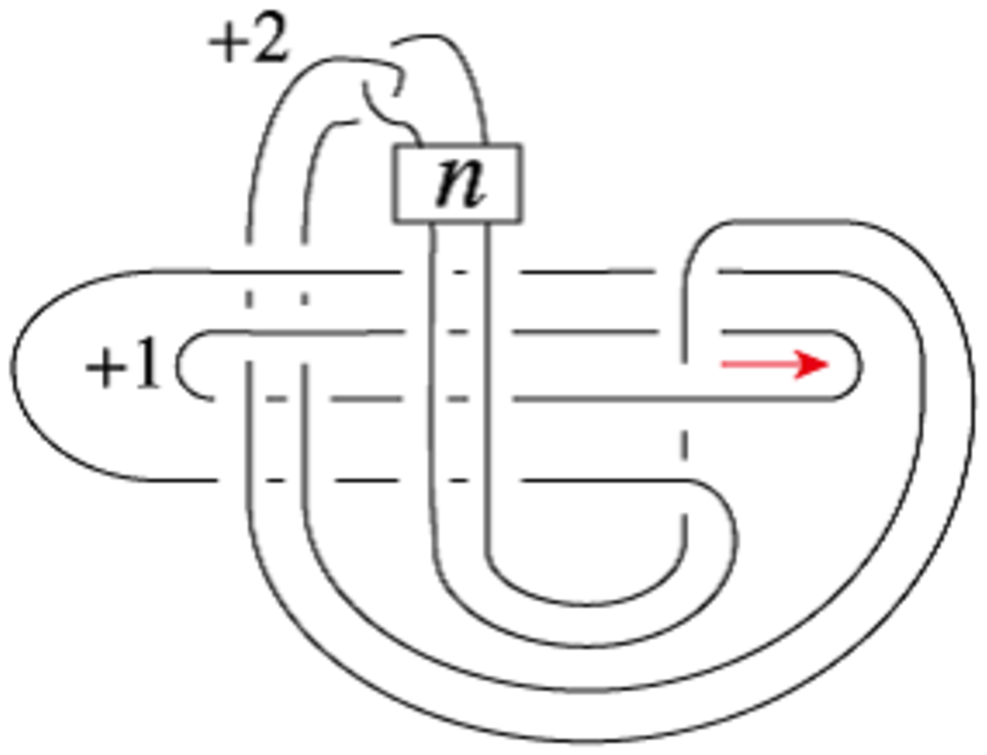}
  \end{center}
  \caption{}
   \label{fig22}
 \end{minipage}%
 \begin{minipage}{0.1\hsize}
  \begin{center}
   \includegraphics[height=10mm]{hds.eps}
  \end{center}
 \end{minipage}%
 \begin{minipage}{0.4\hsize}
  \begin{center}
   \includegraphics[height=35mm]{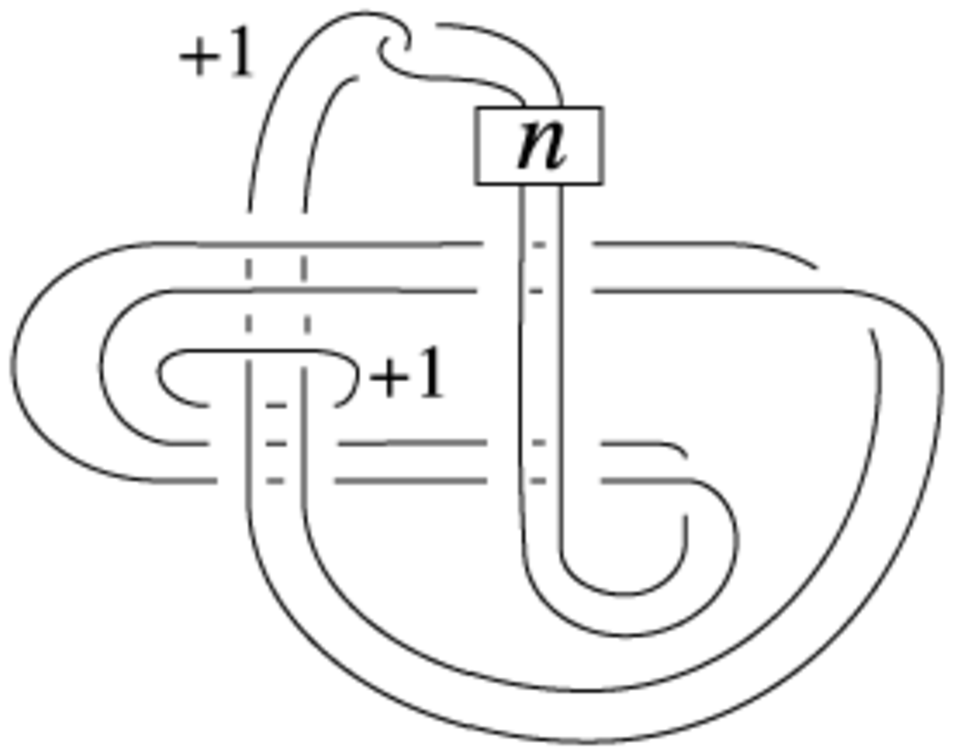}
  \end{center}
  \caption{}
   \label{fig23}
 \end{minipage}
\end{figure}

\begin{figure}[H]
 \begin{minipage}{0.1\hsize}
  \begin{center}
   \includegraphics[height=10mm]{blowdown.eps}
  \end{center}
 \end{minipage}%
 \begin{minipage}{0.4\hsize}
  \begin{center}
   \includegraphics[height=35mm]{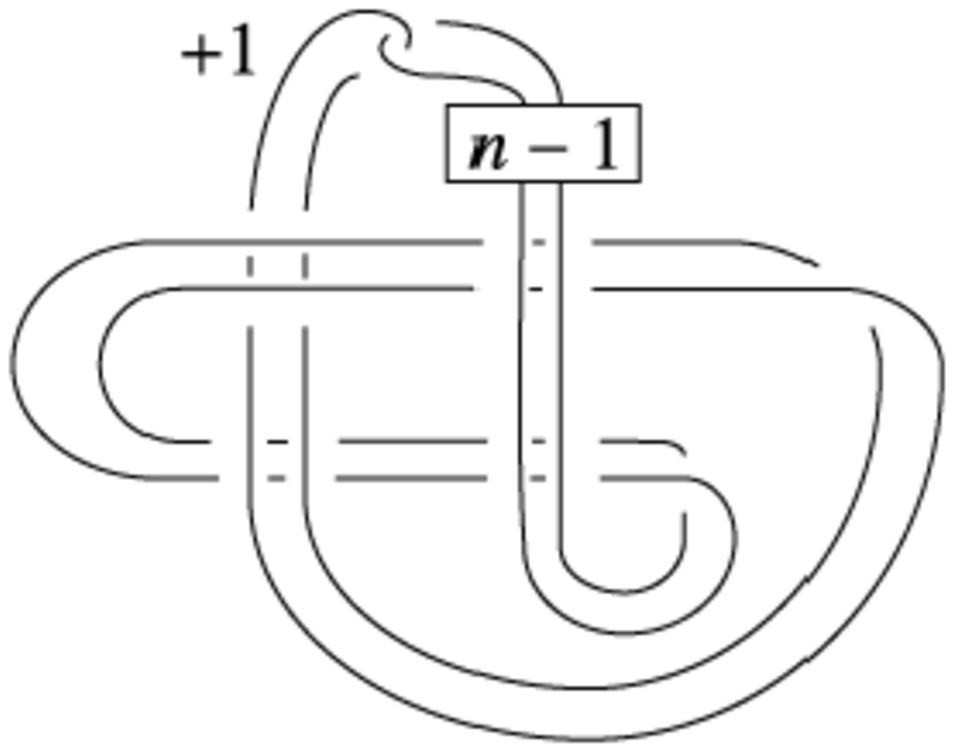}
  \end{center}
  \caption{}
   \label{fig24}
 \end{minipage}%
 \begin{minipage}{0.1\hsize}
  \begin{center}
   \includegraphics[height=10mm]{isotopy.eps}
  \end{center}
 \end{minipage}%
 \begin{minipage}{0.4\hsize}
  \begin{center}
   \includegraphics[height=35mm]{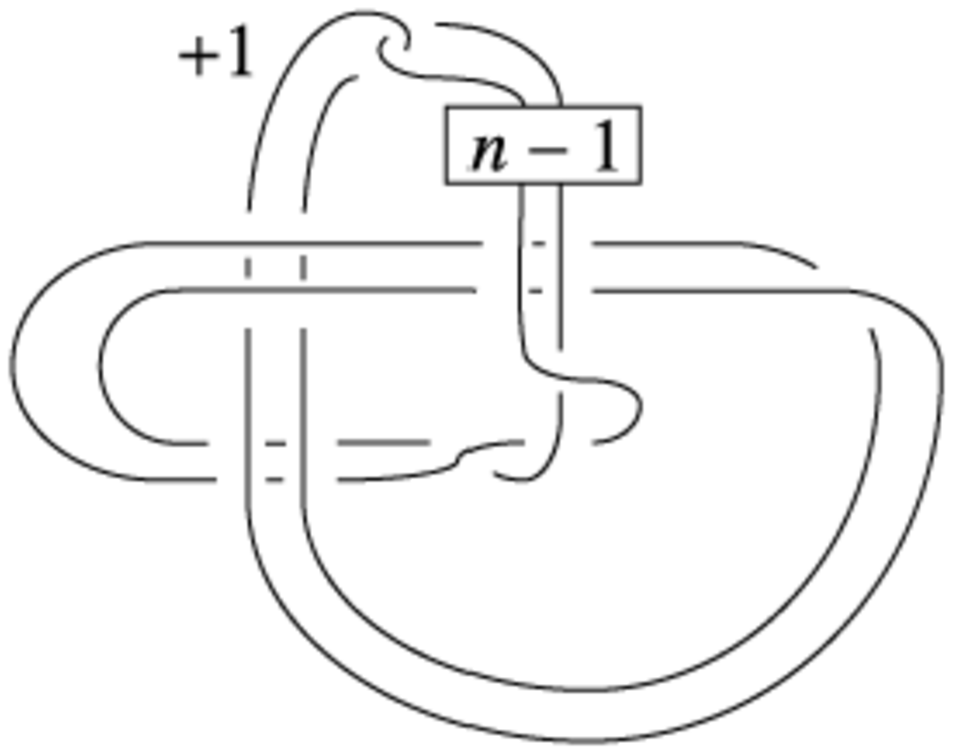}
  \end{center}
  \caption{}
   \label{fig25}
 \end{minipage}
\end{figure}

\begin{figure}[H]
 \begin{minipage}{0.5\hsize}
  \begin{center}
   \includegraphics[height=10mm]{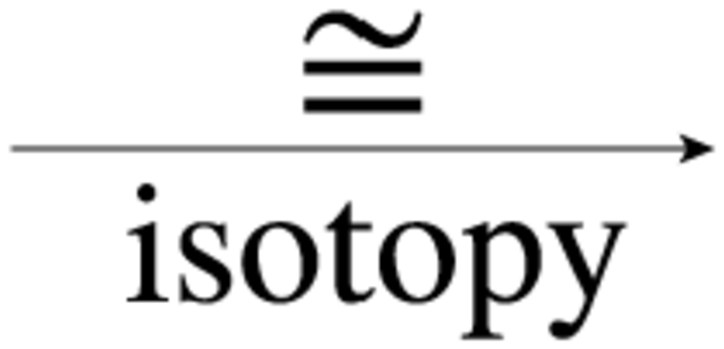}
  \end{center}
 \end{minipage}%
 \begin{minipage}{0.5\hsize}
  \begin{center}
   \includegraphics[height=30mm]{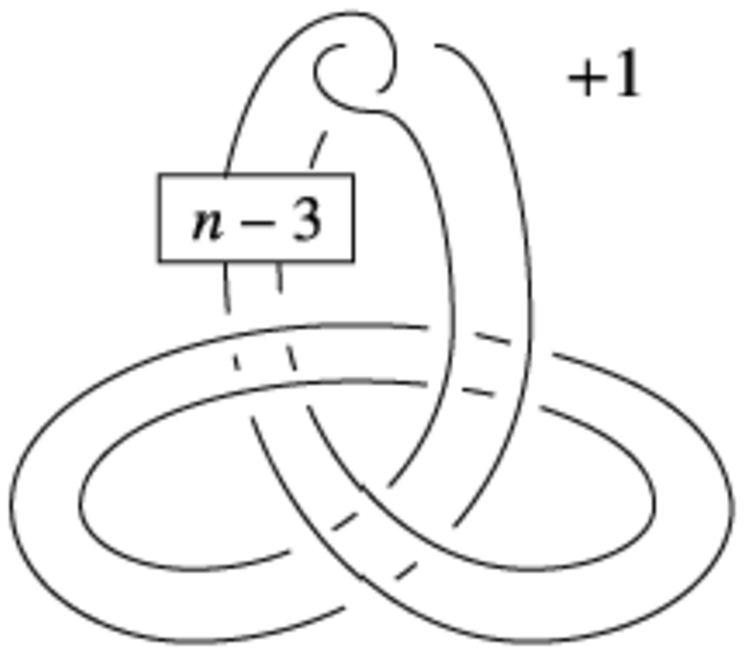}
  \end{center}
  \caption{$S^3_{+1}(D_+(K_2, n))$}
  \label{fig26}
 \end{minipage}
\end{figure}

\endproof

\subsection{Proof of the second row on \thmref{thm:main}'s table}
$K_1$ is a left handed trefoil knot and $K_2$ is a right handed trefoil knot.

\begin{figure}[H]
 \begin{minipage}{0.45\hsize}
  \begin{center}
   \includegraphics[height=20mm]{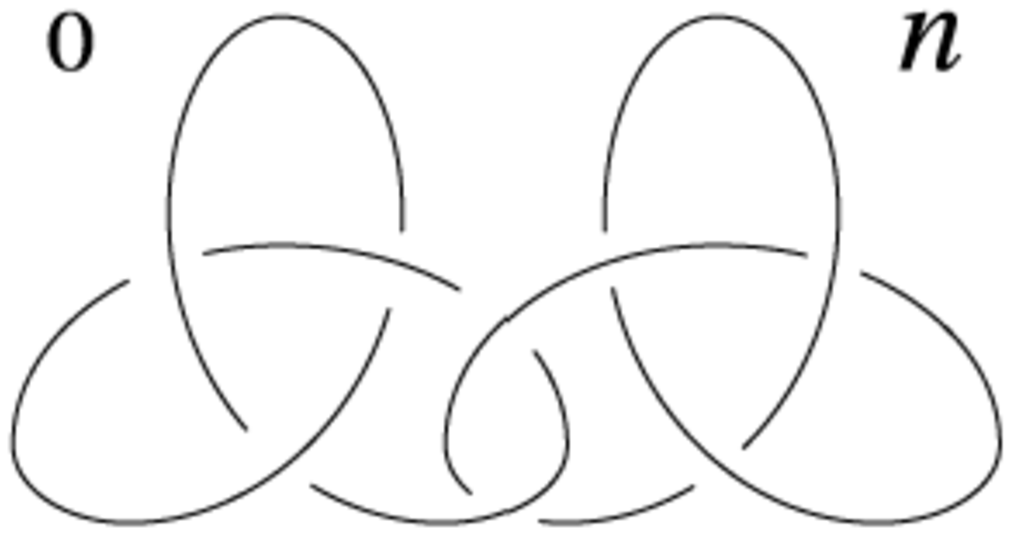}
  \end{center}
  \caption{$M_n(K_1, K_2)$}
   \label{fig27}
 \end{minipage}%
 \begin{minipage}{0.1\hsize}
  \begin{center}
   \includegraphics[height=10mm]{isotopy.eps}
  \end{center}
 \end{minipage}%
 \begin{minipage}{0.45\hsize}
  \begin{center}
   \includegraphics[height=20mm]{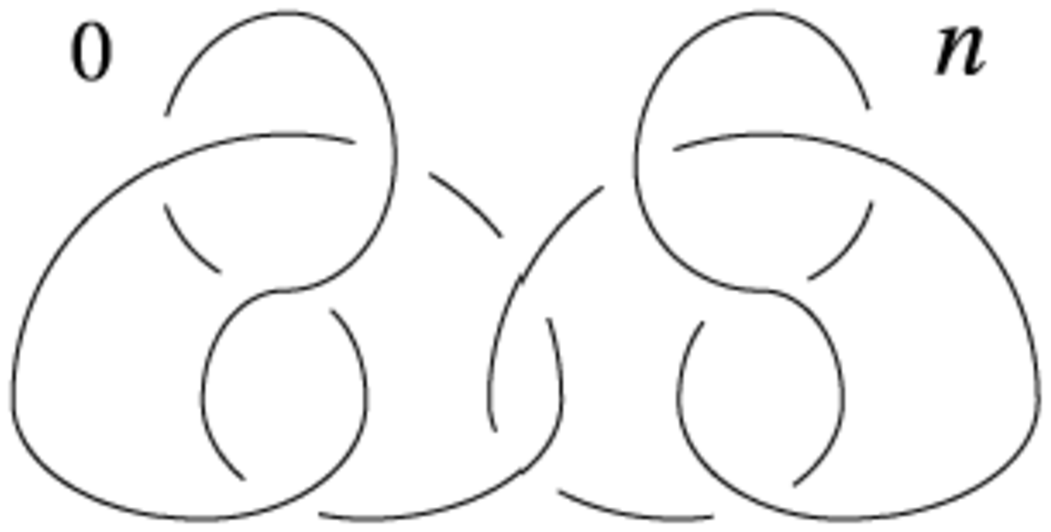}
  \end{center}
 \caption{}
   \label{fig28}
 \end{minipage}
\end{figure}

\proof
We show that the $4$-manifolds represented by Figures \ref{fig27} and \ref{fig32} have the same boundaries.

\begin{figure}[H]
 \begin{minipage}{0.1\hsize}
  \begin{center}
   \includegraphics[height=10mm]{blowup.eps}
  \end{center}
 \end{minipage}%
 \begin{minipage}{0.4\hsize}
  \begin{center}
   \includegraphics[height=20mm]{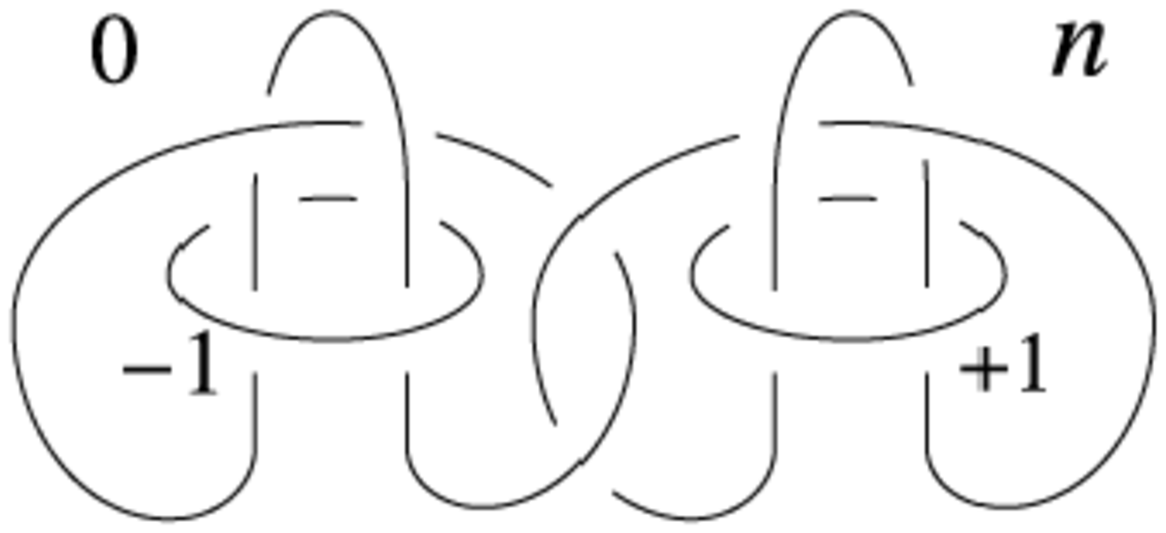}
  \end{center}
  \caption{}
   \label{fig29}
 \end{minipage}%
 \begin{minipage}{0.1\hsize}
  \begin{center}
   \includegraphics[height=10mm]{blowup.eps}
  \end{center}
 \end{minipage}%
 \begin{minipage}{0.4\hsize}
  \begin{center}
   \includegraphics[height=20mm]{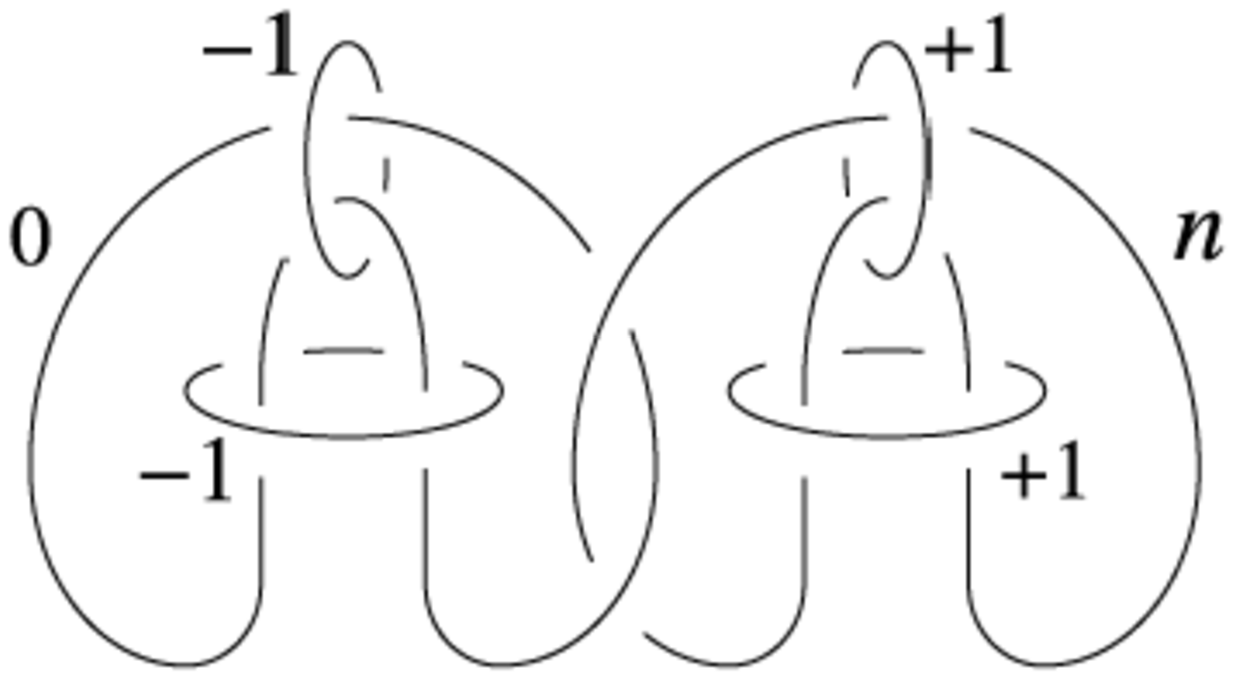}
  \end{center}
  \caption{}
   \label{fig30}
 \end{minipage}
\end{figure}

\begin{figure}[H]
 \begin{minipage}{0.1\hsize}
  \begin{center}
   \includegraphics[height=10mm]{isotopy.eps}
  \end{center}
 \end{minipage}%
 \begin{minipage}{0.3\hsize}
  \begin{center}
   \includegraphics[height=25mm]{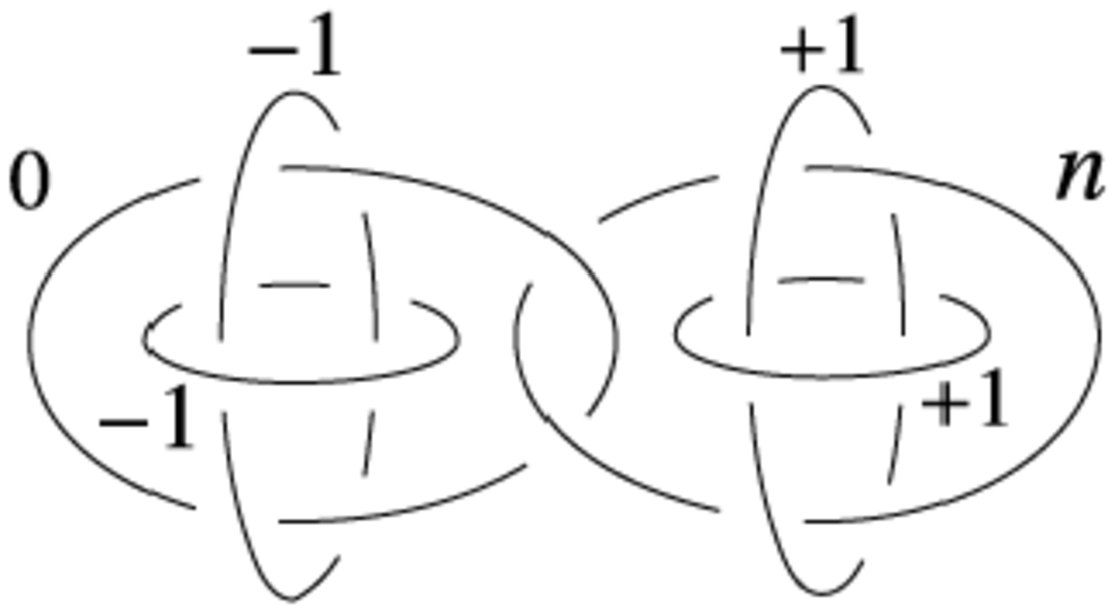}
  \end{center}
  \caption{}
   \label{fig29}
 \end{minipage}%
 \begin{minipage}{0.25\hsize}
  \begin{center}
   \includegraphics[height=25mm]{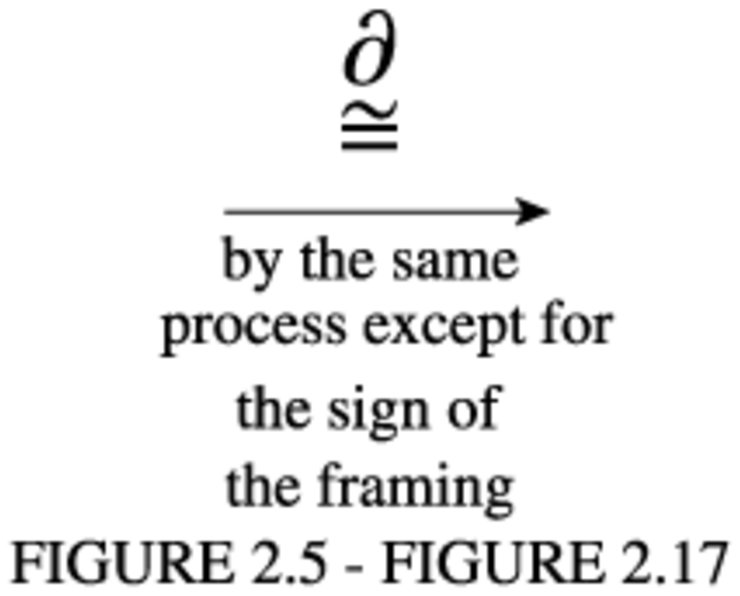}
  \end{center}
 \end{minipage}%
 \begin{minipage}{0.35\hsize}
  \begin{center}
   \includegraphics[height=30mm]{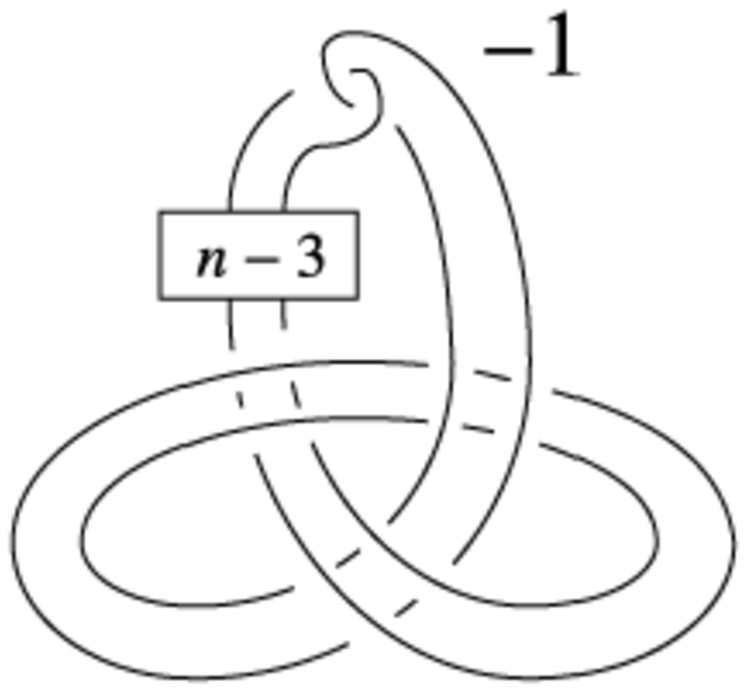}
  \end{center}
  \caption{$S^3_{-1}(D_-(K_2, n))$}
   \label{fig32}
 \end{minipage}
\end{figure}

\endproof

\subsection{Proof of the third row on \thmref{thm:main}'s table}
\begin{figure}[H]
 \begin{minipage}{0.45\hsize}
  \begin{center}
   \includegraphics[height=20mm]{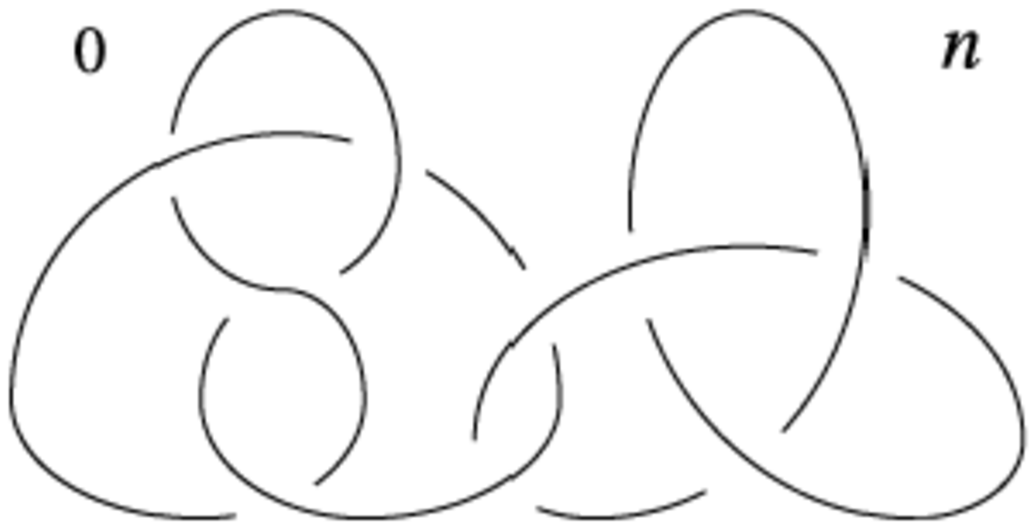}
  \end{center}
  \caption{$M_n(K_1, K_2)$}
   \label{fig33}
 \end{minipage}%
 \begin{minipage}{0.1\hsize}
  \begin{center}
   \includegraphics[height=10mm]{isotopy.eps}
  \end{center}
 \end{minipage}%
 \begin{minipage}{0.45\hsize}
  \begin{center}
   \includegraphics[height=20mm]{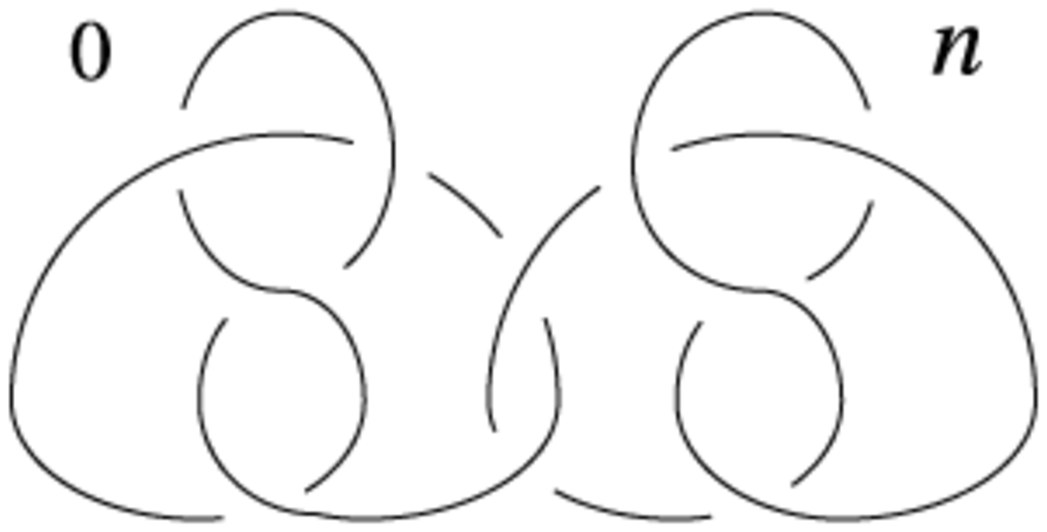}
  \end{center}
 \caption{}
   \label{fig34}
 \end{minipage}
\end{figure}

$K_1$ is a figure eight knot and $K_2$ is a right handed trefoil knot. 
\proof
We show that the $4$-manifolds represented by Figures \ref{fig33}, \ref{fig38} and \ref{fig44} have the same boundaries.

\begin{figure}[H]
 \begin{minipage}{0.1\hsize}
  \begin{center}
   \includegraphics[height=10mm]{blowup.eps}
  \end{center}
 \end{minipage}%
 \begin{minipage}{0.4\hsize}
  \begin{center}
   \includegraphics[height=20mm]{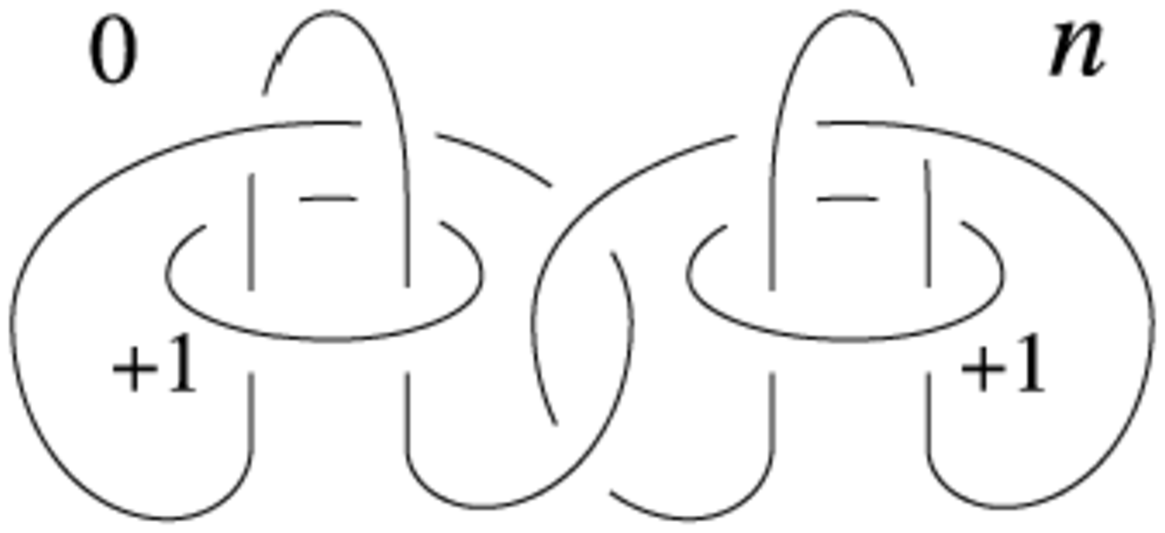}
  \end{center}
  \caption{}
   \label{fig35}
 \end{minipage}%
 \begin{minipage}{0.1\hsize}
  \begin{center}
   \includegraphics[height=10mm]{blowup.eps}
  \end{center}
 \end{minipage}%
 \begin{minipage}{0.4\hsize}
  \begin{center}
   \includegraphics[height=20mm]{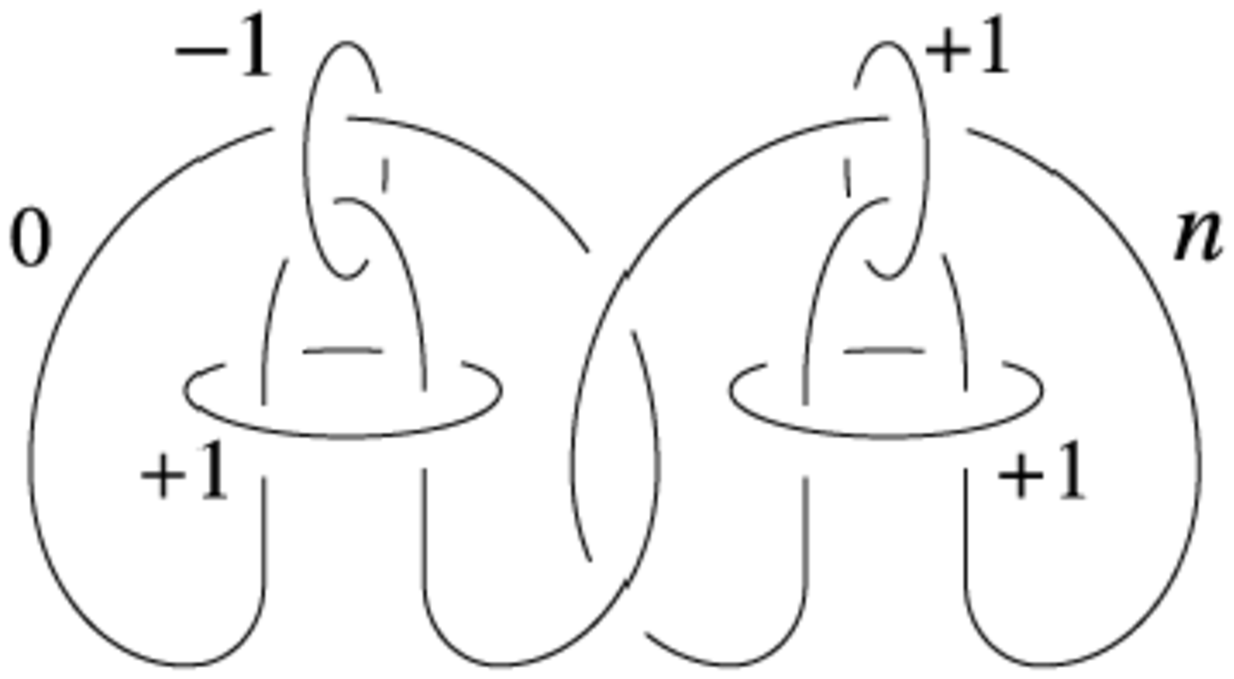}
  \end{center}
  \caption{}
   \label{fig36}
 \end{minipage}
\end{figure}

\begin{figure}[H]
 \begin{minipage}{0.1\hsize}
  \begin{center}
   \includegraphics[height=10mm]{isotopy.eps}
  \end{center}
 \end{minipage}%
 \begin{minipage}{0.3\hsize}
  \begin{center}
   \includegraphics[height=25mm]{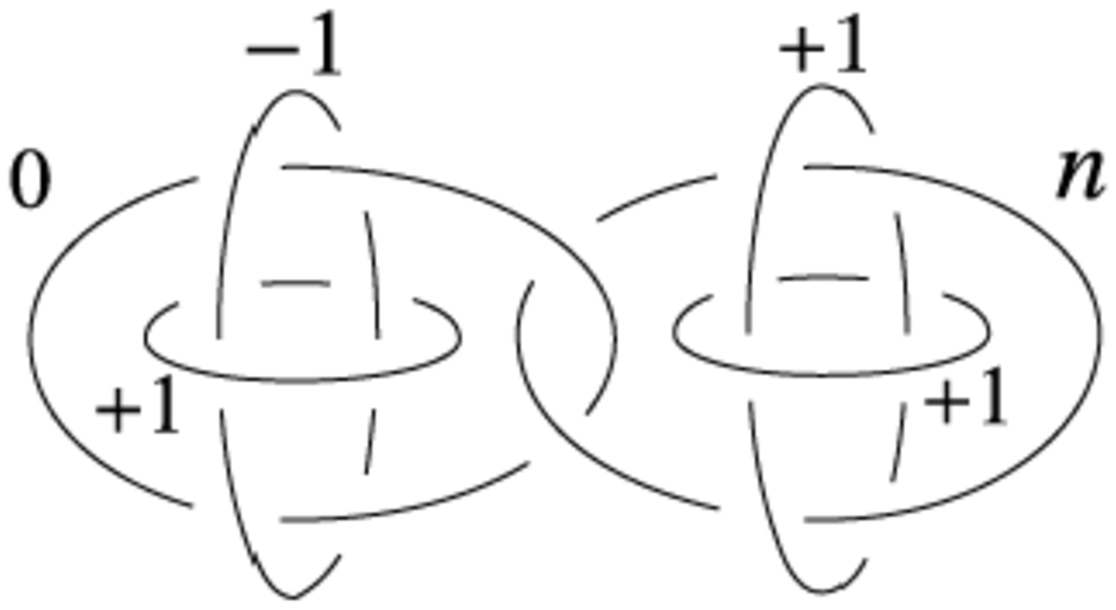}
  \end{center}
  \caption{}
   \label{fig37}
 \end{minipage}%
 \begin{minipage}{0.25\hsize}
  \begin{center}
   \includegraphics[height=25mm]{sproc.eps}
  \end{center}
 \end{minipage}%
 \begin{minipage}{0.35\hsize}
  \begin{center}
   \includegraphics[height=30mm]{Figure8.eps}
  \end{center}
  \caption{$S^3_{-1}(D_+(K_2, n))$}
   \label{fig38}
 \end{minipage}
\end{figure}

\begin{figure}[H]
 \begin{minipage}{0.45\hsize}
  \begin{center}
   \includegraphics[height=20mm]{Figure33.eps}
  \end{center}
  \caption{$M_n(K_1, K_2)$}
   \label{fig39}
 \end{minipage}%
 \begin{minipage}{0.1\hsize}
  \begin{center}
   \includegraphics[height=10mm]{isotopy.eps}
  \end{center}
 \end{minipage}%
 \begin{minipage}{0.45\hsize}
  \begin{center}
   \includegraphics[height=20mm]{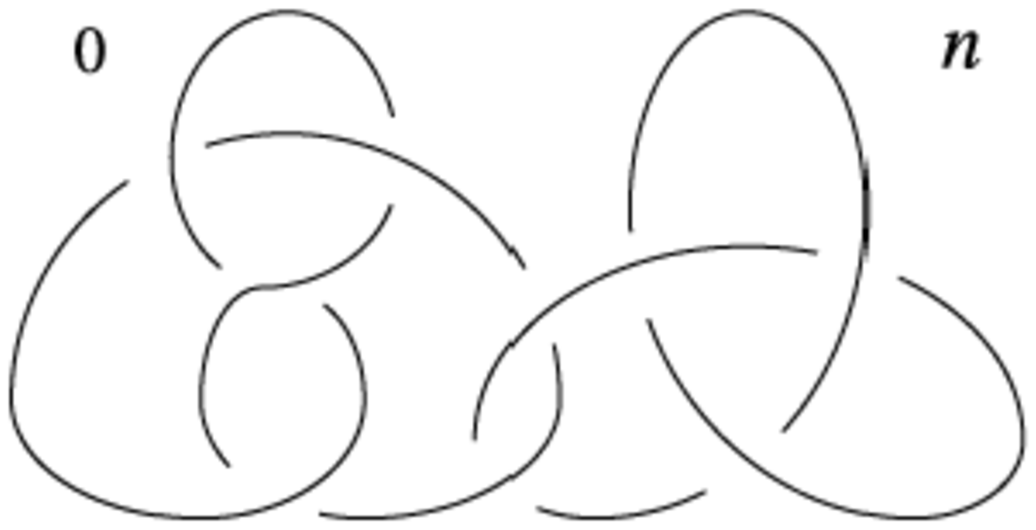}
  \end{center}
 \caption{}
   \label{fig40}
 \end{minipage}
\end{figure}

Figure \ref{fig39} is the same diagram of Figure \ref{fig33}, but by using the invertibility of the figure eight knot, we can show that they can be represented by a different doubled knot.

\begin{figure}[H]
 \begin{minipage}{0.1\hsize}
  \begin{center}
   \includegraphics[height=10mm]{isotopy.eps}
  \end{center}
 \end{minipage}%
 \begin{minipage}{0.4\hsize}
  \begin{center}
   \includegraphics[height=20mm]{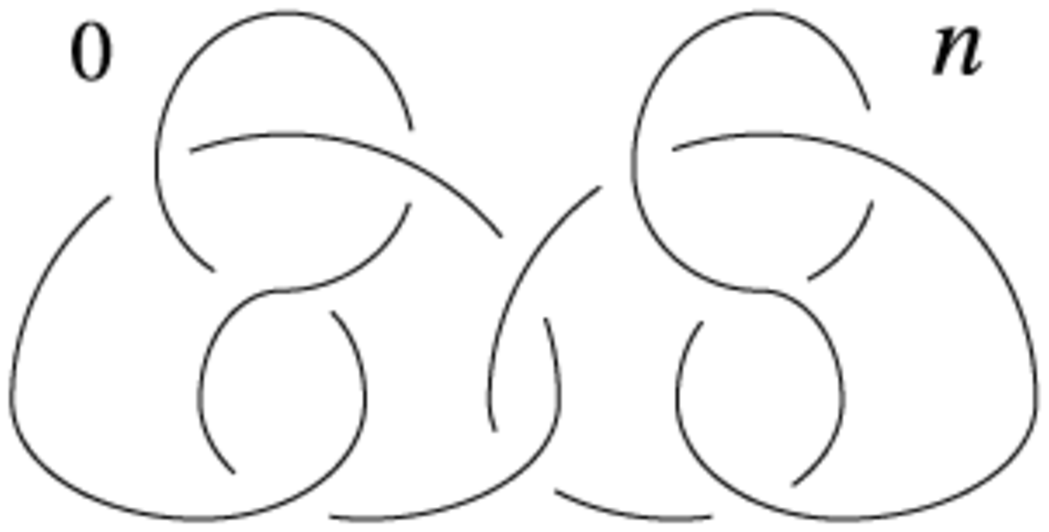}
  \end{center}
  \caption{}
   \label{fig41}
 \end{minipage}%
 \begin{minipage}{0.1\hsize}
  \begin{center}
   \includegraphics[height=10mm]{blowup.eps}
  \end{center}
 \end{minipage}%
 \begin{minipage}{0.4\hsize}
  \begin{center}
   \includegraphics[height=20mm]{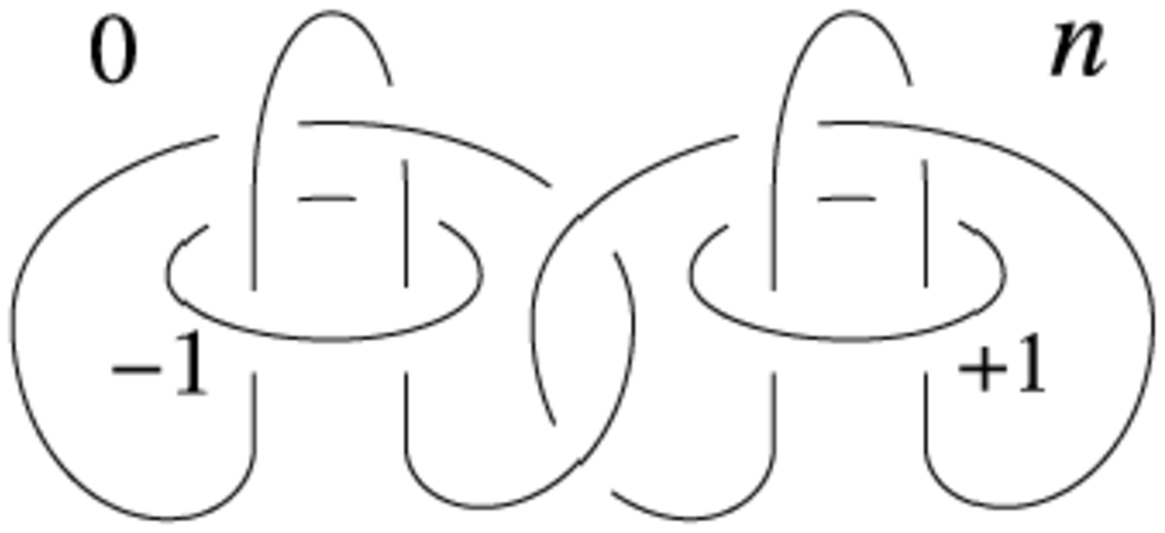}
  \end{center}
  \caption{}
   \label{fig42}
 \end{minipage}
\end{figure}

\begin{figure}[H]
 \begin{minipage}{0.1\hsize}
  \begin{center}
   \includegraphics[height=10mm]{blowup.eps}
  \end{center}
 \end{minipage}%
 \begin{minipage}{0.3\hsize}
  \begin{center}
   \includegraphics[height=25mm]{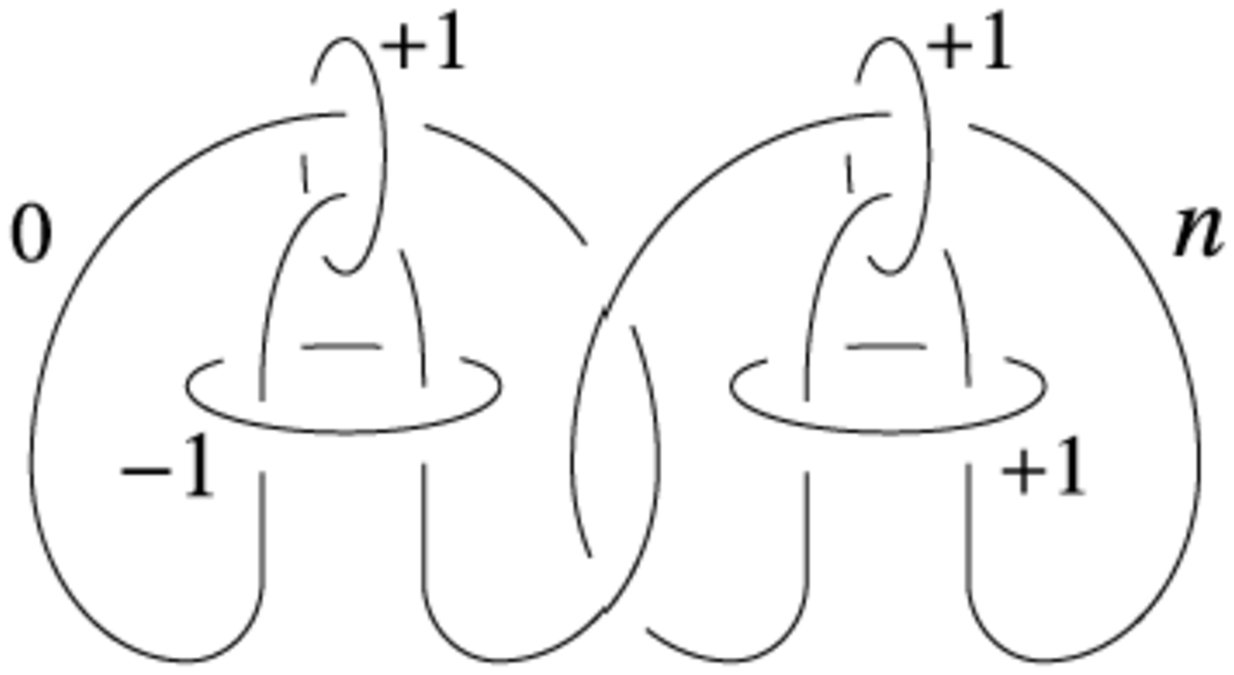}
  \end{center}
  \caption{}
   \label{fig43}
 \end{minipage}%
 \begin{minipage}{0.25\hsize}
  \begin{center}
   \includegraphics[height=25mm]{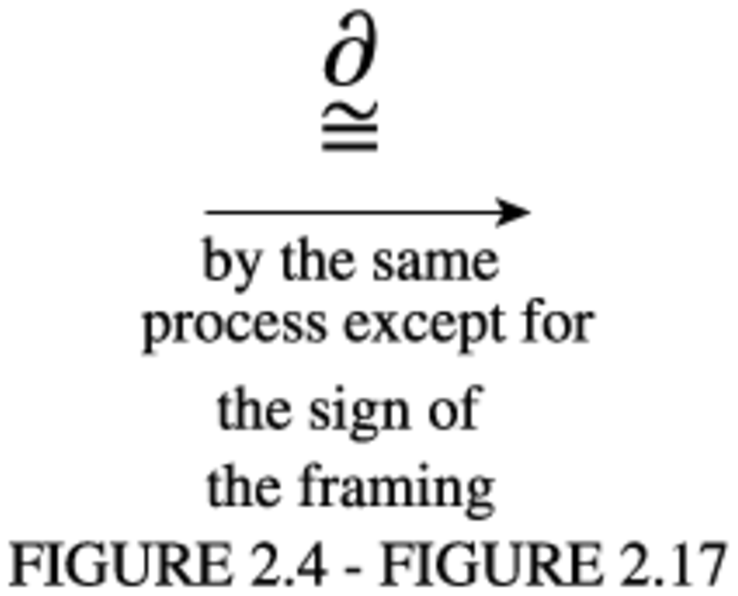}
  \end{center}
 \end{minipage}%
 \begin{minipage}{0.35\hsize}
  \begin{center}
   \includegraphics[height=30mm]{Figure9.eps}
  \end{center}
  \caption{$S^3_{+1}(D_-(K_2, n))$}
   \label{fig44}
 \end{minipage}
\end{figure}
\endproof

\subsection{Proof of the fourth row on \thmref{thm:main}'s table}
\begin{figure}[H]
 \begin{minipage}{0.45\hsize}
  \begin{center}
   \includegraphics[height=20mm]{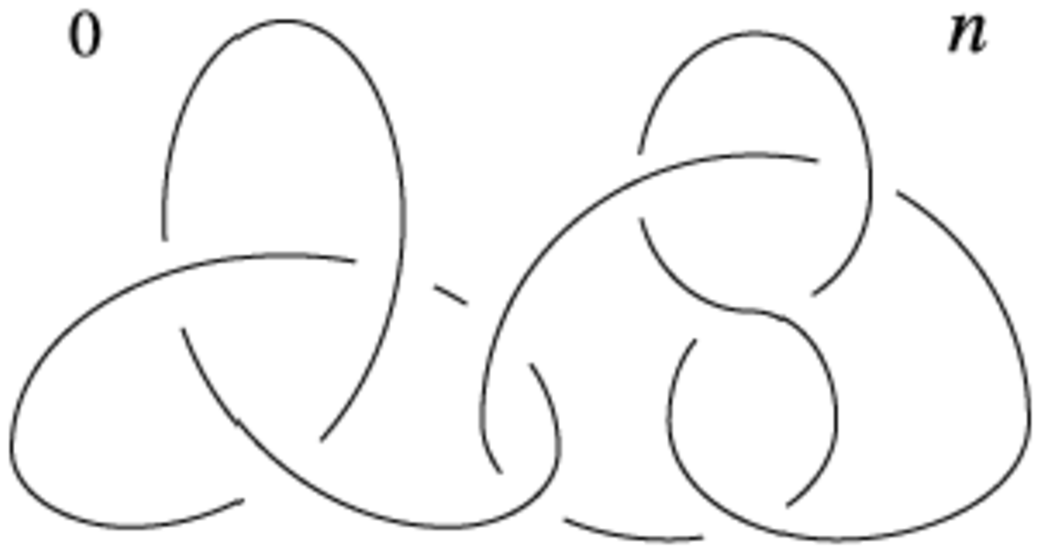}
  \end{center}
  \caption{$M_n(K_1, K_2)$}
   \label{fig45}
 \end{minipage}%
 \begin{minipage}{0.1\hsize}
  \begin{center}
   \includegraphics[height=10mm]{isotopy.eps}
  \end{center}
 \end{minipage}%
 \begin{minipage}{0.45\hsize}
  \begin{center}
   \includegraphics[height=20mm]{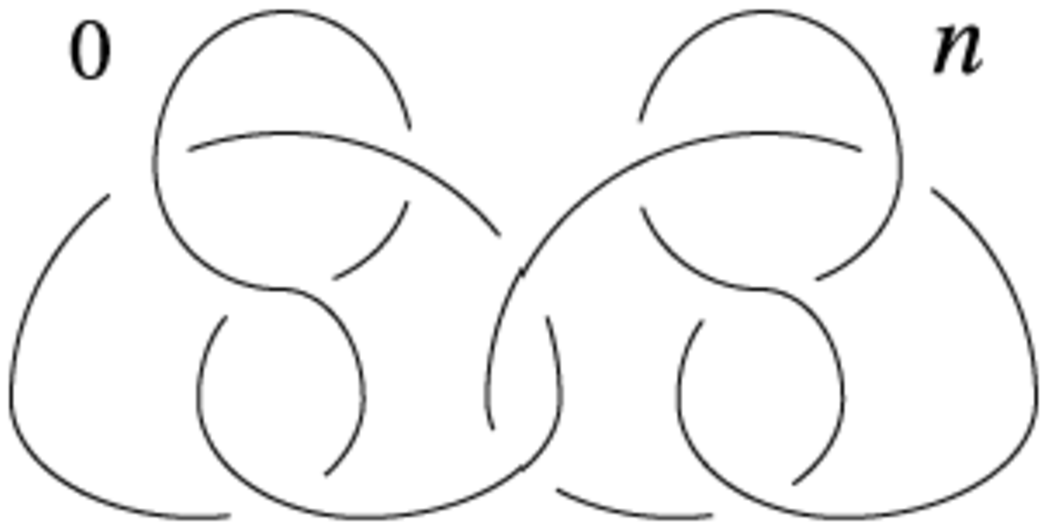}
  \end{center}
 \caption{}
   \label{fig46}
 \end{minipage}
\end{figure}

$K_1$ is a right handed trefoil knot and $K_2$ is a figure eight knot. 
\proof
We show that the $4$-manifolds represented by Figures \ref{fig45} and \ref{fig56} have the same boundaries.

\begin{figure}[H]
 \begin{minipage}{0.1\hsize}
  \begin{center}
   \includegraphics[height=10mm]{blowup.eps}
  \end{center}
 \end{minipage}%
 \begin{minipage}{0.4\hsize}
  \begin{center}
   \includegraphics[height=20mm]{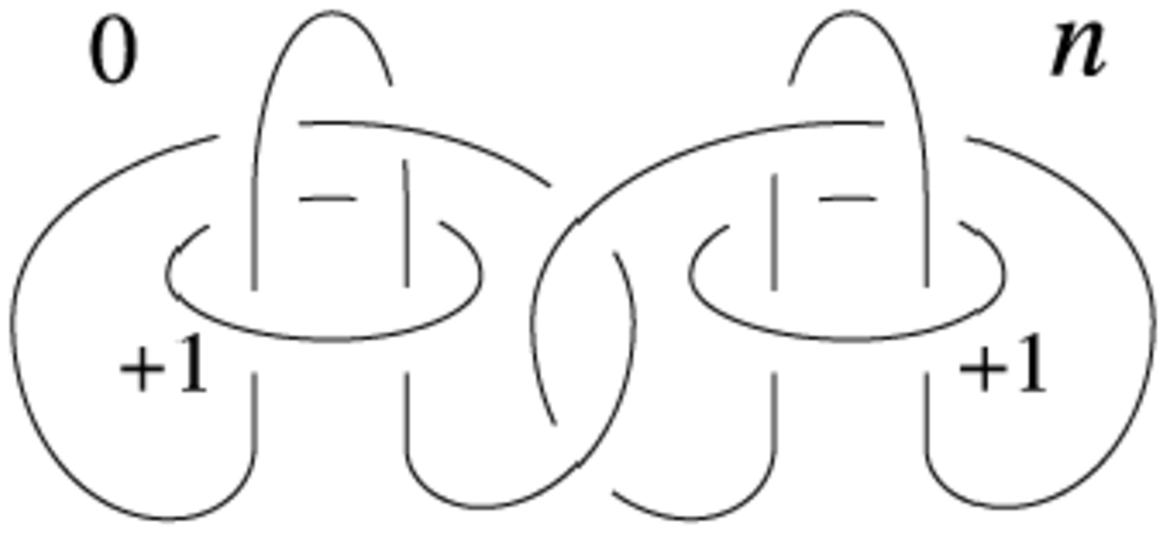}
  \end{center}
  \caption{}
   \label{fig47}
 \end{minipage}%
 \begin{minipage}{0.1\hsize}
  \begin{center}
   \includegraphics[height=10mm]{blowup.eps}
  \end{center}
 \end{minipage}%
 \begin{minipage}{0.4\hsize}
  \begin{center}
   \includegraphics[height=20mm]{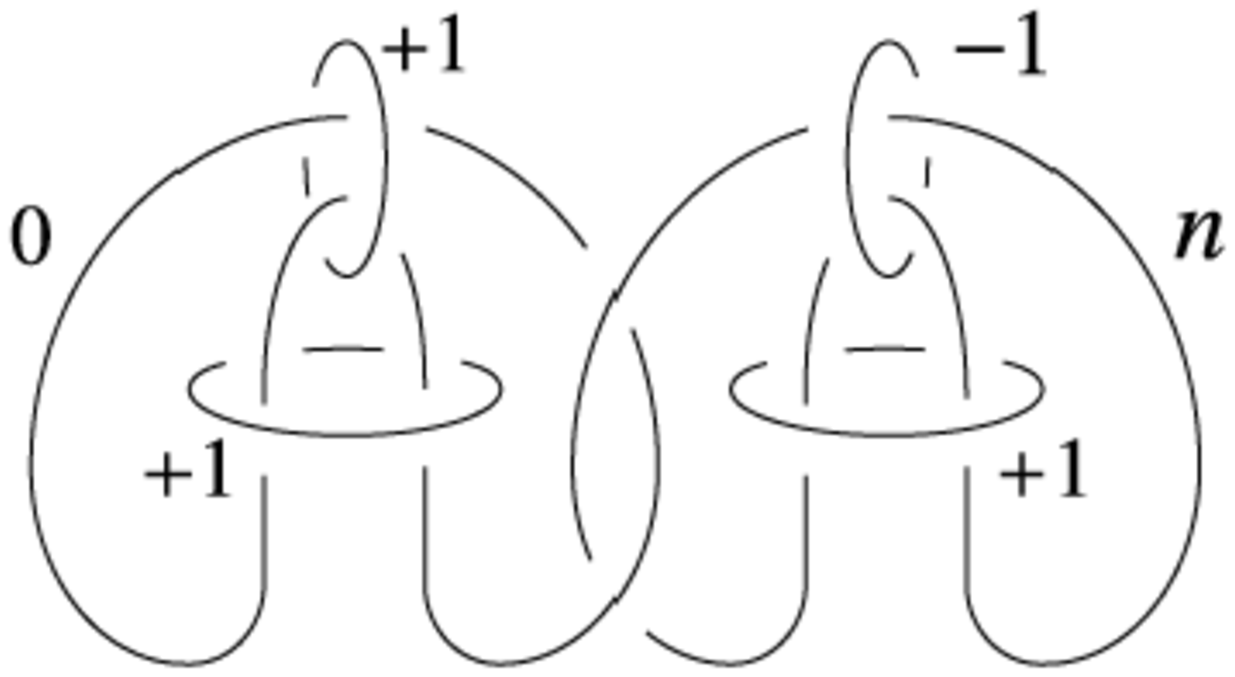}
  \end{center}
  \caption{}
   \label{fig48}
 \end{minipage}
\end{figure}

\begin{figure}[H]
 \begin{minipage}{0.1\hsize}
  \begin{center}
   \includegraphics[height=10mm]{isotopy.eps}
  \end{center}
 \end{minipage}%
 \begin{minipage}{0.35\hsize}
  \begin{center}
   \includegraphics[height=30mm]{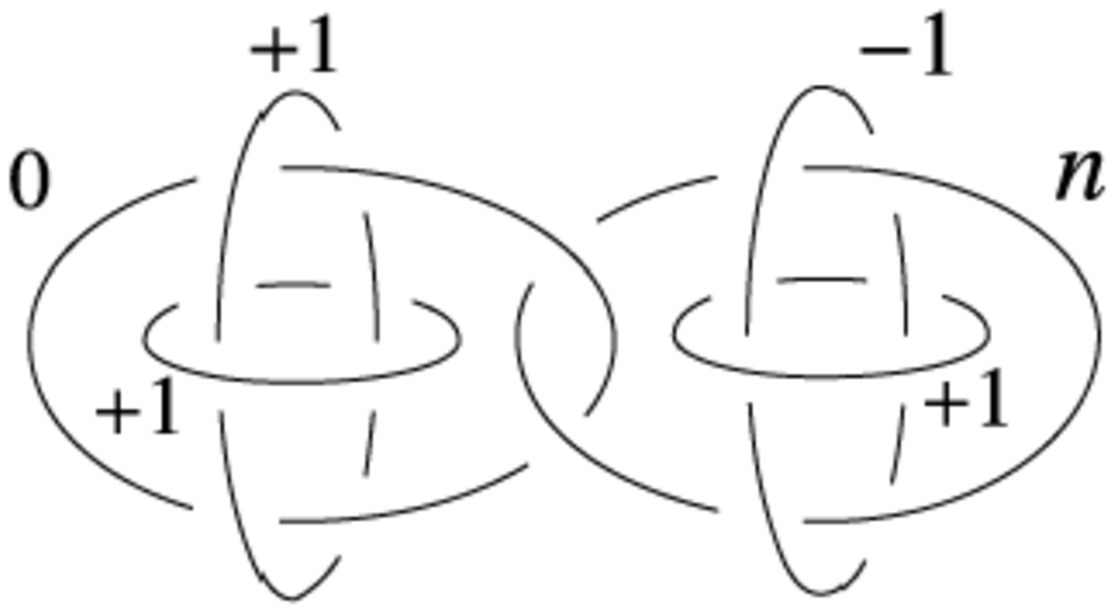}
  \end{center}
  \caption{}
   \label{fig49}
 \end{minipage}%
 \begin{minipage}{0.25\hsize}
  \begin{center}
   \includegraphics[height=30mm]{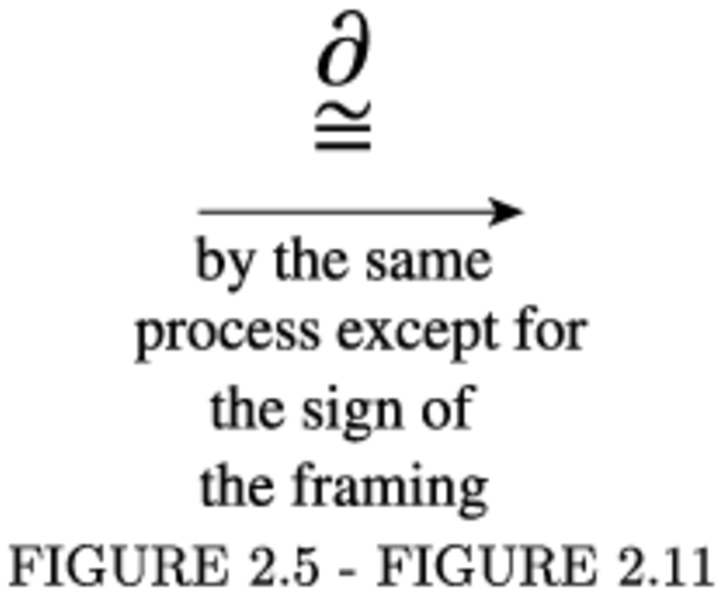}
  \end{center}
 \end{minipage}%
 \begin{minipage}{0.3\hsize}
  \begin{center}
   \includegraphics[height=30mm]{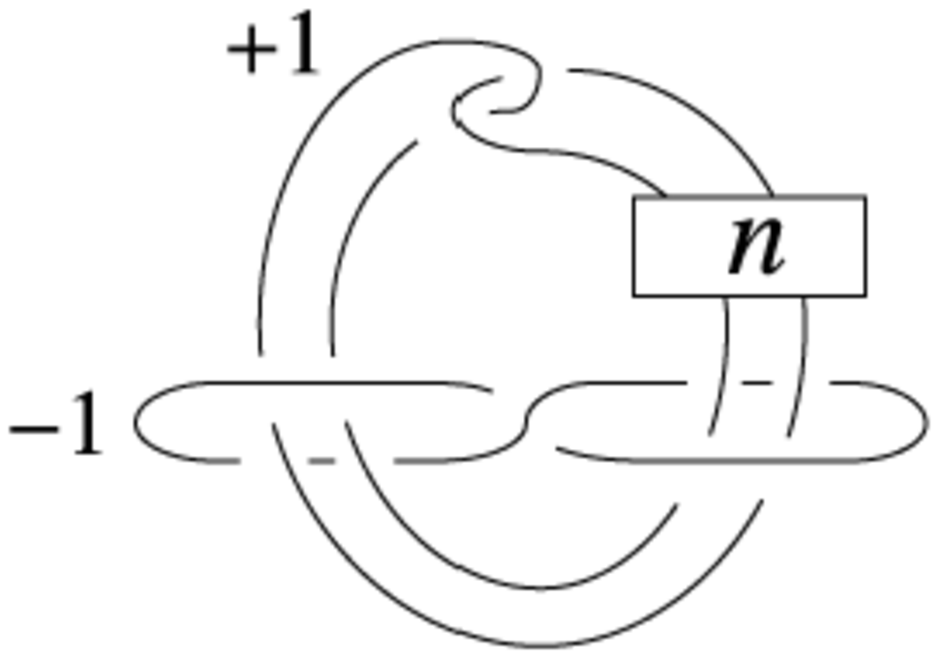}
  \end{center}
  \caption{}
   \label{fig50}
 \end{minipage}
\end{figure}

\begin{figure}[H]
 \begin{minipage}{0.1\hsize}
  \begin{center}
   \includegraphics[height=10mm]{isotopy.eps}
  \end{center}
 \end{minipage}%
 \begin{minipage}{0.4\hsize}
  \begin{center}
   \includegraphics[height=30mm]{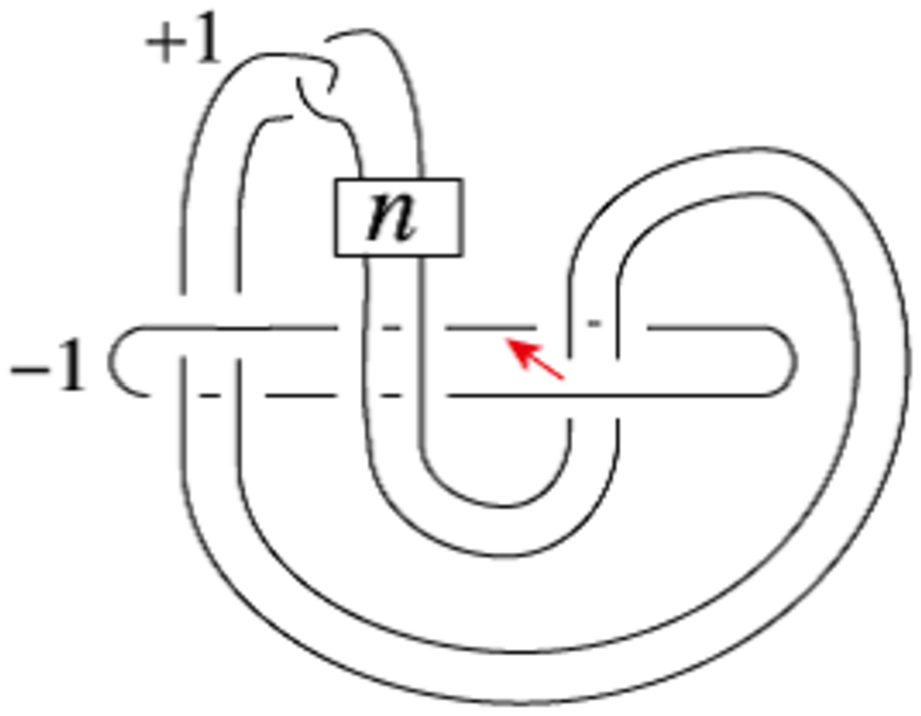}
  \end{center}
  \caption{}
   \label{fig51}
 \end{minipage}%
 \begin{minipage}{0.1\hsize}
  \begin{center}
   \includegraphics[height=10mm]{hds.eps}
  \end{center}
 \end{minipage}%
 \begin{minipage}{0.4\hsize}
  \begin{center}
   \includegraphics[height=30mm]{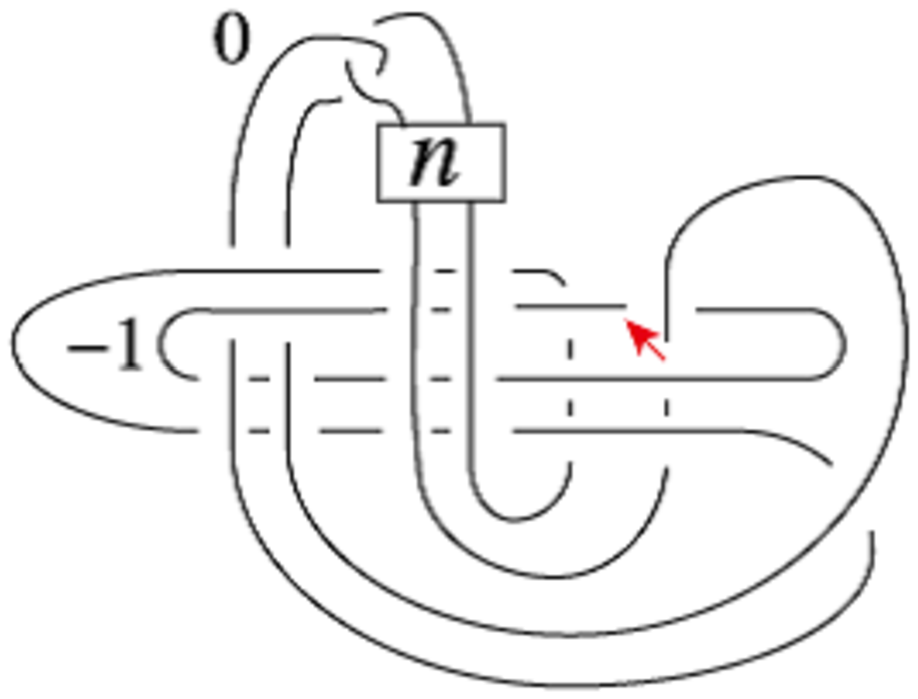}
  \end{center}
  \caption{}
   \label{fig52}
 \end{minipage}
\end{figure}

\begin{figure}[H]
 \begin{minipage}{0.1\hsize}
  \begin{center}
   \includegraphics[height=10mm]{hds.eps}
  \end{center}
 \end{minipage}%
 \begin{minipage}{0.4\hsize}
  \begin{center}
   \includegraphics[height=30mm]{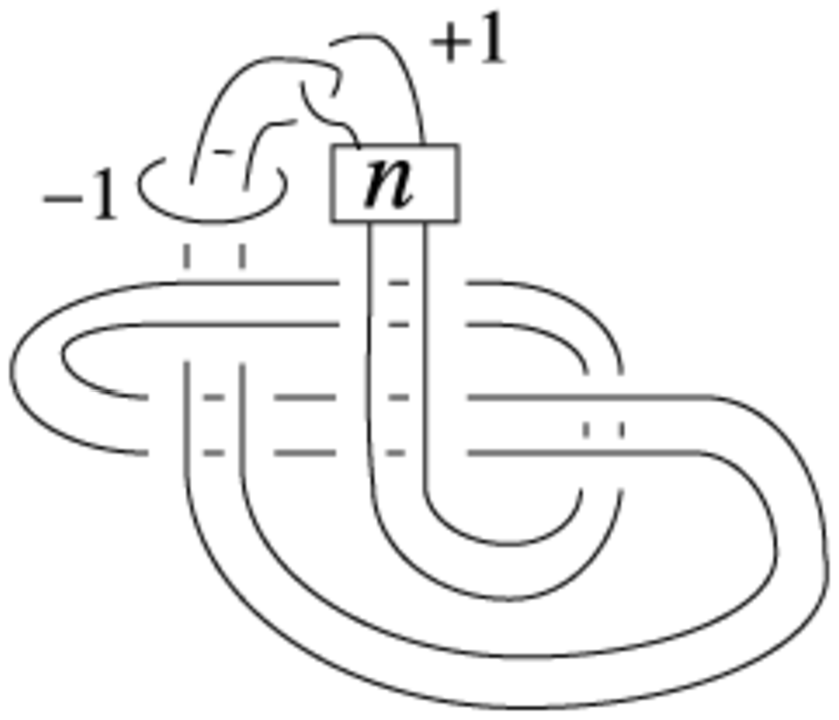}
  \end{center}
  \caption{}
   \label{fig53}
 \end{minipage}%
 \begin{minipage}{0.1\hsize}
  \begin{center}
   \includegraphics[height=10mm]{blowdown.eps}
  \end{center}
 \end{minipage}%
 \begin{minipage}{0.4\hsize}
  \begin{center}
   \includegraphics[height=30mm]{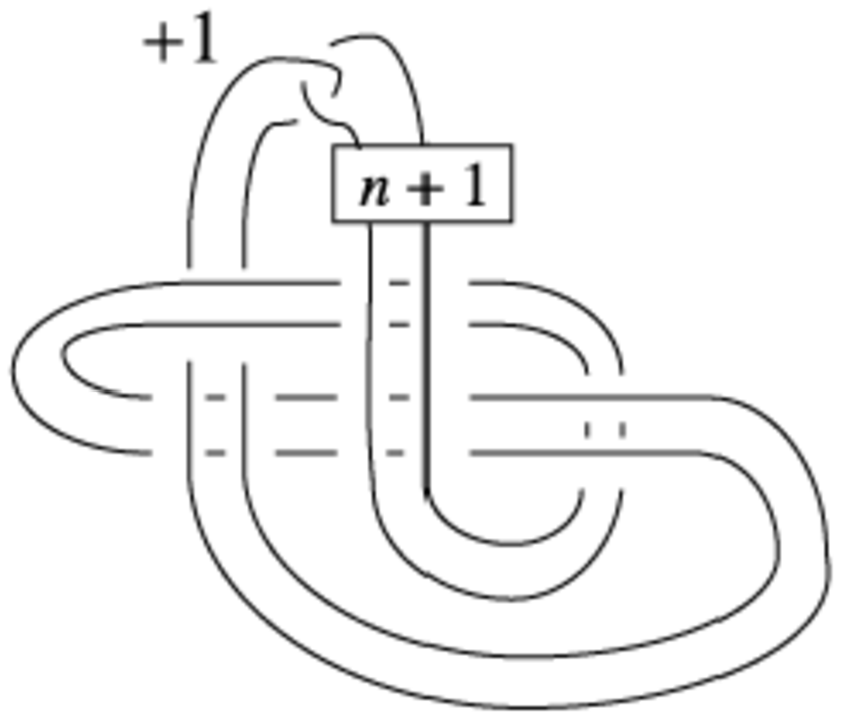}
  \end{center}
  \caption{}
   \label{fig54}
 \end{minipage}
\end{figure}

\begin{figure}[H]
 \begin{minipage}{0.1\hsize}
  \begin{center}
   \includegraphics[height=10mm]{isotopy.eps}
  \end{center}
 \end{minipage}%
 \begin{minipage}{0.4\hsize}
  \begin{center}
   \includegraphics[height=30mm]{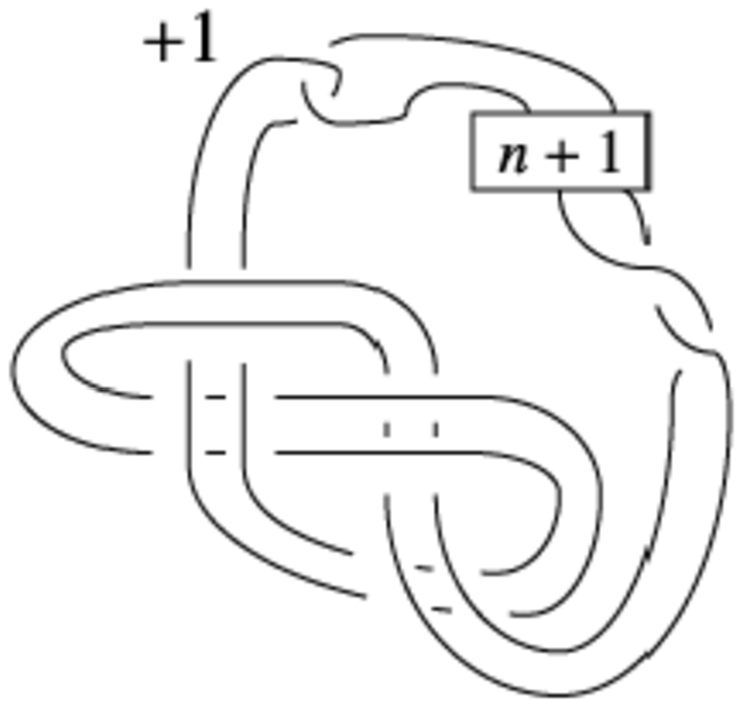}
  \end{center}
  \caption{}
   \label{fig55}
 \end{minipage}%
 \begin{minipage}{0.1\hsize}
  \begin{center}
   \includegraphics[height=10mm]{isotopy.eps}
  \end{center}
 \end{minipage}%
 \begin{minipage}{0.4\hsize}
  \begin{center}
   \includegraphics[height=30mm]{Figure4.eps}
  \end{center}
  \caption{$S^3_{+1}(D_+(K_2, n))$}
   \label{fig56}
 \end{minipage}
\end{figure}
\endproof

\subsection{Proof of the fifth row on \thmref{thm:main}'s table}
\begin{figure}[H]
 \begin{minipage}{0.45\hsize}
  \begin{center}
   \includegraphics[height=20mm]{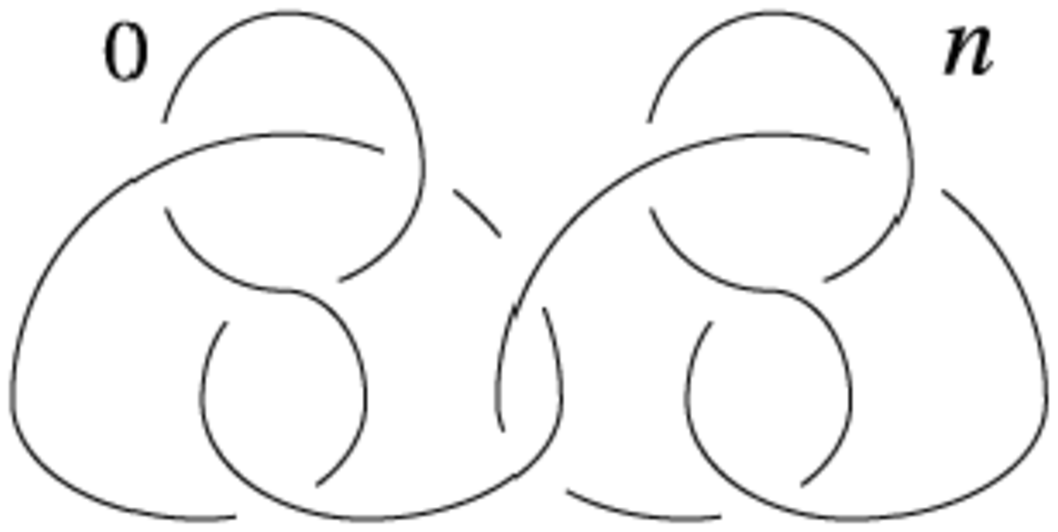}
  \end{center}
  \caption{$M_n(K_1, K_2)$}
   \label{fig57}
 \end{minipage}%
 \begin{minipage}{0.1\hsize}
  \begin{center}
   \includegraphics[height=10mm]{blowup.eps}
  \end{center}
 \end{minipage}%
 \begin{minipage}{0.45\hsize}
  \begin{center}
   \includegraphics[height=20mm]{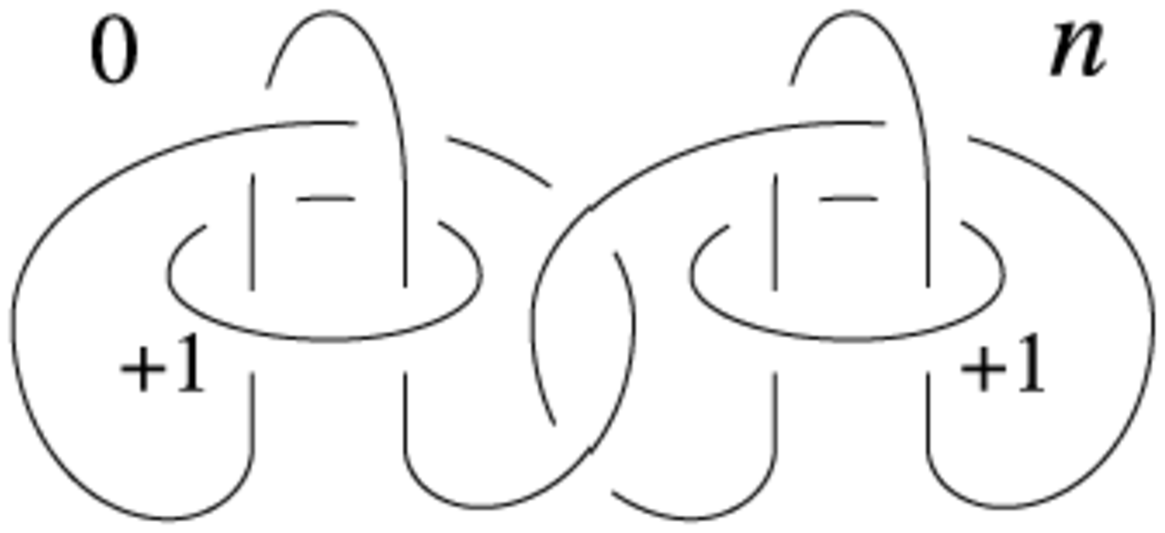}
  \end{center}
 \caption{}
   \label{fig58}
 \end{minipage}
\end{figure}

$K_1$ and $K_2$ are figure eight knots. 
\proof
We show that the $4$-manifolds represented by Figures \ref{fig57} and \ref{fig61-2} have the same boundaries.

\begin{figure}[H]
 \begin{minipage}{0.1\hsize}
  \begin{center}
   \includegraphics[height=10mm]{blowup.eps}
  \end{center}
 \end{minipage}%
 \begin{minipage}{0.4\hsize}
  \begin{center}
   \includegraphics[height=30mm]{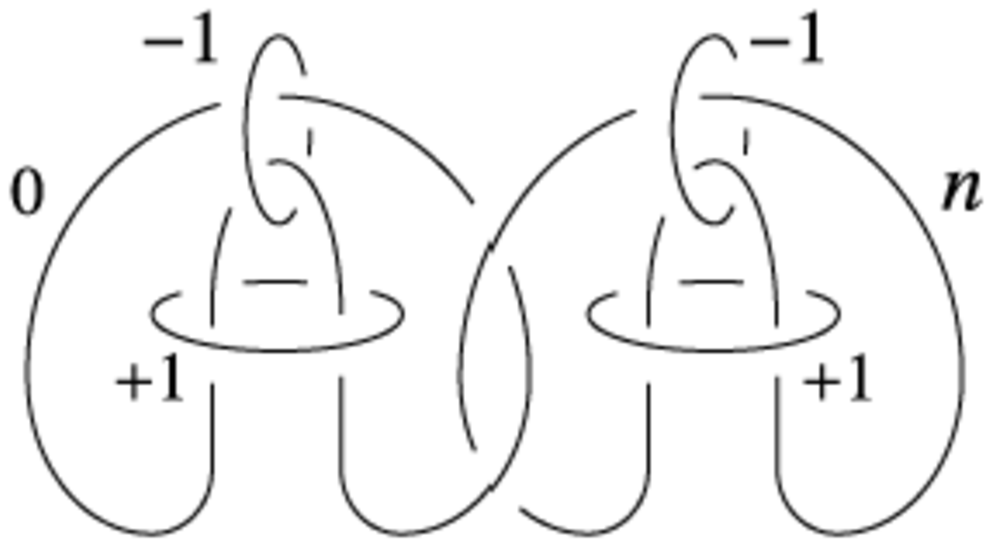}
  \end{center}
  \caption{}
   \label{fig59}
 \end{minipage}%
 \begin{minipage}{0.1\hsize}
  \begin{center}
   \includegraphics[height=10mm]{isotopy.eps}
  \end{center}
 \end{minipage}%
 \begin{minipage}{0.4\hsize}
  \begin{center}
   \includegraphics[height=30mm]{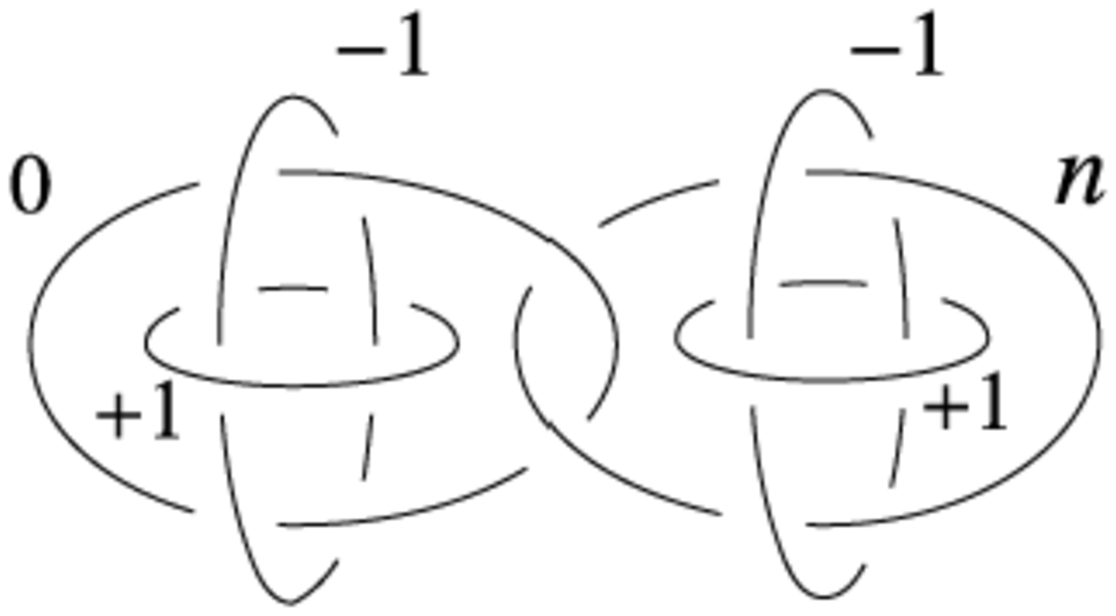}
  \end{center}
  \caption{}
   \label{fig60}
 \end{minipage}
\end{figure}

\begin{figure}[H]
 \begin{minipage}{0.5\hsize}
  \begin{center}
   \includegraphics[height=30mm]{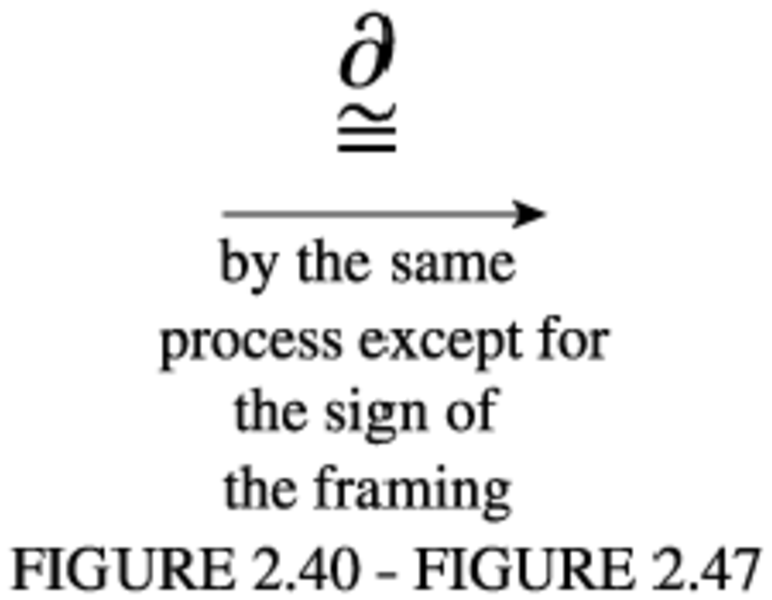}
  \end{center}
 \end{minipage}%
 \begin{minipage}{0.5\hsize}
  \begin{center}
   \includegraphics[height=30mm]{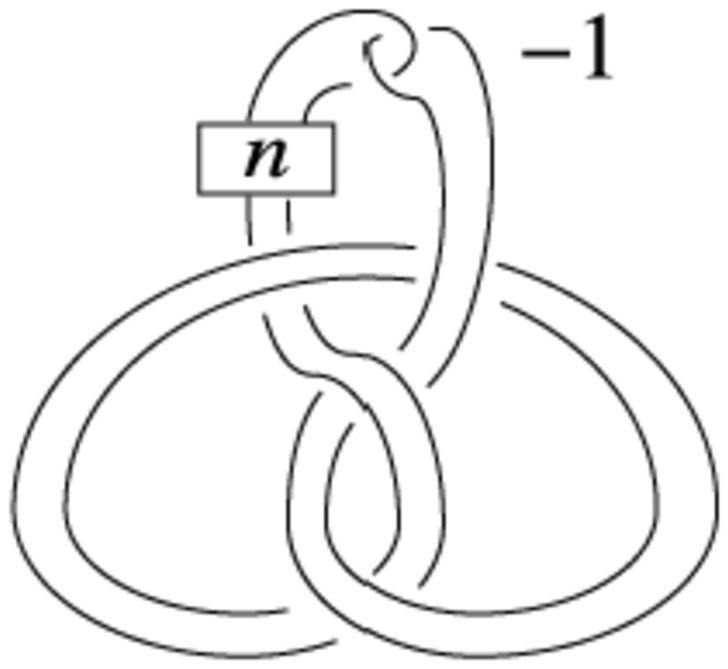}
  \end{center}
  \caption{$S^3_{-1}(D_+(K_2, n))$}
  \label{fig61}
 \end{minipage}
\end{figure}

\begin{figure}[H]
 \begin{minipage}{0.45\hsize}
  \begin{center}
   \includegraphics[height=30mm]{Figure60.eps}
  \end{center}
  \caption{}
   \label{fig60-1}
 \end{minipage}%
 \begin{minipage}{0.1\hsize}
  \begin{center}
   \includegraphics[height=10mm]{isotopy.eps}
  \end{center}
 \end{minipage}%
 \begin{minipage}{0.45\hsize}
  \begin{center}
   \includegraphics[height=30mm]{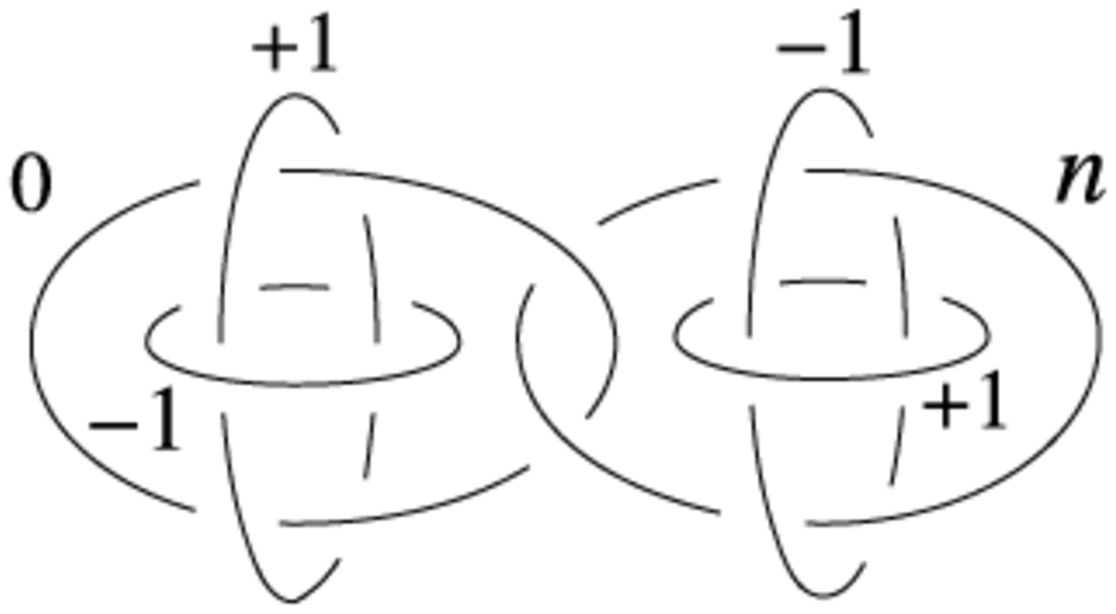}
  \end{center}
 \caption{}
   \label{fig60-2}
 \end{minipage}
\end{figure}

\begin{figure}[H]
 \begin{minipage}{0.5\hsize}
  \begin{center}
   \includegraphics[height=30mm]{sproc3.eps}
  \end{center}
 \end{minipage}%
 \begin{minipage}{0.5\hsize}
  \begin{center}
   \includegraphics[height=30mm]{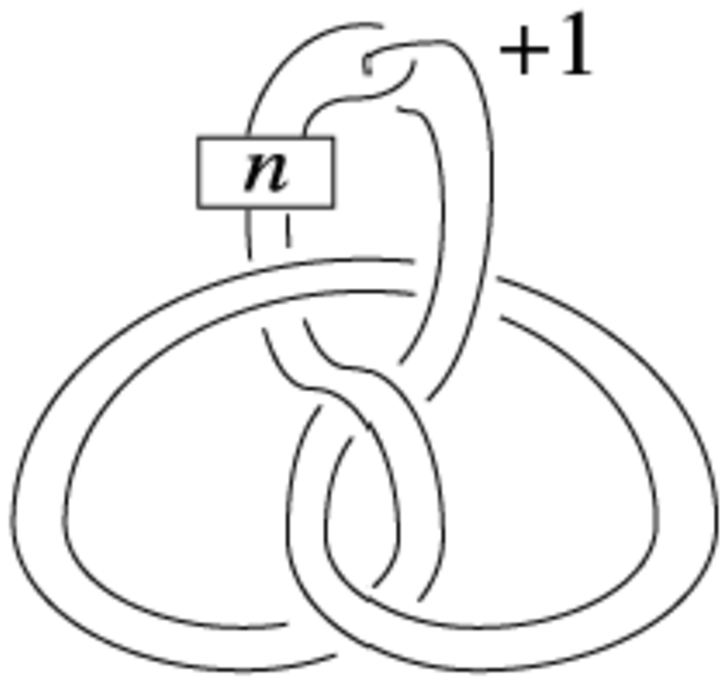}
  \end{center}
  \caption{$S^3_{+1}(D_-(K_2, n))$}
  \label{fig61-2}
 \end{minipage}
\end{figure}

Since Figure \ref{fig60-1} is the same diagram of Figure \ref{fig60}, we can show that they can be represented by a different double knot.

\endproof

Next we compute the Casson invariant $\lambda(M_n(K_1,K_2))$. Now suppose that $K_-$, $K_+$ and $K_0$ are links in $S^3$ which have projections which differ at a single crossing of $K_-$ as depicted below.

\begin{figure}[H]
 \begin{center}
  \includegraphics[height=25mm]{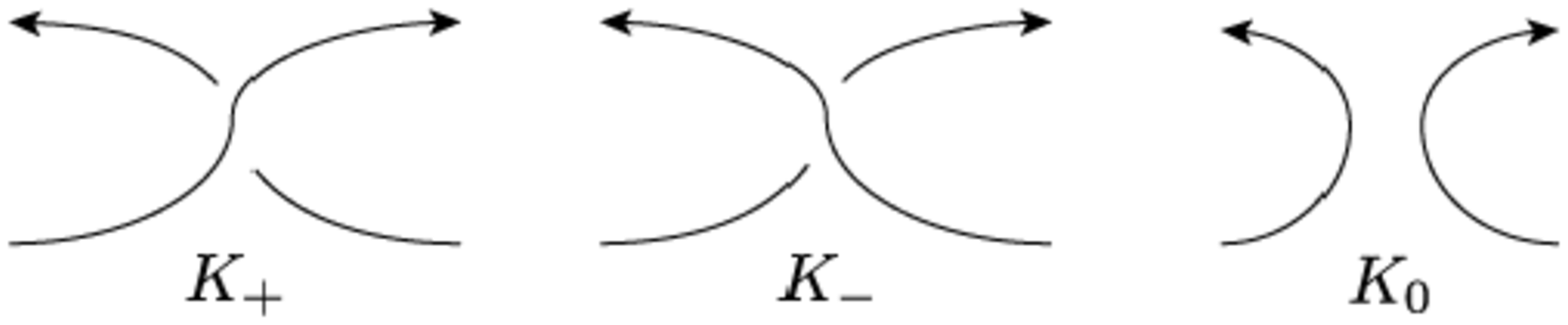}
 \end{center}
 \caption{}
 \label{fig62}
\end{figure}

\noindent\textbf{Remark.}
Our convention in Figure \ref{fig62} is different from that in \cite{AM}. In fact, their $K_+$ (resp.$K_-$) is our $K_-$ (resp.$K_+$). We adopt our convention as Figure \ref{fig62} because by our convention $\lambda'(T_{2,3})$ is computed to be 1, where $T_{2,3}$ is a right handed trefoil knot. While by their convention $\lambda'(T_{2,3})$ is computed to be $-1$, contradicting the normalization  $\lambda'(T_{2,3}) = 1$ (\cite{AM} p148, \cite{S} p52).

\begin{Lemma}[see \cite{AM}, p. 143]\label{lem2}
Let $K_-$ be a knot in $S^3$. Let $K_+$ and $K_0$ be as above. Then $K_0$ is a two component link and:
\begin{center}
$\lambda'(K_+) - \lambda'(K_-) = lk(K_0)$
\end{center}
where $\lambda'(K)$ is the Casson invariant of a knot K.\\
\end{Lemma}

\begin{Lemma}[Surgery formula, see \cite{S}, p. 52]\label{lem3}
Let K be a knot in $S^3$. The Casson invariant $\lambda(S^3_{+1}(K))$ is equal to $\lambda'(K)$.
\end{Lemma}

By Lemmas \ref{lem2}, \ref{lem3} and the second column on \thmref{thm:main}'s table, we can compute the Casson invariant $\lambda(M_n(K_1, K_2))$.

\subsection{The Casson invariant of the first row on \thmref{thm:main}'s table}
\proof
$K_1$ and $K_2$ are right handed trefoil knots. By 2(i), $M_n(K_1, K_2)$ is diffeomorphic to $S^3_{+1}(D_+(K_2, n))$. Therefore $\lambda(M_n(K_1, K_2))$ is equal to $\lambda(S^3_{+1}(D_+(K_2, n)))$. By Lemmas \ref{lem2} and \ref{lem3}, we can compute the Casson invariant $\lambda(S^3_{+1}(D_+(K_2, n)))$ as follows:

\begin{figure}[H]
 \begin{minipage}{0.33\hsize}
  \begin{center}
   \includegraphics[height=30mm]{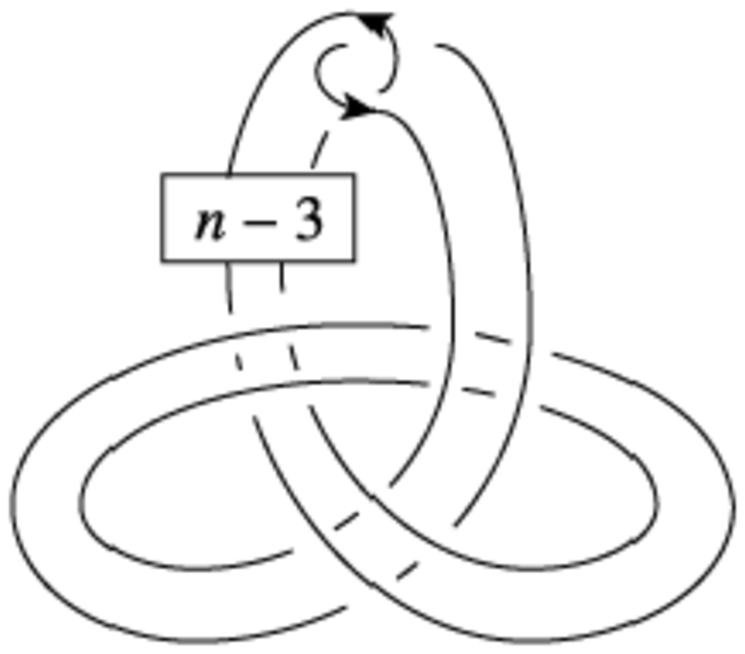}
  \end{center}
  \caption{$K_+$}
   \label{fig63}
 \end{minipage}%
 \begin{minipage}{0.33\hsize}
  \begin{center}
   \includegraphics[height=30mm]{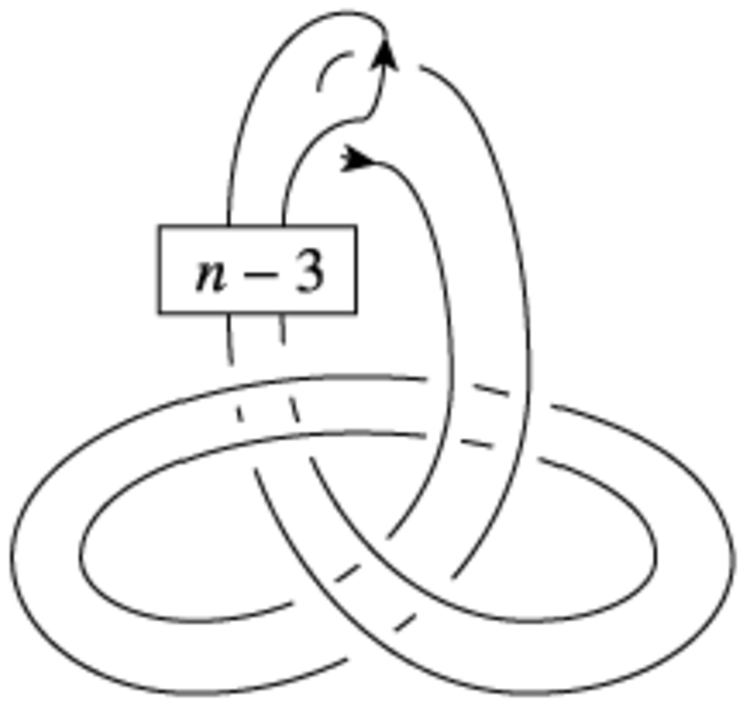}
  \end{center}
 \caption{$K_-$}
   \label{fig64}
 \end{minipage}%
 \begin{minipage}{0.33\hsize}
  \begin{center}
   \includegraphics[height=30mm]{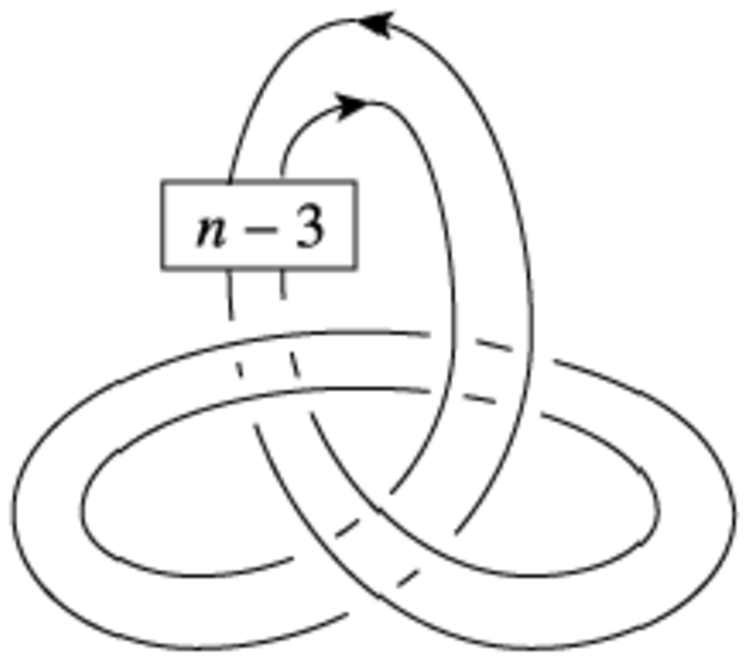}
  \end{center}
 \caption{$K_0$}
   \label{fig65}
 \end{minipage}
\end{figure}

By Lemma \ref{lem2}, $\lambda'(K_+) - \lambda'(K_-) = lk(K_0)$. Since $K_-$ is a trivial knot, $\lambda'(K_-) = 0$. $lk(K_0)$ is $-n$. Therefore, $\lambda'(K_+) = -n$. Then $\lambda(M_n(K_1, K_2)) = \lambda(S^3_{+1}(D_+(K_2, n))) = \lambda'(K_+) = -n$.

\endproof

\subsection{The Casson invariant of the second row on \thmref{thm:main}'s table}
\proof
$K_1$ is a left handed trefoil knot and $K_2$ is a right handed trefoil knot. By 2(ii), $M_n(K_1, K_2)$ is diffeomorphic to $S^3_{-1}(D_-(K_2, n))$. Therefore $\lambda(M_n(K_1, K_2))$ is equal to $\lambda(S^3_{-1}(D_-(K_2, n)))$. Since $\lambda(S^3_{-1}(D_-(K_2, n))) = -\lambda(S^3_{+1}(D_+(K_1, -n)))$ (see \cite{S}, p. 52, Theorem 3.1.), we compute the Casson invariant $\lambda(S^3_{+1}(D_+(K_1, -n)))$ as follows:

\begin{figure}[H]
 \begin{minipage}{0.33\hsize}
  \begin{center}
   \includegraphics[height=30mm]{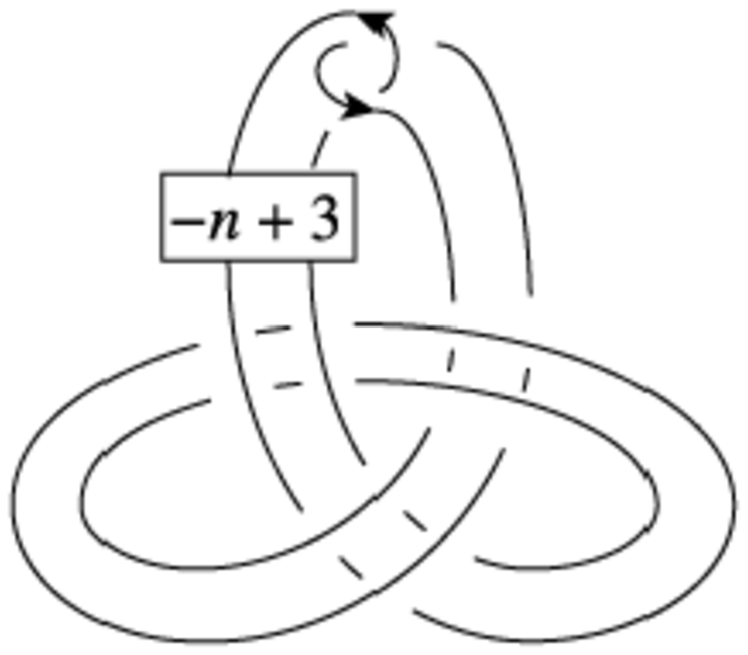}
  \end{center}
  \caption{$K_+$}
   \label{fig66}
 \end{minipage}%
 \begin{minipage}{0.33\hsize}
  \begin{center}
   \includegraphics[height=30mm]{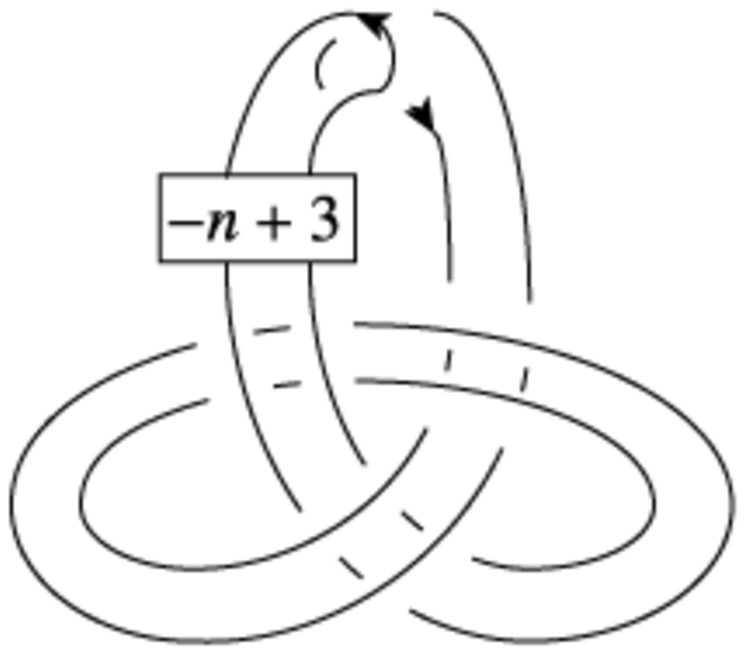}
  \end{center}
 \caption{$K_-$}
   \label{fig67}
 \end{minipage}%
 \begin{minipage}{0.33\hsize}
  \begin{center}
   \includegraphics[height=30mm]{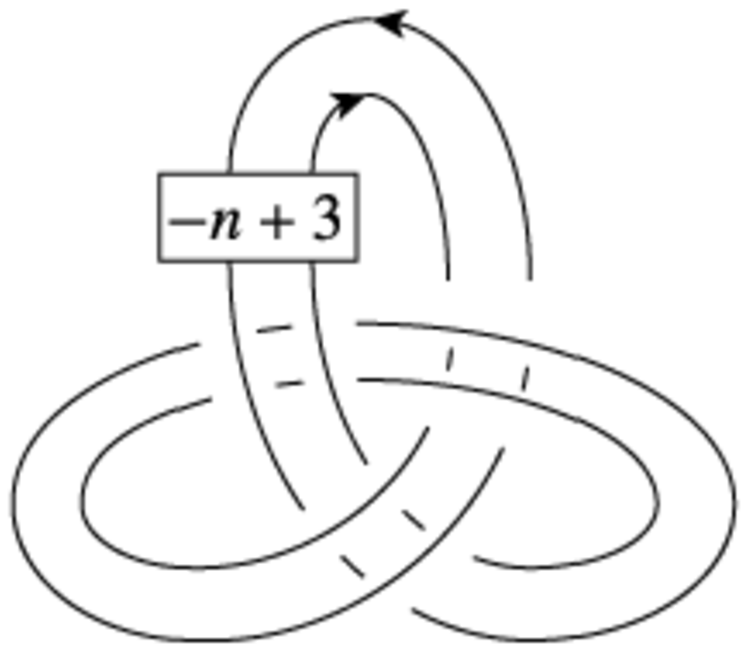}
  \end{center}
 \caption{$K_0$}
   \label{fig68}
 \end{minipage}
\end{figure}

By Lemma \ref{lem2}, $\lambda'(K_+) - \lambda'(K_-) = lk(K_0)$. Since $K_-$ is a trivial knot, $\lambda'(K_-) = 0$. $lk(K_0)$ is $n$. Therefore, $\lambda'(K_+) = n$. Then $\lambda(M_n(K_1, K_2)) = \lambda(S^3_{-1}(D_-(K_2, n))) = -\lambda(S^3_{+1}(D_+(K_1, -n))) = -\lambda'(K_+) = -n$.
\endproof

\subsection{The Casson invariant of the third row on \thmref{thm:main}'s table}
\proof
$K_1$ is a figure eight  knot and $K_2$ is a right handed trefoil knot. By 2(iii), $M_n(K_1, K_2)$ is diffeomorphic to $S^3_{+1}(D_-(K_2, n))$. Therefore $\lambda(M_n(K_1, K_2))$ is equal to $\lambda(S^3_{+1}(D_-(K_2, n)))$. By Lemmas \ref{lem2} and \ref{lem3}, we can compute the Casson invariant $\lambda(S^3_{+1}(D_-(K_2, n)))$ as follows:

\begin{figure}[H]
 \begin{minipage}{0.33\hsize}
  \begin{center}
   \includegraphics[height=30mm]{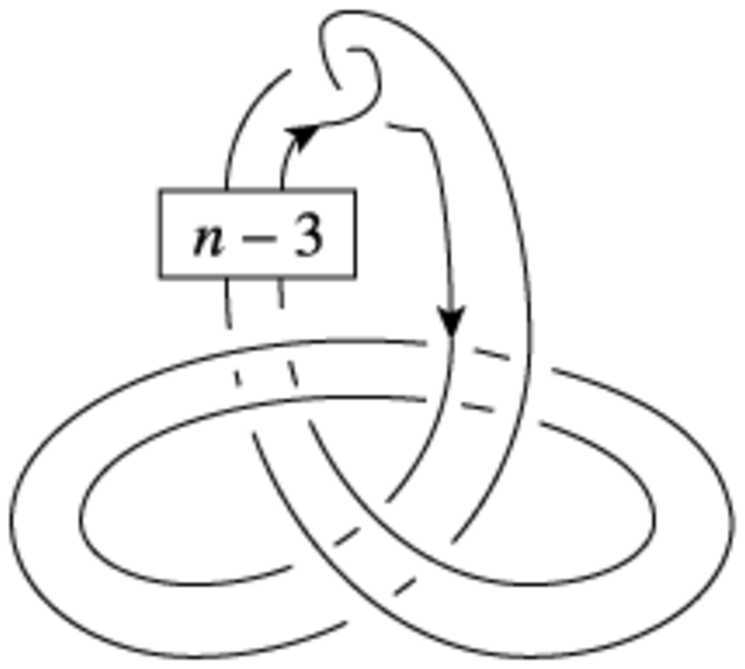}
  \end{center}
  \caption{$K_-$}
   \label{fig69}
 \end{minipage}%
 \begin{minipage}{0.33\hsize}
  \begin{center}
   \includegraphics[height=30mm]{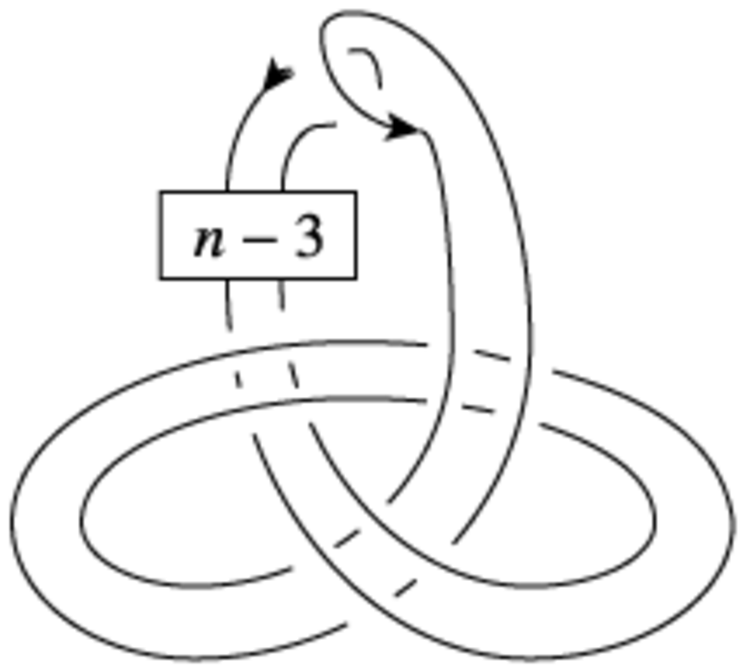}
  \end{center}
 \caption{$K_+$}
   \label{fig70}
 \end{minipage}%
 \begin{minipage}{0.33\hsize}
  \begin{center}
   \includegraphics[height=30mm]{Figure65.eps}
  \end{center}
 \caption{$K_0$}
   \label{fig70-2}
 \end{minipage}
\end{figure}

By Lemma \ref{lem2}, $\lambda'(K_+) - \lambda'(K_-) = lk(K_0)$. Since $K_+$ is a trivial knot, $\lambda'(K_+) = 0$. $lk(K_0)$ is $-n$. Therefore, $\lambda'(K_-) = n$. Then $\lambda(M_n(K_1, K_2)) = \lambda(S^3_{+1}(D_-(K_2, n))) = \lambda'(K_-) = n$.

\endproof

\subsection{The Casson invariant of the fourth row on \thmref{thm:main}'s table}
\proof
$K_1$ is a right handed trefoil knot and $K_2$ is a figure eight knot. By 2(iv), $M_n(K_1, K_2)$ is diffeomorphic to $S^3_{+1}(D_+(K_2, n))$. Therefore $\lambda(M_n(K_1, K_2))$ is equal to $\lambda(S^3_{+1}(D_+(K_2, n)))$. By Lemmas \ref{lem2} and \ref{lem3}, we can compute the Casson invariant $\lambda(S^3_{+1}(D_+(K_2, n)))$ as follows:

\begin{figure}[H]
 \begin{minipage}{0.33\hsize}
  \begin{center}
   \includegraphics[height=30mm]{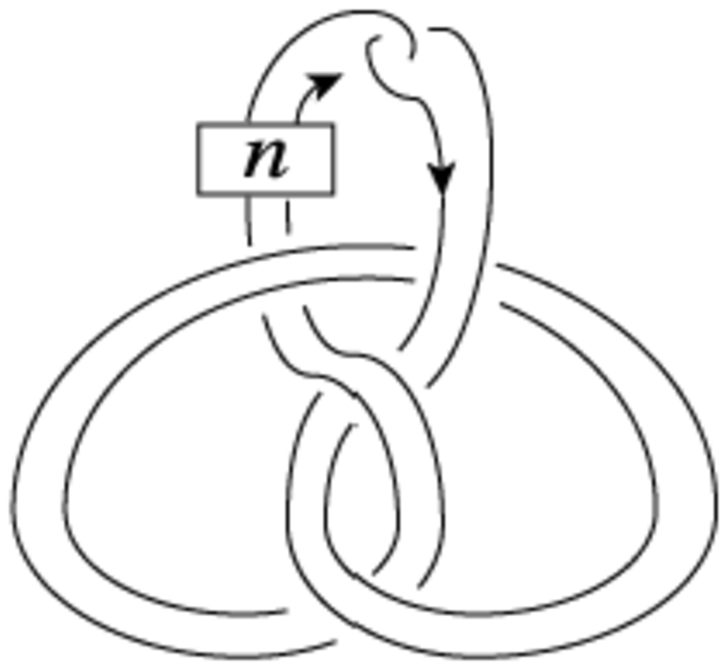}
  \end{center}
  \caption{$K_+$}
   \label{fig71}
 \end{minipage}%
 \begin{minipage}{0.33\hsize}
  \begin{center}
   \includegraphics[height=30mm]{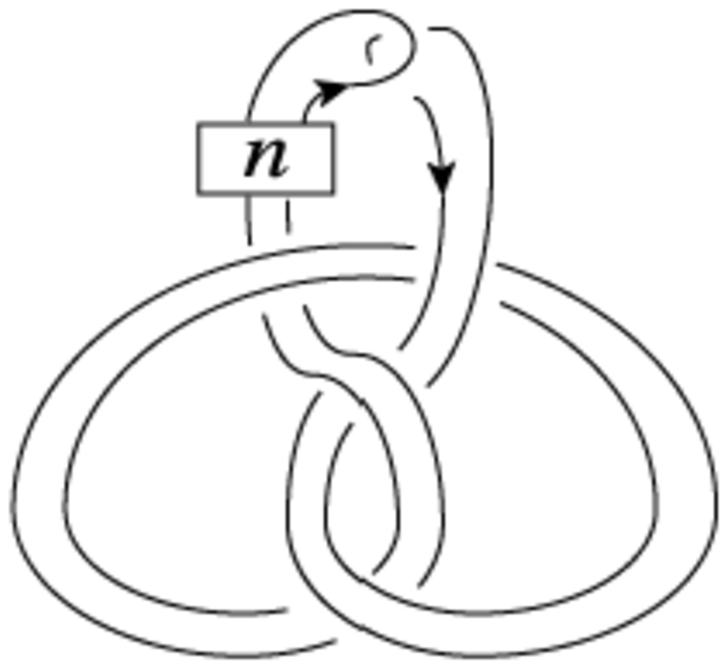}
  \end{center}
 \caption{$K_-$}
   \label{fig72}
 \end{minipage}%
 \begin{minipage}{0.33\hsize}
  \begin{center}
   \includegraphics[height=30mm]{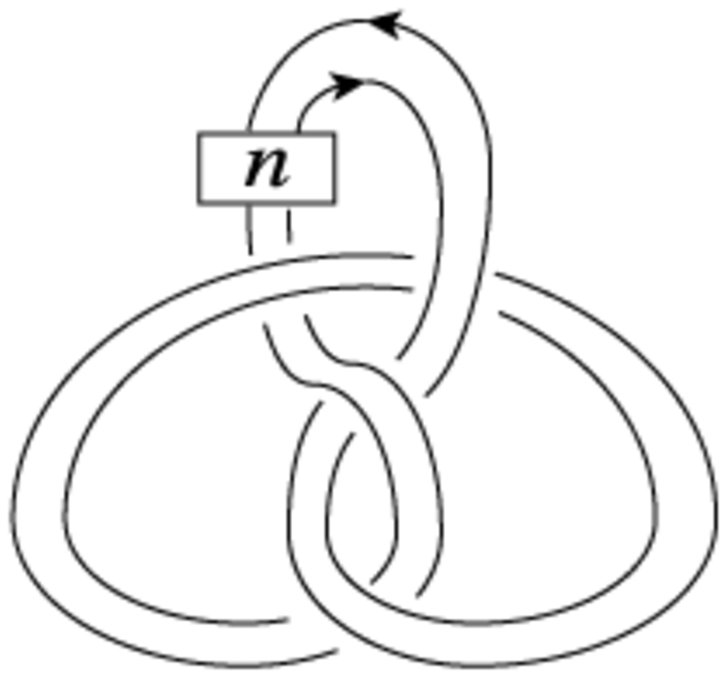}
  \end{center}
 \caption{$K_0$}
   \label{fig73}
 \end{minipage}
\end{figure}

By Lemma \ref{lem2}, $\lambda'(K_+) - \lambda'(K_-) = lk(K_0)$. Since $K_-$ is a trivial knot, $\lambda'(K_-) = 0$. $lk(K_0)$ is $-n$. Therefore, $\lambda'(K_+) = -n$. Then $\lambda(M_n(K_1, K_2)) = \lambda(S^3_{+1}(D_+(K_2, n))) = \lambda'(K_+) = -n$.
\endproof

\subsection{The Casson invariant of the fifth row on \thmref{thm:main}'s table}
\proof
$K_1$ and $K_2$ are figure eight knots. By 2(v), $M_n(K_1, K_2)$ is diffeomorphic to $S^3_{+1}(D_-(K_2, n))$. Therefore $\lambda(M_n(K_1, K_2))$ is equal to $\lambda(S^3_{+1}(D_-(K_2, n)))$. By Lemmas \ref{lem2} and \ref{lem3}, we can compute the Casson invariant $\lambda(S^3_{+1}(D_-(K_2, n)))$ as follows:

\begin{figure}[H]
 \begin{minipage}{0.33\hsize}
  \begin{center}
   \includegraphics[height=30mm]{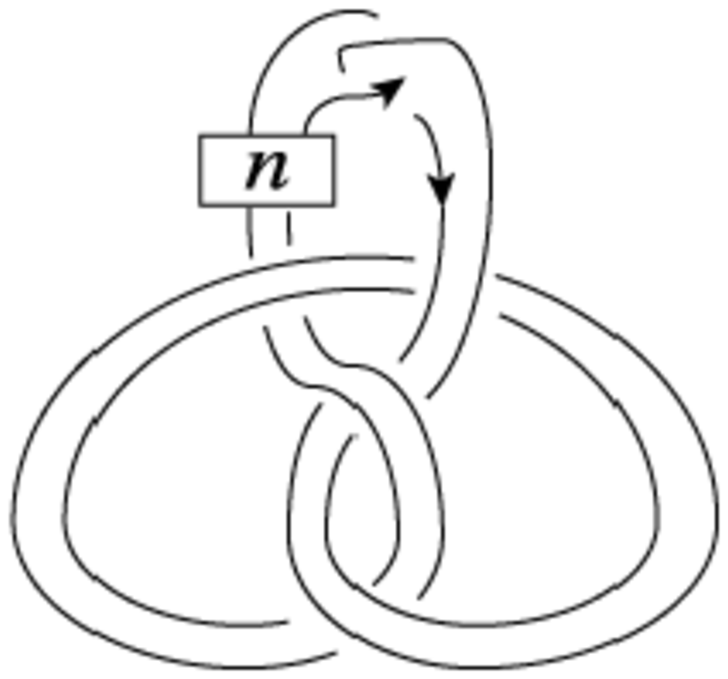}
  \end{center}
  \caption{$K_-$}
   \label{fig74}
 \end{minipage}%
 \begin{minipage}{0.33\hsize}
  \begin{center}
   \includegraphics[height=30mm]{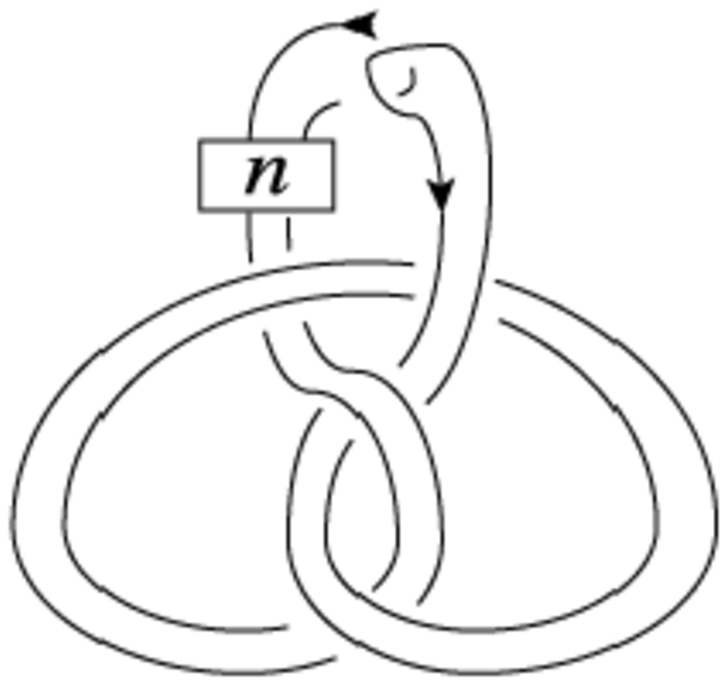}
  \end{center}
 \caption{$K_+$}
   \label{fig75}
 \end{minipage}%
 \begin{minipage}{0.33\hsize}
  \begin{center}
   \includegraphics[height=30mm]{Figure73.eps}
  \end{center}
 \caption{$K_0$}
   \label{fig76}
 \end{minipage}
\end{figure}

By Lemma \ref{lem2}, $\lambda'(K_+) - \lambda'(K_-) = lk(K_0)$. Since $K_+$ is a trivial knot, $\lambda'(K_+) = 0$. $lk(K_0)$ is $-n$. Therefore, $\lambda'(K_-) = n$. Then $\lambda(M_n(K_1, K_2)) = \lambda(S^3_{+1}(D_-(K_2, n))) = \lambda'(K_-) = n$.
\endproof
%
%
%
\section{Proof of Proposition \ref{prop2}}\label{sec:3}
%
%
We show that $V^1_n \cup_{\partial} (-V^2_n)$ is diffeomorphic to $\CP^2 \sharp \CP^2$.

\proof
By Kirby Calculus, we will show that the Kirby diagram of $V^1_n \cup_{\partial} (-V^2_n)$ is represented by Figure \ref{fig93} :

\begin{figure}[H]
 \begin{minipage}{0.33\hsize}
  \begin{center}
   \includegraphics[height=30mm]{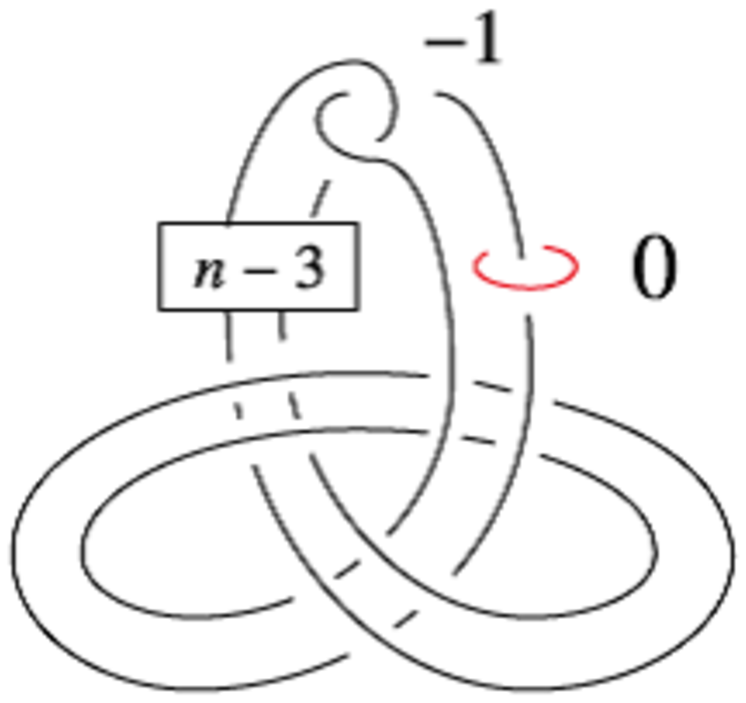}
  \end{center}
  \caption{}
   \label{fig88}
 \end{minipage}%
 \begin{minipage}{0.33\hsize}
  \begin{center}
   \includegraphics[height=30mm]{sproc.eps}
  \end{center}
 \end{minipage}%
 \begin{minipage}{0.34\hsize}
  \begin{center}
   \includegraphics[height=30mm]{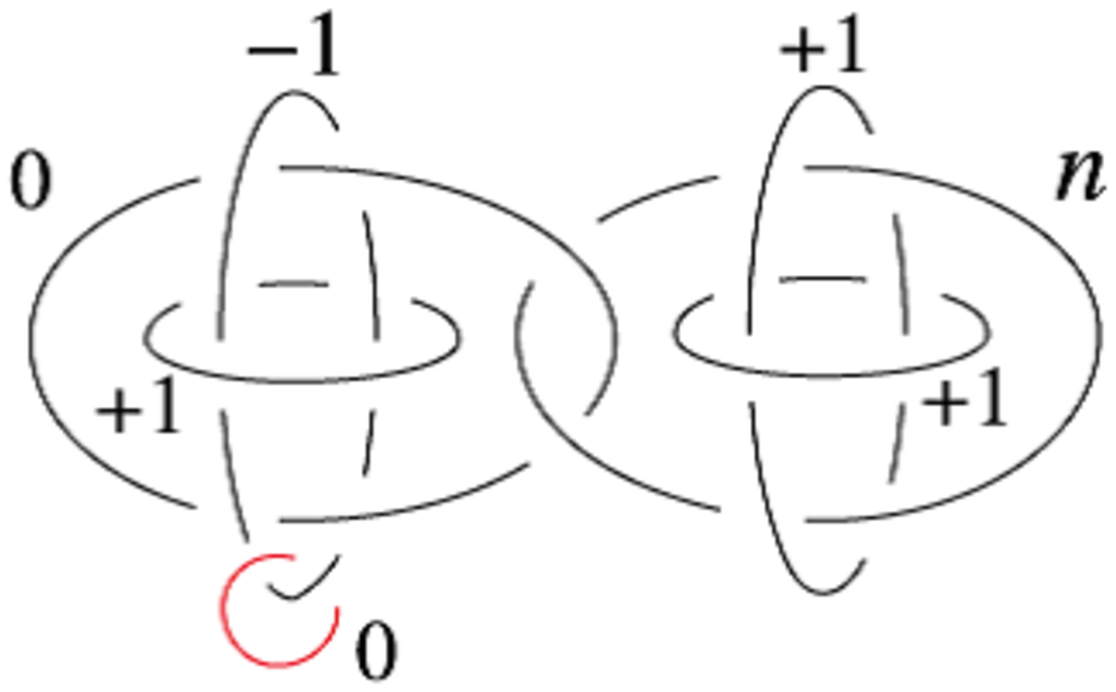}
  \end{center}
 \caption{}
   \label{fig89}
 \end{minipage}
\end{figure}

\begin{figure}[H]
 \begin{minipage}{0.1\hsize}
  \begin{center}
   \includegraphics[height=10mm]{isotopy.eps}
  \end{center}
 \end{minipage}%
 \begin{minipage}{0.3\hsize}
  \begin{center}
   \includegraphics[height=30mm]{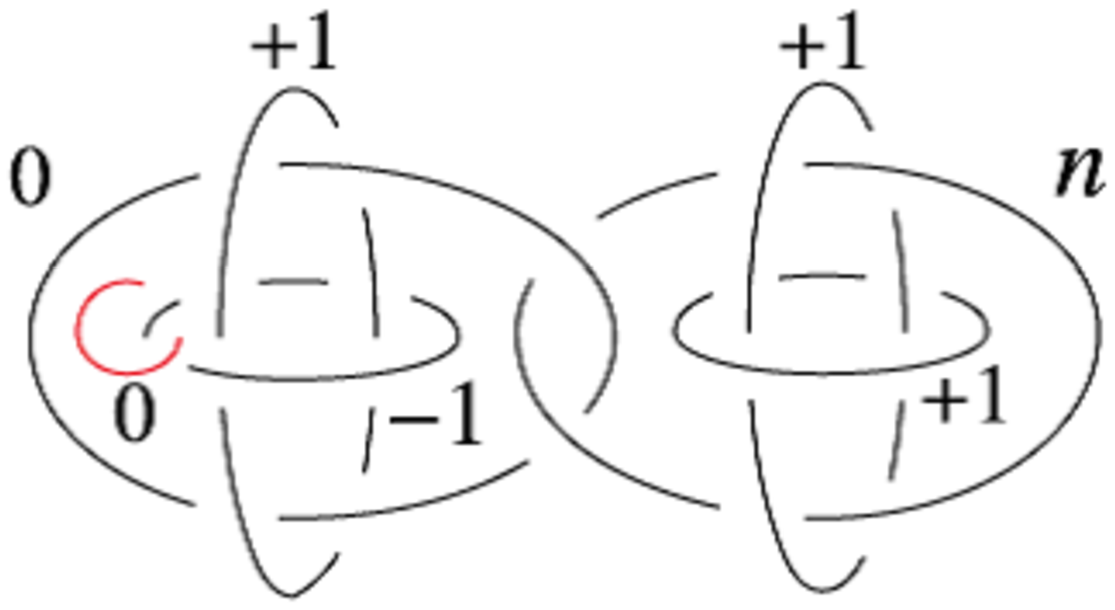}
  \end{center}
  \caption{}
   \label{fig90}
 \end{minipage}%
 \begin{minipage}{0.3\hsize}
  \begin{center}
   \includegraphics[height=30mm]{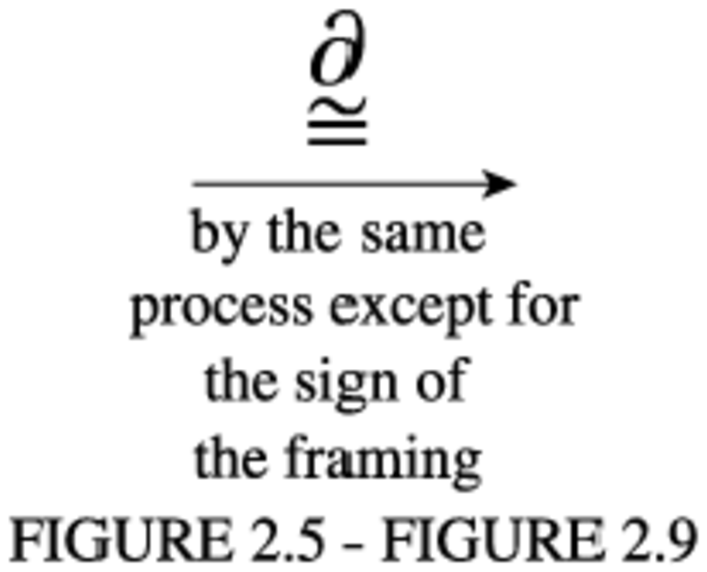}
  \end{center}
 \end{minipage}%
 \begin{minipage}{0.3\hsize}
  \begin{center}
   \includegraphics[height=30mm]{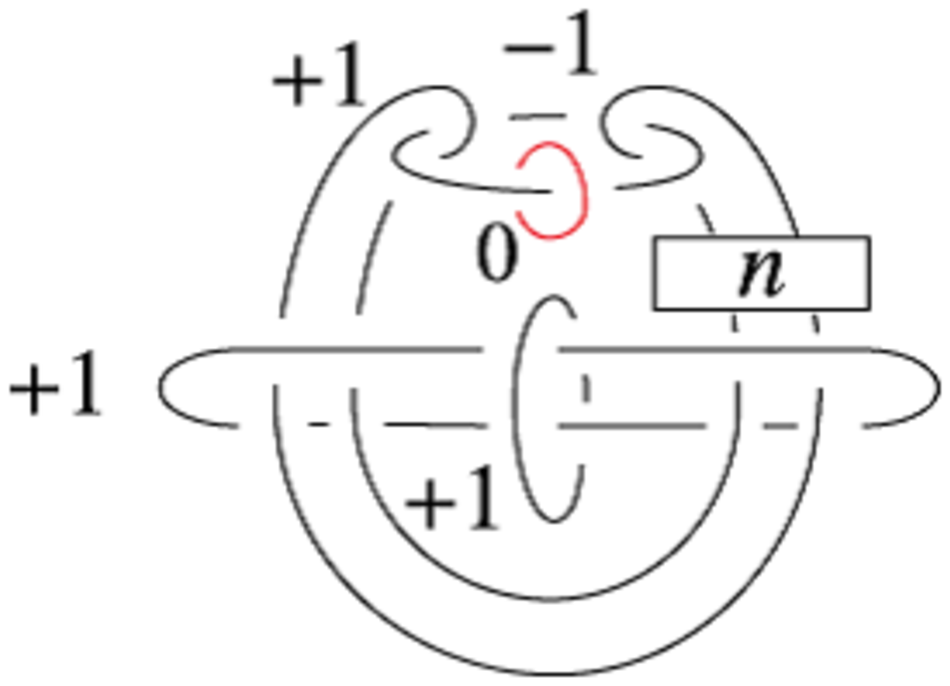}
  \end{center}
  \caption{}
   \label{fig91}
 \end{minipage}
\end{figure}

\begin{figure}[H]
 \begin{minipage}{0.1\hsize}
  \begin{center}
   \includegraphics[height=10mm]{blowdown.eps}
  \end{center}
 \end{minipage}%
 \begin{minipage}{0.4\hsize}
  \begin{center}
   \includegraphics[height=30mm]{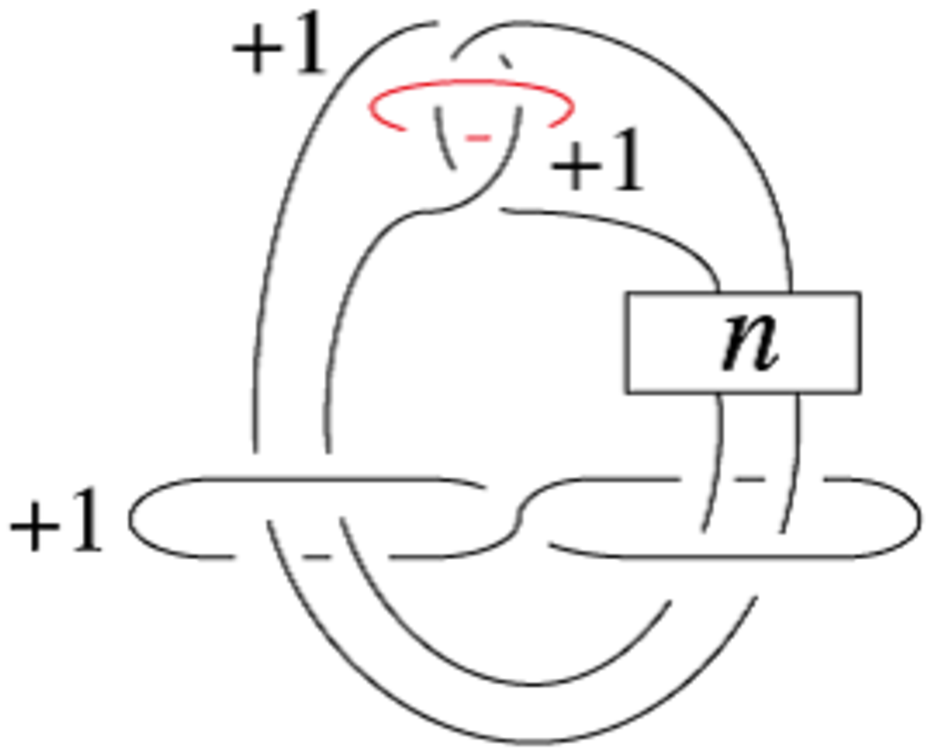}
  \end{center}
  \caption{}
   \label{fig92}
 \end{minipage}%
 \begin{minipage}{0.1\hsize}
  \begin{center}
   \includegraphics[height=10mm]{blowdown.eps}
  \end{center}
 \end{minipage}%
 \begin{minipage}{0.4\hsize}
  \begin{center}
   \includegraphics[height=30mm]{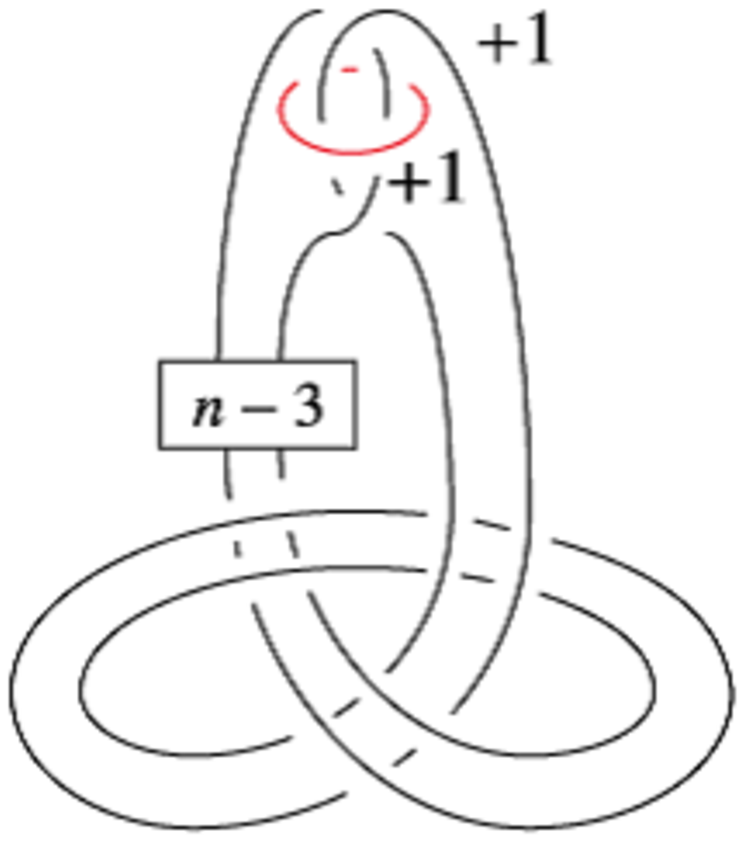}
  \end{center}
  \caption{$V^1_n \cup_{\partial} (-V^2_n)$}
   \label{fig93}
 \end{minipage}
\end{figure}

\begin{figure}[H]
 \begin{minipage}{0.5\hsize}
  \begin{center}
   \includegraphics[height=15mm]{hds.eps}
  \end{center}
 \end{minipage}%
 \begin{minipage}{0.5\hsize}
  \begin{center}
   \includegraphics[height=30mm]{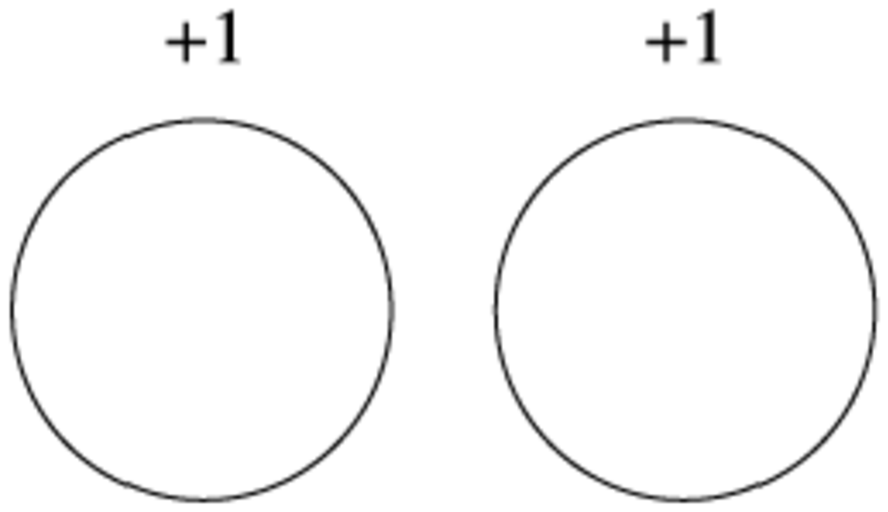}
  \end{center}
  \caption{$\CP^2 \sharp \CP^2$}
  \label{fig94}
 \end{minipage}
\end{figure}
\endproof
%
%
\section{Appendix}\label{sec:Enriques}
%
%
\noindent\textbf{An alternative proof of Corollary \ref{corG}.} By \cite{G}, if $n$ is odd, $M_n(T_{2,3}, T_{2,3})$ does not bound any contractible $4$-manifold. In this Section we will give an alternative proof of this fact. For this purpose, we will prove the following proposition;

\begin{figure}[H]
 \begin{minipage}{0.45\hsize}
  \begin{center}
   \includegraphics[height=20mm]{Figure5.eps}
  \end{center}
  \caption{$M_n(T_{2,3}, T_{2,3})$}
   \label{fig77}
 \end{minipage}%
 \begin{minipage}{0.1\hsize}
  \begin{center}
   \includegraphics[height=10mm]{isotopy.eps}
  \end{center}
 \end{minipage}%
 \begin{minipage}{0.45\hsize}
  \begin{center}
   \includegraphics[height=20mm]{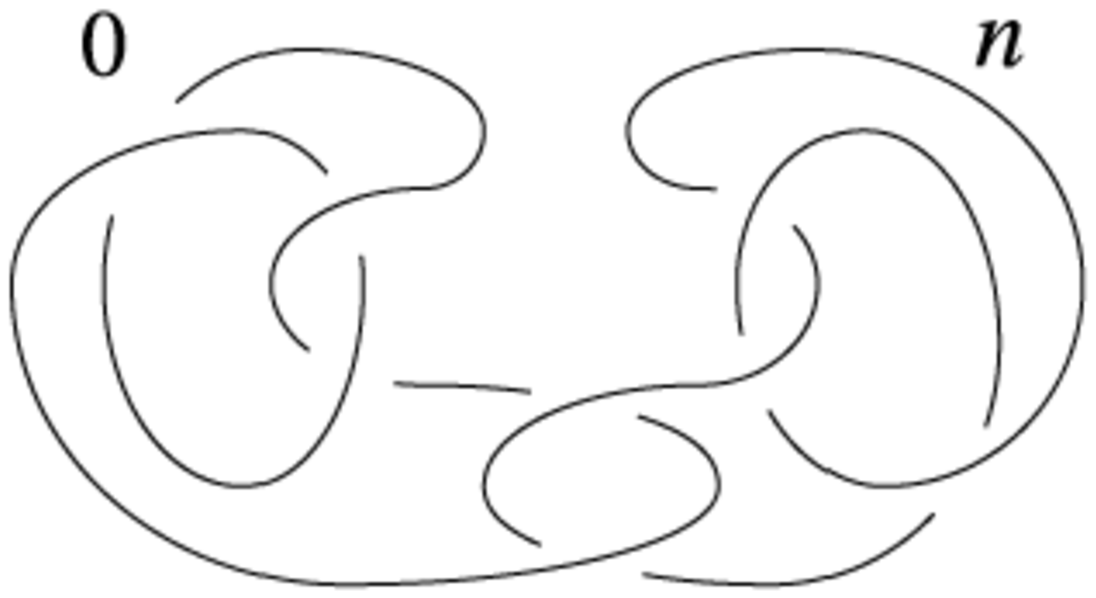}
  \end{center}
 \caption{}
   \label{fig78}
 \end{minipage}
\end{figure}

\begin{Proposition}\label{prop3}
The $4$-dimensional handlebodies represented by Figures \ref{fig77}, \ref{fig84} and \ref{fig85} have the same boundaries.
\end{Proposition}

\proof

\begin{figure}[H]
 \begin{minipage}{0.1\hsize}
  \begin{center}
   \includegraphics[height=10mm]{blowup.eps}
  \end{center}
 \end{minipage}%
 \begin{minipage}{0.4\hsize}
  \begin{center}
   \includegraphics[height=30mm]{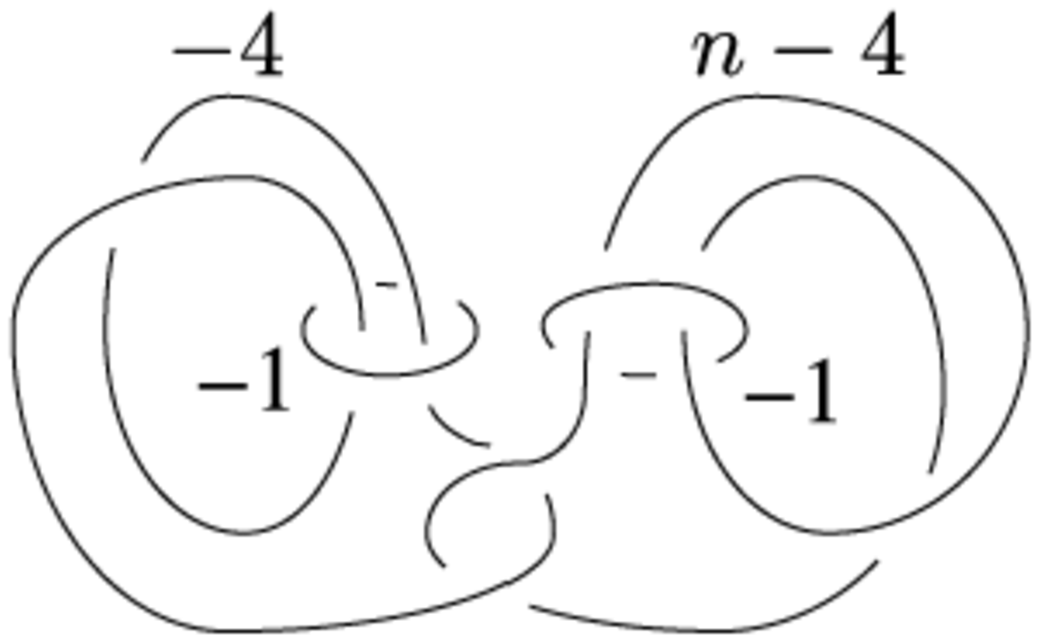}
  \end{center}
  \caption{}
   \label{fig79}
 \end{minipage}%
 \begin{minipage}{0.1\hsize}
  \begin{center}
   \includegraphics[height=10mm]{isotopy.eps}
  \end{center}
 \end{minipage}%
 \begin{minipage}{0.4\hsize}
  \begin{center}
   \includegraphics[height=30mm]{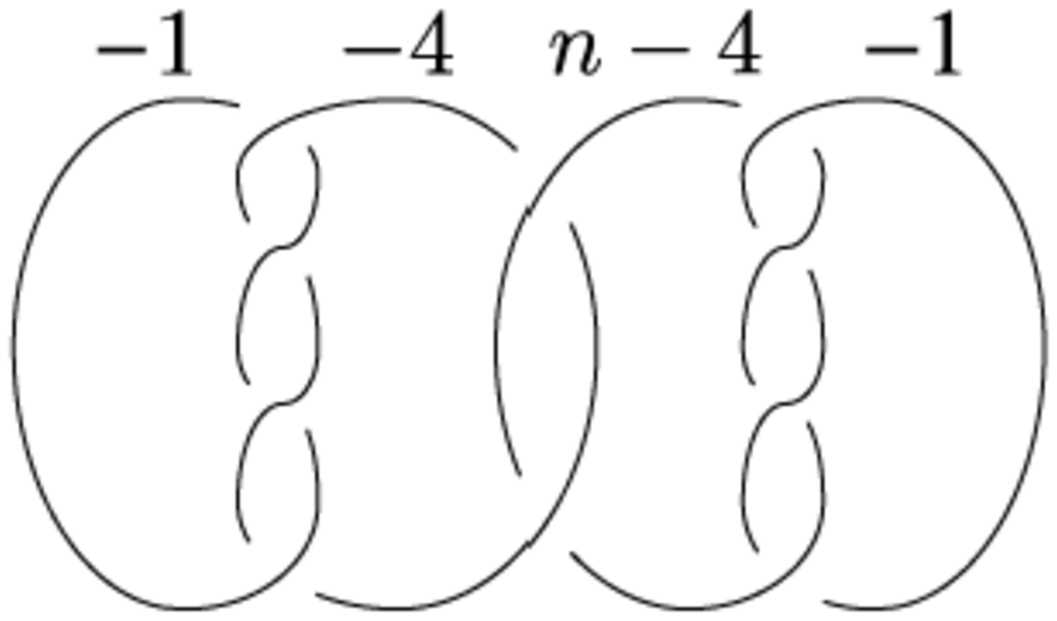}
  \end{center}
  \caption{}
   \label{fig80}
 \end{minipage}
\end{figure}

\begin{figure}[H]
 \begin{minipage}{0.1\hsize}
  \begin{center}
   \includegraphics[height=10mm]{blowup.eps}
  \end{center}
 \end{minipage}%
 \begin{minipage}{0.4\hsize}
  \begin{center}
   \includegraphics[height=30mm]{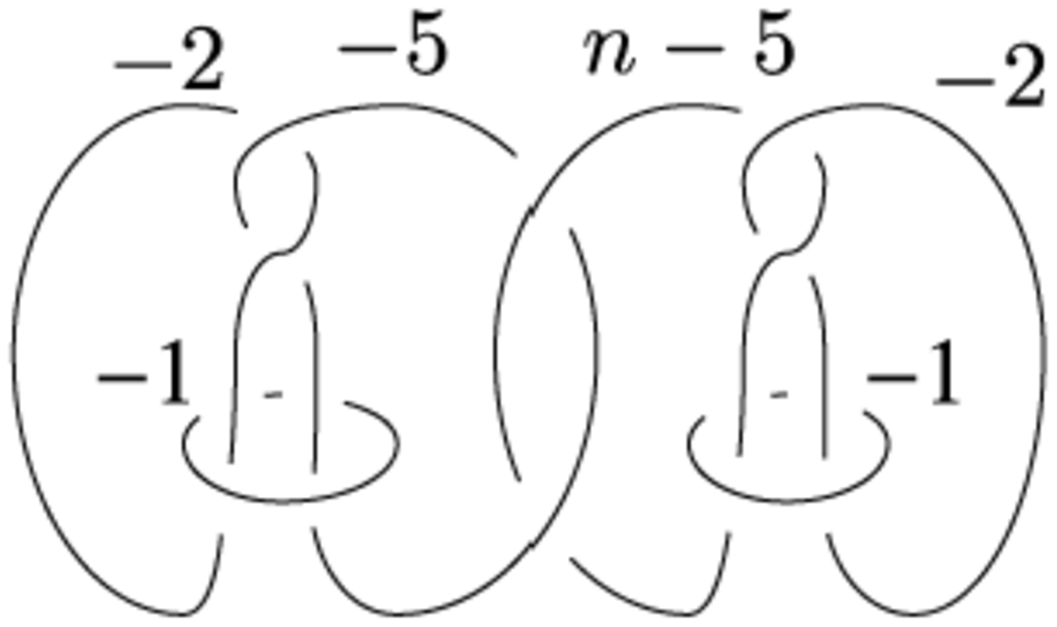}
  \end{center}
  \caption{}
   \label{fig81}
 \end{minipage}%
 \begin{minipage}{0.1\hsize}
  \begin{center}
   \includegraphics[height=10mm]{blowup.eps}
  \end{center}
 \end{minipage}%
 \begin{minipage}{0.4\hsize}
  \begin{center}
   \includegraphics[height=30mm]{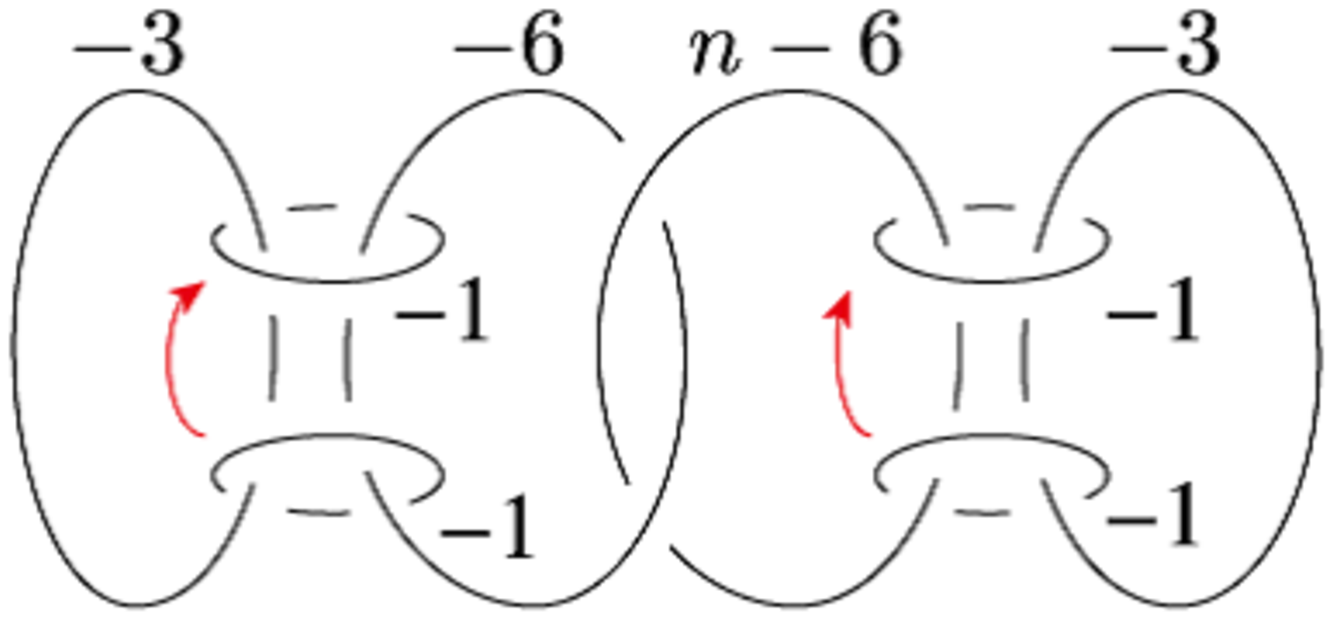}
  \end{center}
  \caption{}
   \label{fig82}
 \end{minipage}
\end{figure}

\begin{figure}[H]
 \begin{minipage}{0.1\hsize}
  \begin{center}
   \includegraphics[height=15mm]{hds.eps}
  \end{center}
 \end{minipage}%
 \begin{minipage}{0.9\hsize}
  \begin{center}
   \includegraphics[height=30mm]{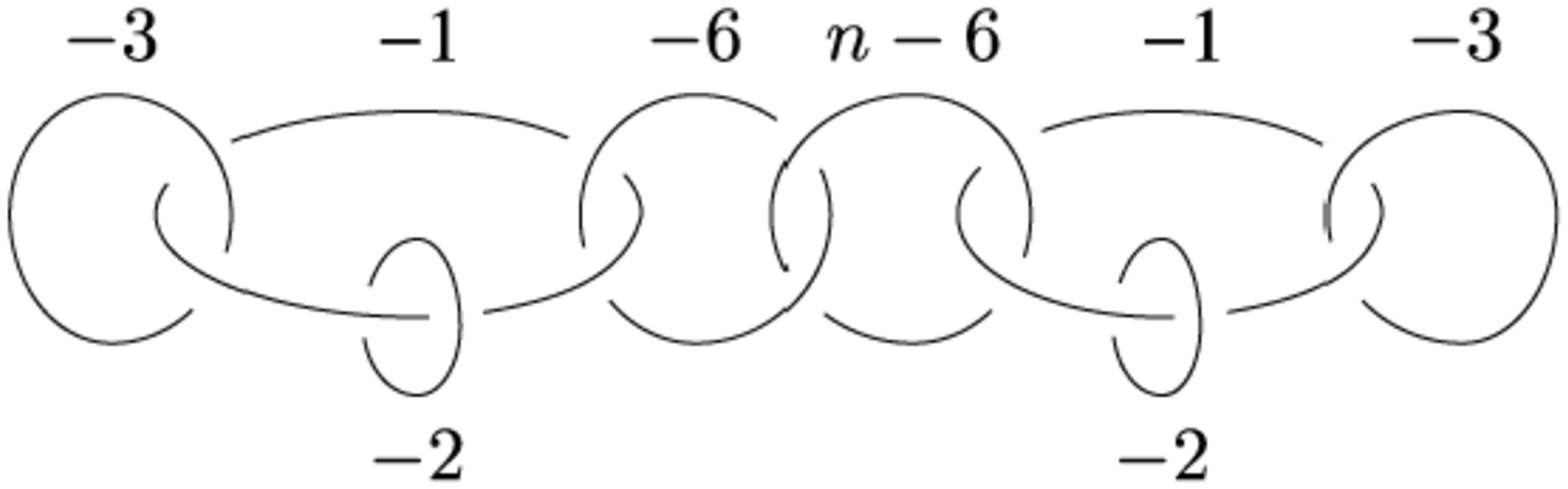}
  \end{center}
  \caption{}
  \label{fig83}
 \end{minipage}
\end{figure}

\begin{figure}[H]
 \begin{minipage}{0.1\hsize}
  \begin{center}
   \includegraphics[height=20mm]{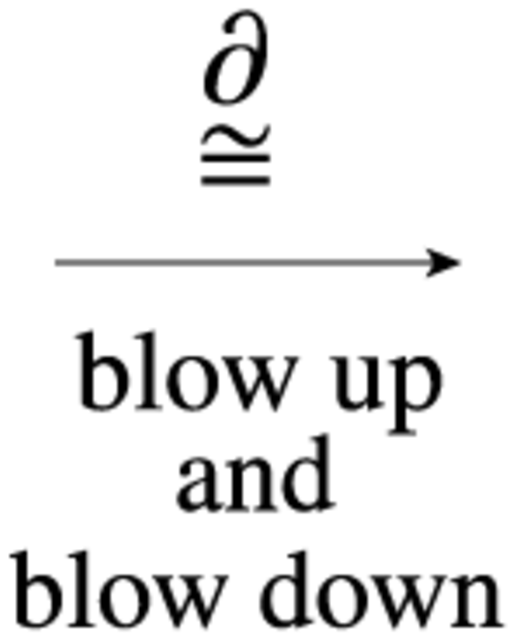}
  \end{center}
 \end{minipage}%
 \begin{minipage}{0.9\hsize}
  \begin{center}
   \includegraphics[height=30mm]{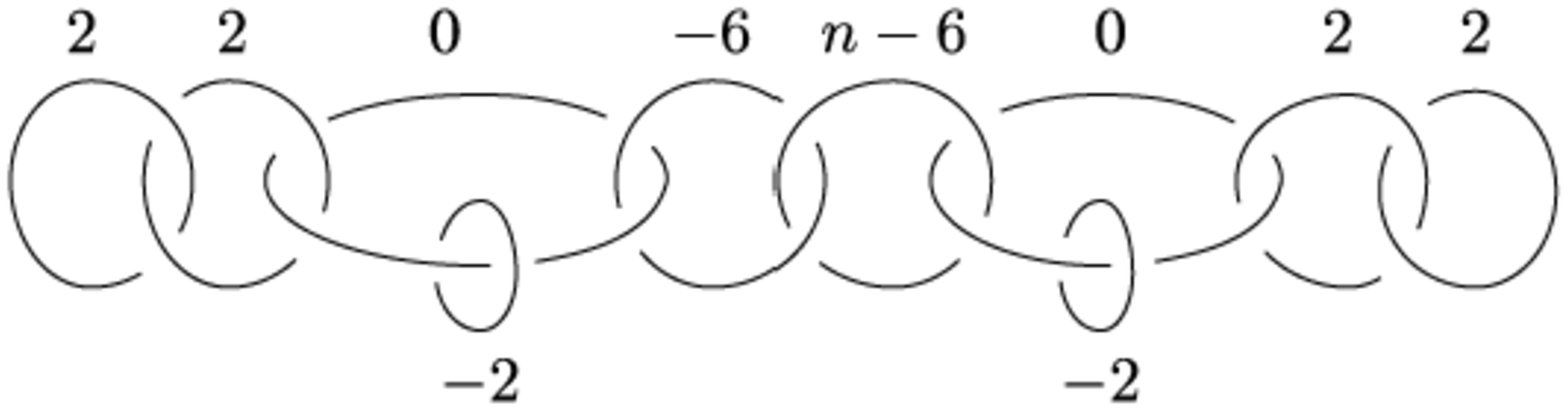}
  \end{center}
  \caption{}
  \label{fig84}
 \end{minipage}
\end{figure}

\begin{figure}[H]
 \begin{minipage}{0.1\hsize}
  \begin{center}
   \includegraphics[height=20mm]{blowupdown.eps}
  \end{center}
 \end{minipage}%
 \begin{minipage}{0.9\hsize}
  \begin{center}
   \includegraphics[height=30mm]{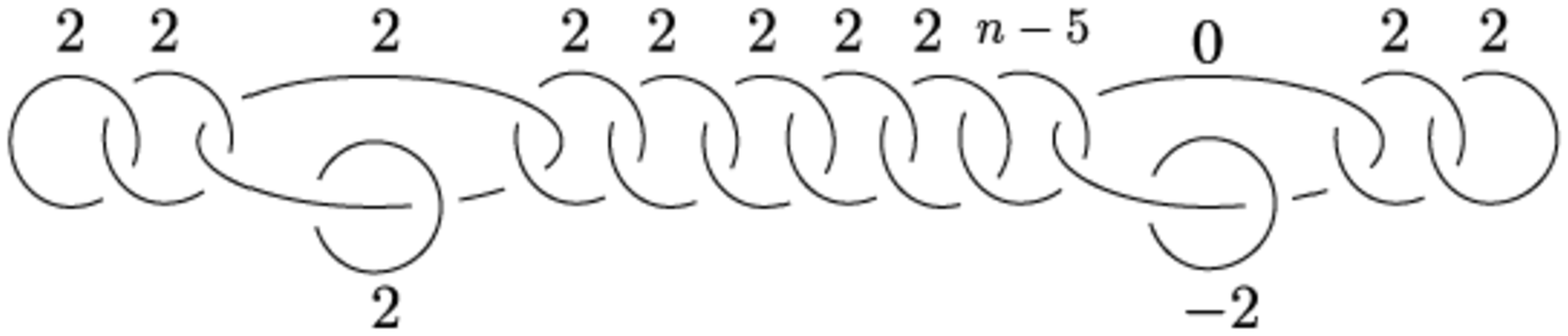}
  \end{center}
  \caption{}
  \label{fig85}
 \end{minipage}
\end{figure}
\endproof

We give an alternative proof of Corollary \ref{corG}.
\proof
Figure \ref{fig84} gives a smooth $4$-manifold $Q_1$ with intersection form $A$.

\begin{figure}[H]
 \begin{center}
  \includegraphics[height=30mm]{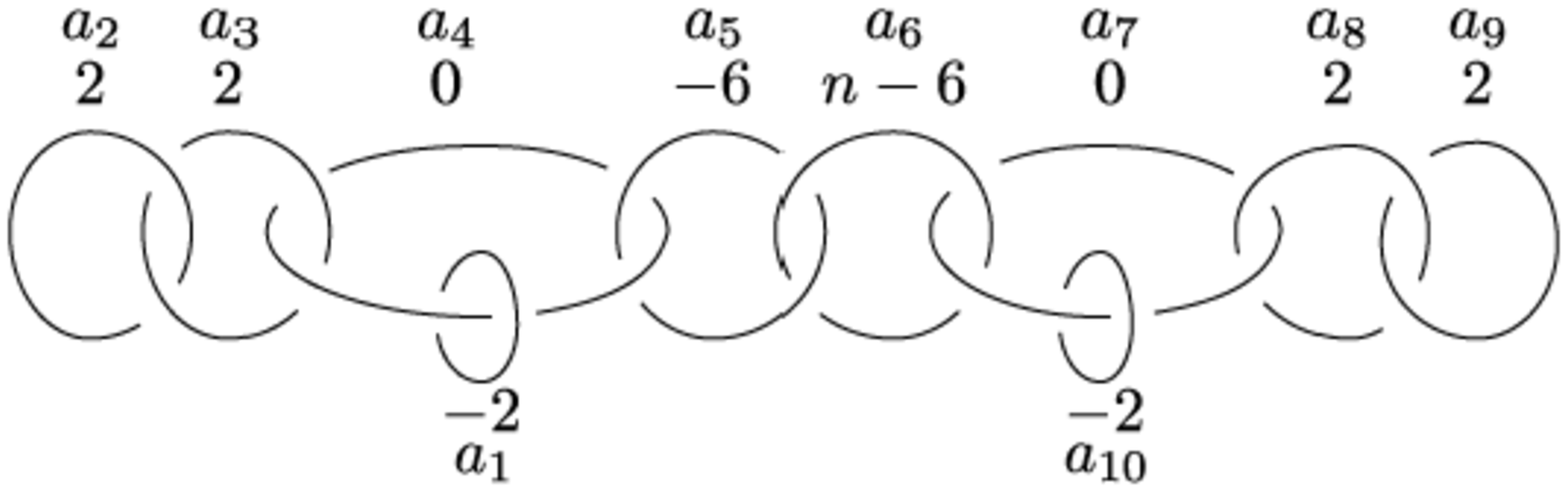}
 \end{center}
 \caption{$Q_1$}
 \label{fig86}
\end{figure}

\begin{center}
$A = (\alpha_{ij}), \; \alpha_{ij} = a_i \cdot a_j$, \; $1 \leq i, j \leq 10$
\end{center}
\begin{center}
$A$ $=$ $\left( \begin{array}{cccccccccc}
$-2$ & & & 1 & & & & & &\\
 & 2 & 1 & & & & & & &\\
 & 1 & 2 & 1 & & & & & &\\
1 &  & 1 & 0 & 1 & & & & &\\
 & & & 1 & $-6$ &1 & & & &\\
 & & & & 1 & $n-6$ & 1 & & &\\
 & & & & & 1 & 0 & 1 & & 1\\
 & & & & & & 1 & 2 & 1 &\\
 & & & & & & & 1 & 2 &\\
 & & & & & & 1 & & & $-2$
\end{array} \right)$
\end{center}

\noindent We have Index($A$) $= 0$. Note that $A$ is an even type matrix if $n$ is even.

Figure \ref{fig85} gives a smooth $4$-manifold $Q_2$ with intersection form $B$.

\begin{figure}[H]
 \begin{center}
  \includegraphics[height=30mm]{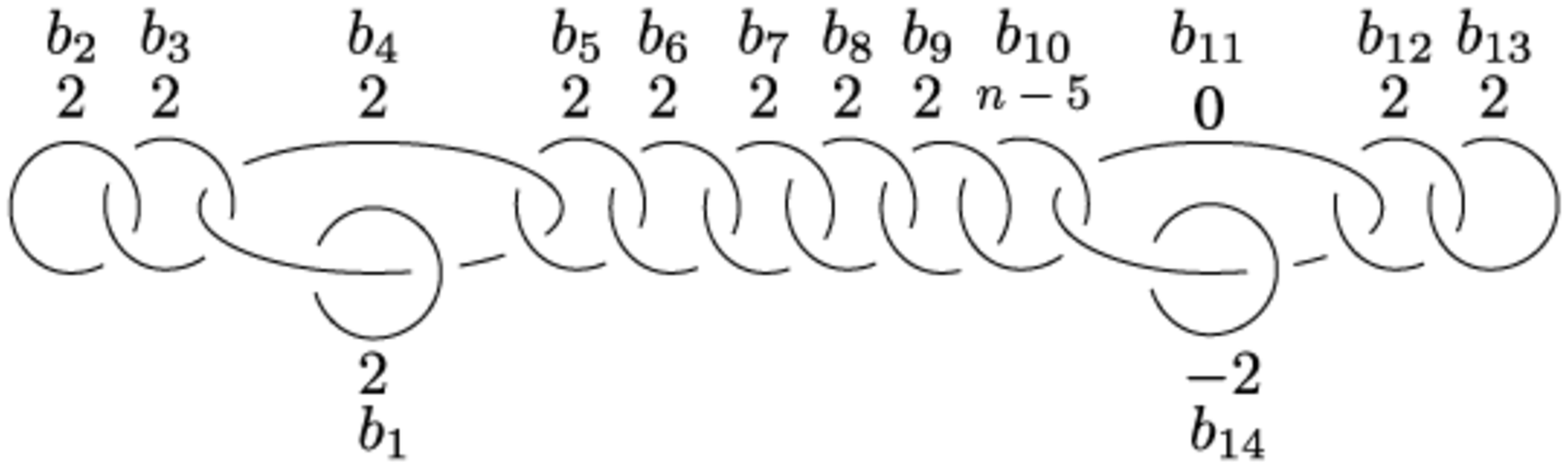}
 \end{center}
 \caption{$Q_2$}
 \label{fig87}
\end{figure}

\begin{center}
$B = (\beta_{ij}), \; \beta_{ij} = b_i \cdot b_j$, \; $1 \leq i, j \leq 14$
\end{center}
\begin{center}
$B$ $=$ $\left( \begin{array}{cccccccccccccc}
2 & & & 1 & & & & & & & & & &\\
 & 2 & 1 & & & & & & & & & & &\\
 & 1 & 2 & 1 & & & & & & & & & &\\
1 &  & 1 & 2 & 1 & & & & & & & & &\\
 & & & 1 & 2 &1 & & & & & & & &\\
 & & & & 1 & 2 & 1 & & & & & & &\\
 & & & & & 1 & 2& 1 & & & & & &\\
 & & & & & & 1 & 2 & 1 & & & & &\\
 & & & & & & & 1 & 2 & 1 & & & &\\
 & & & & & & & & 1 & $n-5$ & 1 & & &\\
 & & & & & & & & & 1 & 0 & 1 & & 1\\
 & & & & & & & & & & 1 & 2 & 1 & \\
 & & & & & & & & & & & 1 & 2 & \\
 & & & & & & & & & & 1 & & & $-2$
\end{array} \right)$
\end{center}

\noindent We have Index($B$) $= 8$. Note that $B$ is an even type matrix if $n$ is odd.

By Proposition \ref{prop3}, we have the Rohlin invariant $\mu(M_n(T_{2,3}, T_{2,3}))$ as follows:

\begin{equation*}
\mu(M_n(T_{2,3}, T_{2,3})) \equiv
\begin{cases}
\textrm{Index}(B) \equiv 1\quad (n \textrm{ is odd})\\
\textrm{Index}(A) \equiv 0\quad (n \textrm{ is even})
\end{cases}
\mod2
\end{equation*}

Therefore if $n$ is odd, $M_n(T_{2,3}, T_{2,3})$ does not bound any contractible $4$-manifold.
\endproof

\end{document}